\documentclass[3p,times]{elsarticle}
\usepackage{float}
\usepackage[cmex10]{amsmath}
\usepackage[title]{appendix}
\usepackage{listings}
\usepackage{array}    % in preamble
\usepackage{multirow} % in preamble
\usepackage{enumitem}
\usepackage{soul}
\usepackage{algorithm,algorithmic}
\usepackage[font={footnotesize,it}]{caption}

\usepackage{array}
\usepackage[utf8x]{inputenc}
\usepackage{csquotes}
\usepackage{natbib}
\usepackage{eurosym}
\usepackage{appendix}
\usepackage{bm}
\usepackage{upgreek}

\usepackage{amsthm}

\usepackage{mathrsfs}

\usepackage{setspace}
\usepackage{times}

\usepackage{multicol}
\usepackage{booktabs,caption,threeparttable}
\usepackage{exscale,relsize}
\usepackage{varwidth}

\usepackage{siunitx}
\usepackage{amssymb}
\usepackage{subcaption}
\usepackage{afterpage}
\usepackage{longtable}

\usepackage{amssymb,cases}

\usepackage{eqnarray}
\usepackage{commath}
\usepackage{tabularx}
\usepackage{graphicx}

\DeclareUnicodeCharacter{03B5}{\ensuremath{\epsilon}}

\usepackage{multirow}

\usepackage{amsthm}
\usepackage [long]{optidef}
\usepackage{epstopdf}
\usepackage{fancyhdr}
\usepackage{mleftright}

\newcommand{\bi}{\begin{itemize}}
	\newcommand{\ei}{\end{itemize}}

\newcommand*{\authorimg}[1]%
{ \raisebox{-1\baselineskip}{\includegraphics[width=\imagesize]{#1}}}
\newlength\imagesize    % new lwngth for determining image size and label width

\newcommand{\Iota}{\mathrm{I}}

\DeclareRobustCommand{\BilevelOpt}{\textnormal{$P_B$}}
\DeclareRobustCommand{\BilevelOptwithepsilon}{\textnormal{$P_B(\epsilon)$}}

\DeclareRobustCommand{\Coupled}{\textnormal{Alternating Bilevel Optimization-MILP}}
\DeclareRobustCommand{\Coupledacronym}{\textnormal{ABO-MILP}}

\DeclareRobustCommand{\Exhaustive}{\textnormal{Alternating Bilevel Optimization-Exhaustive Search}}
\DeclareRobustCommand{\ExhaustiveAcronym}{\textnormal{ABO-ES}}

\DeclareRobustCommand{\Surrogate}{\textnormal{Alternating Bilevel Optimization-Surrogate Solver}}
\DeclareRobustCommand{\SurrogateAcronym}{\textnormal{ABO-SS}}

\DeclareRobustCommand{\Approximate}{\textnormal{Approximate User Equilibrium}}
\DeclareRobustCommand{\ApproximateAcronym}{\textnormal{AUE}}

\renewcommand{\epsilon}{\varepsilon}

\biboptions{authoryear}

\usepackage[figuresright]{rotating}
\usepackage[bookmarks=false]{hyperref}
\hypersetup{colorlinks,
linkcolor=blue,
citecolor=blue,
urlcolor=blue}

\begin{document}
\begin{frontmatter}
	
	
	
	%\title{Coordinated School Scheduling for Morning Congestion Mitigation: A User Equilibrium Approach}
	
	%\title{Coordinated school scheduling with departure time choice for congestion mitigation in multi-region urban networks}
	
	%\title{\textcolor{blue}{On the integration of departure choice and school scheduling in a bi-level optimization framework for urban congestion mitigation}}
	
	\title{Mitigating Regional Traffic Congestion via School Start Time Scheduling: A Bilevel Alternating Optimization Approach}
	
	%\title{A bi-level optimization framework integrating departure time choice and school scheduling in urban networks}
	
	%\title{\textcolor{blue}{A bi-level mathematical programming framework for optimizing school start times in multi-region urban networks}}
	
	%\title{\textcolor{blue}{Mitigating urban peak congestion through bi-level optimization of school start times in multi-region urban networks.}}
	
	%\title{\textcolor{blue}{On the coupling between school start time optimization and departure time choice for congestion mitigation}}

	
	\author[]{A. Georgantas\corref{cor1}}
	\author[]{S. Timotheou}
	\author[]{C. G. Panayiotou}
	
	\address[a]{KIOS Research Center of Innovation and Excellence, Department of Electrical and Computer Engineering, University of Cyprus, 1 Aglantzia Avenue, Nicosia 2109, Cyprus}
	
	%\fntext[myfootnotetel]{The work of Antonios Georgantas, Charalambos Menelaou, Stelios Timotheou and Christos G. Panayiotou is supported in part by the European Union’s Horizon 2020 research and innovation programme under Grant 739551 (KIOS CoE), in part by the Government of the Republic of Cyprus through the Directorate General for European Programmes, Coordination, and Development.}

	\begin{abstract}
		
		This paper addresses the challenge of morning commute congestion caused by concentrated school-related trips in urban networks. While traditional approaches to school start-time coordination often focus on isolated corridors or system-optimal perspectives, they frequently neglect network-wide dynamics, congestion patterns, and the behavioral responses of commuters. We propose a novel bi-level optimization framework for regulating school start times within a multi-region urban network characterized by Macroscopic Fundamental Diagrams (MFDs). This framework explicitly couples system-level regulation with multi-class user-equilibrium-based departure time choices. The Upper-Level problem is formulated as a non-convex optimization program that aims to jointly minimize total time spent and deviations from current school schedules, while the Lower-Level problem models commuter behavior through a deterministic dynamic multi-class user equilibrium formulation, incorporating $\alpha$-$\beta$-$\gamma$ preferences to account for travel time, earliness, and lateness costs. To overcome the computational challenges posed by the bilevel and large-scale nature of the problem, the non-convex traffic dynamics, and endogenous demand responses, an iterative algorithm that alternates between the Upper- and Lower-Level problems is developed. For the Upper-Level problem, a relaxed upper-bound formulation is constructed that can be solved with standard mathematical programming solvers. For the Lower-Level problem, an iterative algorithm is developed that yields an approximate solution to the deterministic dynamic user equilibrium problem. Numerical results demonstrate that the proposed framework significantly reduces congestion, providing transportation authorities with a robust, practical tool to evaluate the trade-offs between school start-time flexibility and traffic efficiency. Sensitivity analyses are further conducted to evaluate the impact of MFD uncertainty and varying scheduling preferences on network performance.

	\end{abstract}
	
	\cortext[cor1]{A. Georgantas. e-mail: georgantas.antonios@ucy.ac.cy}
\end{frontmatter}

%\correspondingauthor[*]{Corresponding author. Tel.: +0-000-000-0000 ; fax: +0-000-000-0000.}
%\email{georgantas.antonios@ucy.ac.cy}

% \linenumbers


%\textcolor{blue}{\textbf{14/02/2024 Be careful which parts to put in the ITSC paper and the journal. Must not have same parts that I cannot later put to the journal}}\\

%\textcolor{blue}{\textbf{put document in Grammarly before sending to the Professors}}\\
%\textcolor{blue}{\textbf{10/11/2023 as future work we could use a nonlinear MINLP solver such as BARON or FICO Xpress 9.2 to compare with our linear approximation model in terms of optimality gap}}\\
%\textcolor{blue}{\textbf{28/09/2023 Put in the density plots $\rho_1^C$ and not $\rho_r^C$}}\\
%\textcolor{blue}{\textbf{07/02/2024 Check where I need to use \text{cite} and where \text{citep} in the introduction. Put some part of the abstract from the proposal}}\\
%\textcolor{blue}{\textbf{07/02/2024 put literature from the proposal}}\\

%\textcolor{blue}{\textbf{03/11/2023 I have old papers in the introduction. Update them}}\\
%\textcolor{blue}{\textbf{I only shift in time in this paper, not in both in time and space. Correct this}}\\

\section{Introduction}
\label{ch:intro}

%Traffic congestion-particularly during morning rush hour is increasingly driven by the proliferation of private vehicles, leading to reduced travel speeds, extended queues, and longer trip durations \citep{papageorgiou2003review}. The congestion shortcomings pose a severe threat to the economy and the environment~\citep{khan2013estimating,levy2010evaluation}. In principle, this phenomenon typically arises when several classes of commuters simultaneously utilize the road network infrastructure to meet different goals \citep{SU2020334}.\par 

%\hspace{4cm}\textcolor{red}{\textbf{School-work commute: Congestion statistics}}\\

Morning congestion in urban networks is strongly influenced by school-related trips. A substantial share of vehicle commuters have to perform an intermediate stop by a private car to drop off their children at school before going to their final destination~\citep{JIA2016173, He2022}. This intermediate school stop forces \textit{school-work} commuters to enter the network at times and on routes that coincide with other trips. The resulting overlap with non-school-related traffic contributes to increased congestion levels across the network \citep{zheng2013distribution}.
%In modern cities, the majority of commuters are given the possibility to arrive at their workplaces within a flexible time range \citep{200135}. Based on data from the national household travel surveys for Paris and San Francisco,~\cite{MUNCH2023103712} argue that flexible-working hours are associated with a significantly higher probability of late arrival at work and may not directly lead to the mitigation of morning commute congestion.
The problem is further exacerbated when all schools start simultaneously. This is verified from three studies conducted in Beijing, China. Firstly, \cite{LU2017280} highlighted that workdays during school holidays result in a traffic congestion index $20\%$ lower than that of non-school holidays. Secondly, \cite{SUN2021103606} observed an increase in the probability of road congestion proliferation by 4.5$\%$ when schools start simultaneously. Thirdly, \cite{kang2024association} found that school runs significantly exacerbate traffic congestion around schools, reducing the likelihood of free-flow conditions by 8.34$\%$ during school run times and that school-related traffic congestion is more severe in areas with multiple schools, bus stops, and areas in close proximity with businesses. Consequently, synchronized school start times tend to generate measurable congestion peaks, especially during narrow drop‑off windows.
At the same time, the requirement imposed to school-work commuters to arrive punctually at their children’s school, highlights the need to investigate how school schedules influence morning congestion. One practical direction for tackling peak-period traffic is to adjust school start times, thereby implicitly managing the temporal distribution of travel demand within the road network.

%\hspace{4cm}\textcolor{red}{\textbf{Staggered work hours and Staggered school schedules}}\\

A commonly discussed approach to deal with the school-related morning congestion, involves altering institutional start times so as to redistribute travel demand across the peak period. One such approach is \textit{staggered work hours}, in which commuters are allowed to arrive to their work premises within a specified time range \citep{HENDERSON1981349}. While such measures can reduce peak concentration in some contexts, their effectiveness is substantially limited when school constraints are present. Commuters with flexible work hours often continue to depart during the peak hour, due to household coordination needs, such as school drop-offs \citep{munch2019irresistible, e3f7588f91244cf8b3766759f722ba33}. Behavioral analyses further reveal that when employees are given flexible work hours, they may choose to start later in the morning, which can shift, rather than reduce traffic demand and potentially further worsen late morning congestion \citep{MUNCH2023103712}. Moreover, modifying work schedules on its own may have only limited impact on the congestion patterns \citep{yang2014evaluating}. Overall, these findings suggest that the interplay between school and work schedules should be captured to achieve meaningful improvements in traffic performance. In the same context, \textit{staggered school schedules} act as a mechanism, where school start times are modified to redistribute morning traffic demand \citep{bertsimas2019optimizing}. This strategy was examined in Link\"oping, aiming to assess the sustainability gains that can be accompanied when changing the school start time~\citep{ljungberg2009staggered}.
\cite{Bertsimas2022} propose an optimization framework that staggers the schools' start times, which relies on a school bus routing algorithm, allowing school districts to explore trade-offs between competing priorities and choose times that best fulfill community needs. \cite{He2022} develop an optimization model coupling work and school start times, showing that under high demand, shifting school start times can significantly improve social welfare, although under low demand it may be counterproductive, yielding a Pareto frontier of policy options. In a different study, the authors devise an integer linear programming formulation to pinpoint school start times and bus operation times so that the transportation cost is minimized \citep{banerjee2019incorporatingequityschoolbus, doi:10.1287/trsc.2022.1130}. One limitation of these works is that the dynamic evolution of congestion is not adequately captured, which could potentially deliver valuable insights for constructing efficient congestion mitigation policies. In addition, most of the studies mentioned above examine their staggered schedule-based methodologies, assuming that buses alone are used to take children to schools. However, over $64\%$ of children in certain areas are being driven to school by a private car, often within a short distance, necessitating a more in-depth investigation of the impact of school start time changes on the car-based trips of commuters \citep{lay2010influencing}.

%\hspace{3cm}\textcolor{red}{\textbf{Description of core idea and justification for MFD}}\\

This study aims to determine school start times in an urban road network of arbitrary architecture, assuming regional traffic dynamics such that the \textit{school-related congestion} is mitigated. To model the mobility pattern involving the transition of commuters towards different destinations, we consider two distinct classes of commuters. Specifically, commuters of class S must make an intermediate stop to drop off children at school before proceeding to work, whereas commuters belonging to class W travel directly to their workplace. To capture these interactions at an aggregated spatial scale, the Macroscopic Fundamental Diagram (MFD) framework serves as a natural fit. The MFD characterizes the relationship between density and flow within an urban region, enabling a tractable yet realistic representation of congestion evolution across interconnected regions \citep{daganzo2007urban}. Interestingly, most empirical MFDs used in the literature are estimated from traffic data that inherently include school-related travel demand and associated temporal and spatial variations. Despite the fact that school trips may give rise to localized demand patterns at finer spatial scales, extensive empirical evidence shows that, when traffic states are sufficiently aggregated, the resulting MFDs remain stable and unimodal, thereby supporting a regional-level representation of network dynamics \citep{ji2012spatial}. This has been confirmed in urban regions with a high concentration of schools, such as Lyon’s 6\textsuperscript{th} district with approximately 40 schools \citep{BATISTA2019192} and central Barcelona with nearly 300 schools \citep{geroliminis2008existence}. Exploiting the concept of the MFD,~\cite{10185822} proposed an optimization-based formulation to find the optimal school start times that lead to the minimization of the total time spent of vehicles for a specialized network architecture, where commuters pursue only a forward movement. A major assumption of this work is that the solution retrieved from the optimization can only belong to the free-flow regime. Nonetheless, this can only be achieved if the temporal redistribution of traffic resulting from the change of the start time of each school, yields free-flow conditions. To deal with this limitation, the authors extended the operation of their school start time selection framework to account for a network of arbitrary size, utilizing a derivative-free optimization framework and a successive convexification-based algorithm, respectively \citep{10919691,11059330}. However, a common limitation identified in these three works involving MFD-based traffic dynamics is the absence of an explicit departure time choice model within the context of the school start time selection problem.

A growing body of literature indicates that system-level interventions that reshape temporal demand, such as school start time regulation, cannot be reliably assessed without explicitly accounting for commuters’ departure time choice behavior. Even modest schedule adjustments may induce substantial shifts in departure times, thereby altering congestion patterns and potentially offsetting the intended benefits of the intervention. While early studies examined departure time choice under heterogeneous user preferences and route choice considerations \citep{liu2011morning, thorhauge2021heterogeneity}, such behavioral responses become particularly critical when school and work-related constraints interact at the network level. In this context, regional MFD-based traffic models provide a suitable framework for embedding departure time choice behavior, as they allow congestion dynamics and temporal demand redistribution to be captured at an aggregated spatial scale. Several studies have demonstrated the applicability of multi-region MFD-based formulations in representing heterogeneous departure time decisions under varying congestion conditions. For instance, \cite{doi:10.1287/trsc.2022.1147} incorporated stochastic user preferences into a dynamic multi-region MFD model using a game-theoretic formulation to capture realistic departure patterns under varying congestion levels. Similarly, \cite{BATISTA2025104980} developed a comprehensive modeling approach integrating MFD dynamics with behavioral insights to assess the impact of departure time shifts on system-wide efficiency via a dynamic user equilibrium optimization formulation. In the same context, \cite{DUNCAN2025105008} proposed a stochastic user equilibrium formulation under a dynamic multi-region MFD setting and validated it using a large-scale real-world case study, highlighting the effectiveness of the MFD in modeling temporal and spatial heterogeneity in traffic flows. Despite these advances, existing MFD-based departure time choice models have largely been developed independently of the school start time selection problem. As a result, the behavioral feedback between school schedule regulation, heterogeneous departure time adjustments, and the ensuing evolution of congestion remains insufficiently explored. This gap motivates the explicit integration of departure time choice behavior into MFD-based frameworks for school start time optimization.\par

%\hspace{4cm}\textcolor{red}{\textbf{Proposed solution approach: Bi-level mathematical programming framework}}\\

%\textcolor{blue}{Importantly, we only control the school start times; work start times are known beforehand, following either a fixed or a flexible scheduling scheme. Nonetheless, the departure times of both classes of commuters are subject to change in view of the prospective school start time shifts.}

This study addresses these gaps by integrating school start time regulation with commuters’ departure time choices within a multi-region MFD framework. To this end, we develop a bilevel optimization framework, consisting of two individual problems. First, the Upper-Level problem is responsible for the derivation of the school start times. We develop a macroscopic approach that incorporates the traffic dynamics associated with each class of commuters (S and W), such that the coupling of the dynamics concerning each class corresponds to the region-level dynamics obtained from the MFD. The resulting formulation leads to a bi-objective optimization program, aiming to jointly minimize i) the total time spent (TTS) of all vehicles inside the network and ii) the resulting overall school start time change (STC). Then at the Lower-Level problem, commuters' responses are modeled through a deterministic dynamic user equilibrium (DUE) formulation, in which departure-time decisions follow the $\alpha$–$\beta$–$\gamma$ preference model that captures the aspects of \textit{travel time}, \textit{earliness cost} and \textit{lateness cost}, respectively \citep{vickrey1969congestion}. The preference model is extended to account for the co-existence of different classes of commuters in a multi-regional MFD-based setting. The resulting bilevel optimization problem is computationally challenging to solve due to its large scale, the nonlinear and nonconvex traffic dynamics embedded in the Upper-Level formulation, and the need to compute equilibrium-based demand responses at the Lower-Level formulation. To address these challenges, we develop an iterative solution approach that alternates between the Upper and Lower-Level problems, enabling coordination between scheduling decisions and behavioral responses while preserving their mutual feedback. Within this framework, the Upper-Level problem is highly nonconvex and nonlinear. To deal with this, we develop a convexification approximation approach that offers relatively fast and close-to-optimal solutions. In this approach, each nonconvex and nonlinear constraint is relaxed to a convex approximate version. The convexification procedure yields a bi-objective mixed-integer linear program (MILP). To further deal with the bi-objective nature of the problem, we utilize a scalarization technique known as \textit{$\epsilon$-constraint method} that transforms the two objective functions included in the optimization problem into a single one by inserting the second objective as a constraint. In this regard, we obtain a tractable MILP that can be optimally solved using standard mathematical optimization solvers. At the Lower-Level problem, the deterministic dynamic user equilibrium is formulated as a variational inequality (VI) problem, consistent with instantaneous equilibrium assumptions in large-scale MFD-based networks \citep{friesz2011}. Solving the exact VI problem is computationally demanding due to the endogenous coupling between departure-time decisions and nonlinear congestion dynamics. To this end, we develop a tailored iterative procedure that approximates the user-equilibrium conditions while remaining computationally tractable.

%\hspace{4cm}\textcolor{red}{\textbf{Contributions of this work}}\\

To the best of our knowledge, this is the first study that systematically integrates multi-region MFD-based congestion modeling with user-equilibrium-based departure time choices to optimize school start times. Specifically, the contributions of this work are the following: 

\begin{itemize}
	\item Design of a bilevel optimization problem for coordinated school start time regulation in multi-region MFD-based networks. The Upper-Level problem is formulated as a mixed-integer nonlinear program to determine system-optimal school schedules under multi-class traffic dynamics. The Lower-Level problem designs a novel deterministic dynamic multi-class user equilibrium model to capture heterogeneous departure-time responses, accounting for regional traffic flow interactions and class-specific perceived travel costs under congestion.
	\item Introduction of an iterative solution methodology with respect to the bilevel problem that alternates between the Upper and Lower-Level problems to consistently capture the feedback between school scheduling decisions and equilibrium demand responses.
	\item Construction of a linear relaxation solution method for the Upper-Level problem that approximates nonlinear and nonconvex constraints governing the traffic dynamics with linear ones, thereby yielding a tractable optimization formulation with reduced computational complexity that leads to high-quality feasible solutions.
	\item Development of an iterative algorithm that yields approximate solutions to the deterministic dynamic multi-class user equilibrium problem.
\end{itemize}

%\hspace{6cm}\textcolor{red}{\textbf{Benchmarking}}\\

To assess the performance of the proposed optimization framework, two benchmarking approaches are also developed. First, an Exhaustive Search (ES) procedure is designed to determine the globally optimal solution by evaluating all possible combinations of school start times and their corresponding impact on the TTS. However, due to the exponential growth in computational complexity with increasing problem size, the ES approach becomes intractable for large-scale instances. To enable comparison under such conditions, a derivative-free optimization method is also implemented, providing computationally efficient and high-quality solutions that approximate the optimal benchmark. \\

%\hspace{5cm}\textcolor{red}{\textbf{Structure of the paper}}\\

The remainder of the paper is organized as follows. Section \ref{sec:problem_description} presents the general problem description. Section \ref{sec:statement} introduces the formal problem statement for school start time selection, integrating class-specific regional traffic dynamics with commuters’ departure-time choice behavior. Section \ref{sec:solution} describes the proposed solution approach with respect to the bilevel optimization framework, while Section \ref{ch:simulation} presents the simulation results, highlighting the performance of the developed solution approaches. %Section \ref{sec:discussion} provides a comprehensive discussion of the results, interpreting the impact of coordinated school start times on network congestion, examining the sensitivity to commuter preferences and MFD uncertainty, and highlighting implications for urban traffic management and policy interventions.
Section \ref{ch:conclusions} concludes this work and discusses future research directions. Appendix~\ref{sec:es_contribution} details the operation of the proposed Exhaustive Search algorithm, while Appendix~\ref{sec:derivative_free} presents the structure of the solution methodology involving a derivative-free optimization approach. Finally, Appendix \ref{sec:appendixc} provides a comprehensive table of all parameters used in the study.

\textit{Notation}: All boldface letters indicate vectors (lowercase), whereas calligraphic letters denote sets. $|\mathcal{B}|$ denotes the cardinality of set $\mathcal{B}$, while $|x|$ denotes the absolute value of variable $x$. Positive real numbers are denoted by $\mathbb{R}^+$, and positive integer numbers are described by $\mathbb{Z}^+$, while non-negative real numbers are denoted by $\mathbb{R}_{\ge 0}$. The quantity $\binom{n}{b}$ denotes the binomial coefficient on $n$ and $b$, defined as $\frac{n!}{b!(n-b)!}$. The operator $\text{col}(\cdot)$ represents the collection of its arguments stacked into a column vector.

\section{General Problem Description}
\label{sec:problem_description}

This section describes the considered planning problem and the decision-making
structure governing school scheduling and commuter travel behavior in a regional urban
road network. We consider two distinct classes of vehicle commuters that utilize the network during the morning period:
\begin{itemize}
	\item \textbf{Class W:} \textit{Single-destination} commuters who travel
	directly from their origin to their destination.
	\item \textbf{Class S:} \textit{Dual-destination} commuters who make an
	intermediate stop at a school before proceeding to their final destination.
\end{itemize}

When schools start simultaneously, the aggregate school-related demand (class S) tends to be temporally concentrated so as to arrive at schools as close as possible to the prescribed start time, reflecting prevailing scheduling preferences at the OD pair level. Following the school drop-off, this demand continues toward final destinations, thereby temporally aligning with work-related demand (class W). This synchronization of school- and work-oriented travel induces pronounced peak demand periods, during which congestion intensifies and overall network performance deteriorates.

From a system planning perspective, in this paper the regulation of school start times in a multi-region urban network is performed such that Specifications (R1)-(R3) below are satisfied

\noindent\textit{(R1)}: Minimize the Total Time Spent (TTS) of all the vehicles in the network.\\
\textit{(R2)}: Minimize the Start Time Change (STC) between the initial and the shifted start time of each school.\\
\textit{(R3)}: Reflect the demand shifting flexibility of commuters to respect user equilibrium conditions for specific school start times.\\

\begin{figure}[t]
	\centering
	\includegraphics[width=12cm]{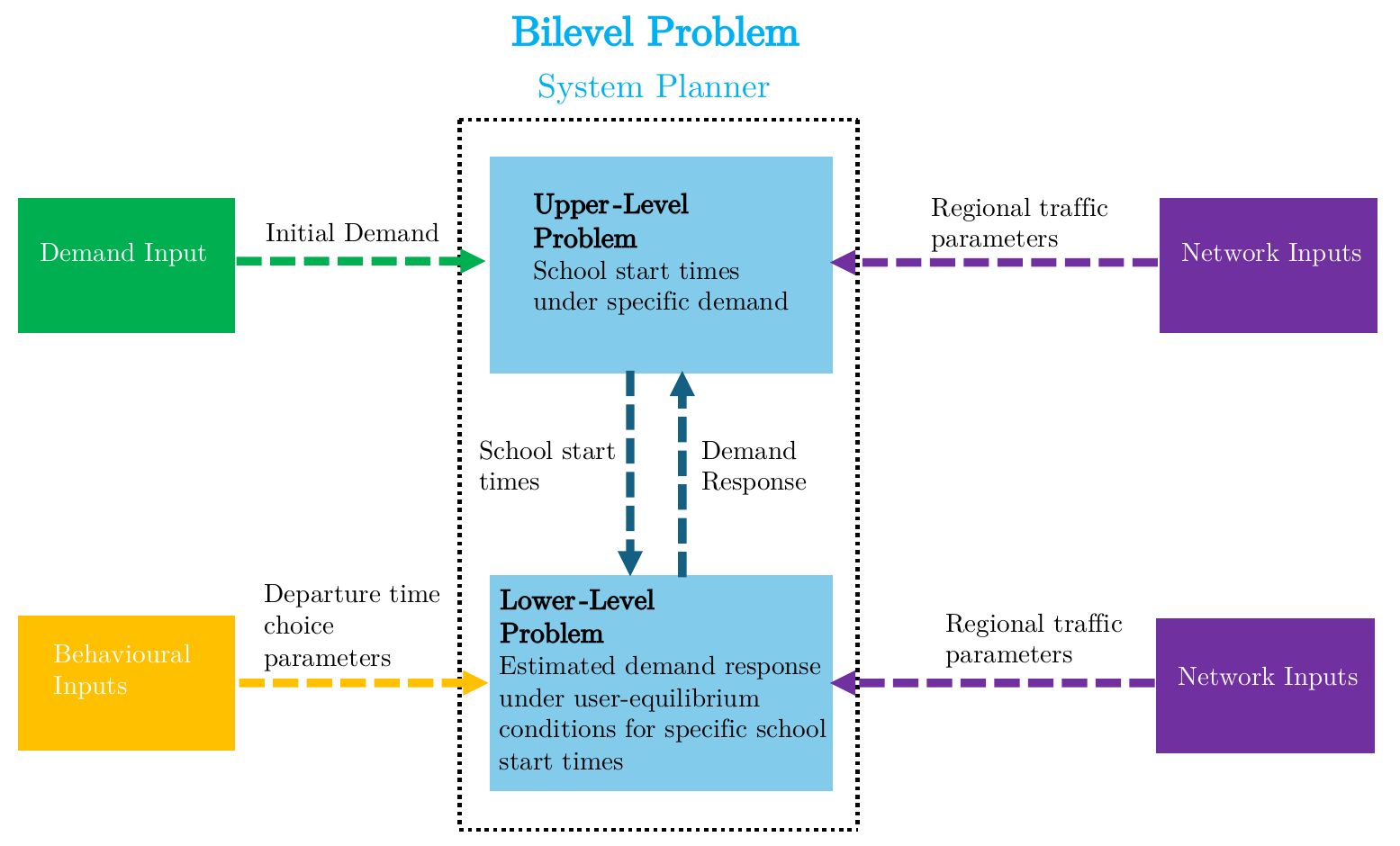} 
	\caption{Bilevel optimization problem for the design of the school start time selection structure.}
	\label{fig:overview_fig}
\end{figure}

Crucially, the effectiveness of any school start time configuration cannot be assessed independently of aggregate travel demand responses. Changes in school schedules modify congestion patterns and travel costs, which in turn affect the departure time distribution of OD-based demand. As a result, the temporal distribution of traffic demand emerges endogenously from congestion-induced, OD-level equilibrium adjustments. This hierarchical interaction between system-level scheduling decisions and aggregate demand responses is modeled as a bilevel optimization problem, as illustrated in Fig.~\ref{fig:overview_fig}. The entire problem (consisting of the Lower and Upper-Level) is solved by the system planner, who determines the school scheduling decisions while explicitly accounting for the resulting equilibrium adjustments in travel demand.

\section{Problem Statement}
\label{sec:statement}

The bi-level optimization problem aims to satisfy Specifications (R1)-(R3) through the regulation of school start times, as defined in Section \ref{sec:problem_description}. To this end, we first describe the Upper-Level problem that aims to provide the mathematical framework for selecting the school start times through the linkage of the regional traffic dynamics of the aforementioned classes of commuters with their associated demands. \par

\subsection{Upper-Level problem}
\label{sec:upper_level}

To capture the movement of vehicles in the network we consider macroscopic regional traffic dynamics. First, let $\mathcal{R}$ denote the set that contains all the regions inside the road network, $\mathcal{O}\subseteq \mathcal{R}$ the set of origins from where commuters start their trip, and $\mathcal{D}\subseteq \mathcal{R}$ the set of destinations, where commuters terminate their trip. Also, let $\mathcal{B}\subseteq \mathcal{R}$ denote the set of regions containing schools and $\mathcal{S}_b$ the set of schools located in region $b\in\mathcal{B}$.\par

The duration of each discrete school shifting interval is equal to $V$ (min), where $m\in\mathcal{M}$ denotes the number of discrete shifting intervals, with $\mathcal{M} = \{-M, -M+1, \ldots, 0, 1, \ldots, M\}$ representing the set of possible shifting intervals. Here, $M$ denotes the maximum number of backward or forward shifting intervals that can be applied to the start time of each school, resulting in $2M+1$ different intervals.\par

Let us also denote the start time of school $s\in\mathcal{S}_b$ with $\tau_s$. The selected start time for each school is represented through the binary variable $\xi_{m,s}, m\in\mathcal{M}, s\in\mathcal{S}_b, b\in\mathcal{B}$, defined as
\begin{equation*}
	\label{ggxda}
	\displaystyle   \xi_{m,s} = \left\{
	\begin{array}{lll}
		\displaystyle	\hspace{-0.2cm}1, \quad &\textrm{if}~ \text{school $s$ starts at time $\tau_s + mV$}, \\[6pt]
		\displaystyle	\hspace{-0.2cm}0, &\textrm{otherwise}.   
	\end{array} 
	\right. 
\end{equation*}

\noindent Binary variable $\xi_{m,s}$ is equal to 1 if the shifted start time of school $s$ becomes equal to $\tilde{\tau}_s = \tau_s + m V$, otherwise $\xi_{m,s}$ is equal to 0. %The aforementioned definitions link the shifted distribution of demand concerning commuters of class S, $\mathbf{d}_{o,s,m}^{S}$ with the candidate shifted school start time $\tau_s + m V$.  \par

The Upper-Level decision vector
\(
\boldsymbol{\xi} = \{\xi_{m,s}\}_{\forall m, \forall s}
\)
selects the start time shift applied to each school.
The feasible set is given by
\begin{equation}
	\Xi_U =
	\left\{
	\boldsymbol{\xi} :
	\sum_{m\in\mathcal{M}} \xi_{m,s} = 1,
	\;
	\xi_{m,s}\in\{0,1\},
	\;
	\forall s\in\mathcal{S}_b,\; \forall b\in\mathcal{B}
	\right\}.
\end{equation}

\subsubsection{Demand Modeling}
\label{sec:demand_modeling}

First, time is discretized into time-steps of duration $T$ (hour), such that the $k$th time-step expresses the time period $[kT, (k+1)T)$, for $k\in\mathcal{K}$ and $\mathcal{K}=\{0,1,\ldots,K\}$. For the sake of clarity, in this paper, we consider that commuters perform their trips using their private cars. Single-destination commuters are associated with a demand $d_{o,d}^W(k;\boldsymbol{\xi})$ (veh), which denotes the aggregate number of vehicles that enter the network at time-step $k$ from origin $o\in \mathcal{O}$ towards destination $d\in\mathcal{D}$\footnote{Demand $d_{o,d}^W(k;\boldsymbol{\xi})$ is regarded in this work as the \textit{background-traffic} term.}. 
On the other hand, dual-destination commuters are associated with two demand terms, i) $d_{o,s}^S(k;\boldsymbol{\xi})$ (veh), which expresses the number of vehicles that enter the network at time-step $k$ from origin $o\in\mathcal{O}$ towards school $s\in\mathcal{S}_b, b\in\mathcal{B}$, ii) $d_{b,d}^{S\rightarrow W}(k)$ (veh), representing the aggregate demand of dual-destination commuters that arrive at a school in $b\in\mathcal{B}$, at time step $k$ and are later heading towards the corresponding destination $d\in\mathcal{D}$. Note that the demand terms $d_{o,s}^S(k;\boldsymbol{\xi})$ and $d_{o,d}^W(k;\boldsymbol{\xi})$ explicitly depend on the Upper-Level decision vector $\boldsymbol{\xi}$, whereas the demand term $d_{b,d}^{S\rightarrow W}(k)$ does not explicitly appear as a function of $\boldsymbol{\xi}$; however, it is implicitly influenced by the Upper-Level vector $\boldsymbol{\xi}$, since the departures toward work destinations depend on the realized arrivals at schools. \par

The Upper-Level problem cannot be evaluated independently of aggregate travel demand responses, since the temporal distribution of traffic demand is an endogenous outcome of OD-level departure time adjustments. In particular, changes in school start times affect congestion patterns and travel costs, which in turn determine the departure time choices of both classes of commuters. Accordingly, for any given Upper-Level decision vector $\boldsymbol{\xi}$, the demand profiles $d_{o,s}^S(k;\boldsymbol{\xi})$ and $d_{o,d}^W(k;\boldsymbol{\xi})$ are implicitly defined by the solution of the Lower-Level problem described in Section~\ref{sec:lower_level}.\par

Commuters traveling to destination $d$ are assumed to either follow a fixed work start time $\bar{t}_d$ or operate within a flexible work arrival time window $[\bar{t}_d^l, \bar{t}_d^u]$, where $\bar{t}_d^l$, $\bar{t}_d^u$ denote the lower and upper-bounds on acceptable arrival times, respectively. For each destination $d \in \mathcal{D}$, the work-related demand $d_{o,d}^W(k;\boldsymbol{\xi})$, $o \in \mathcal{O}, k \in \mathcal{K}$, is decomposed into a i) fixed-schedule and a ii) flexible-schedule component such that
\begin{align}
	\label{eq:fixednew}
	d_{o,d}^{W,fix}(k;\boldsymbol{\xi}) &= \eta_d\, d_{o,d}^W(k;\boldsymbol{\xi}), \\ 
	\label{eq:fixednew2}
	d_{o,d}^{W,flex}(k;\boldsymbol{\xi}) &= (1-\eta_d)\, d_{o,d}^W(k;\boldsymbol{\xi}),
\end{align}

\noindent where $\eta_d \in [0,1]$ represents the proportion of commuters at destination $d$ who adhere to fixed work start times. It is important to highlight that, although school start time decisions directly reshape the school-related demand $d_{o,s}^S(k;\boldsymbol{\xi})$ by altering the temporal distribution of arrivals at schools, their impact is not limited to school-bound trips. In particular, work-related demand $d_{o,d}^W(k;\boldsymbol{\xi})$ is indirectly affected by the school start times, as well as from congestion-induced Lower-Level demand responses. Consequently, both demand components are treated as endogenous functions of the school start time decision vector $\boldsymbol{\xi}$. 

\subsubsection{Macroscopic Traffic Flow Model}
\label{sec:tbyjhgt}

\begin{figure}[t]
	\centering
	\includegraphics[width=9cm]{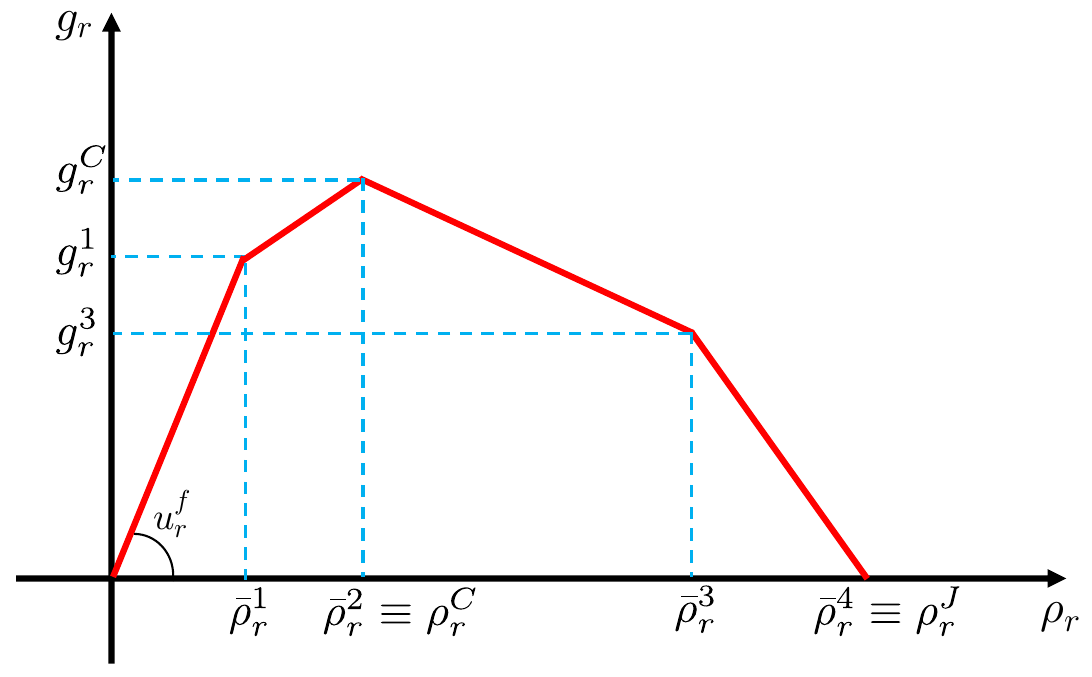} 
	\caption{Piecewise linear Macroscopic Fundamental Diagram containing four different traffic regimes.}
	\label{fig:piecewise_MFD}
\end{figure}

Each region $r\in\mathcal{R}$ in the transportation network is described by a piecewise linear MFD consisting of $N$ linear segments\footnote{The piecewise linear MFD is an appealing alternative to a unimodal MFD curve that can approximate the latter to an arbitrary accuracy, and it has been used in the past for different urban traffic control applications \citep{5625161, doi:10.3141/2390-01, 6847695}.}. Fig. \ref{fig:piecewise_MFD} shows a typical case with $N=4$ linear segments\footnote{Although we show a 4-segment MFD in Fig. \ref{fig:piecewise_MFD}, our solution methodology is applicable to any number of piecewise linear segments.}. The flow of vehicles, $g_r(\rho_r(k))$ is linked with the traffic density, $\rho_{r}(k)$ (veh/km) associated with the $l$-th linear segment of region $r$ at time-step $k$ through the relationship
\begin{align}
	\label{mfd}
	g_r(\rho_r(k)) =~ c_{r,l}\rho_r(k) + b_{r,l}, \quad \rho_r(k)\in[\bar{\rho}_{r}^l, \bar{\rho}_{r}^{l+1}], ~l=1,\ldots,N+1,
\end{align}

\noindent where $c_{r,l}, b_{r,l}$ represent the slope and the intercept of $g_r(\rho_r(k))$ for the $l$-th linear segment, respectively. Parameters $\bar{\rho}_{r}^l, l=1,\ldots,N+1$ represent the breakpoints that partition the density domain into $N$ piecewise-linear segments. Furthermore, parameters $\rho_r^J$, $\rho_r^C$ and $g_r^C$ denote the \textit{jam density}, the \textit{critical density} and the \textit{maximum outflow} generated from the MFD of region $r$, respectively. The critical density $\rho_r^C$ separates the traffic operation into the \textit{free-flow} and the \textit{congested} regime of the MFD. Leveraging the flow-density MFD relationship together with the fundamental property, in which $g_r(\rho_r(k)) = \rho_r(k) u_r(\rho_r(k))$, the speed is obtained as
\begin{equation}
	\label{speed}
	u_r(\rho_r(k))=g_r(\rho_r(k))/\rho_r(k).
\end{equation}

Another important quantity is the \textit{intended outflow} $q_r(\rho_r(k))$, which denotes the total flow of vehicles that transit to the neighbouring regions of region $r$ when the interchanged flows between regions are not restricted from their inter-boundary capacity limitations. The intended outflow $q_r(\rho_r(k))$ is linked with $g_r(\rho_r(k))$ through the relationship
\begin{equation}\label{poeee}
	q_r(\rho_r(k)) = g_r(\rho_r(k))\frac{L_r}{l_r},
\end{equation}

\noindent where $L_r$ and $l_r$ (km) denote the total length and the average trip length of vehicles in region $r$, respectively. 

Let $\mathcal{J}_r^+\subseteq \mathcal{R}$ denote the set of neighbouring regions that can receive flow from region $r\in\mathcal{R}$ and similarly let $\mathcal{J}_r^-\subseteq \mathcal{R}$ denote the set of regions that can transmit flow to their neighbouring region. The inter-boundary capacity, $C_{r,j}(\rho_j(k))$, represents the maximum flow that can be exchanged between region $r$ and its neighbouring region $j\in \mathcal{J}_r^-$. Mathematically, it is expressed as
\begin{equation}
	\label{upoptos}
	\displaystyle   C_{r,j}(\rho_j(k)) = \left\{
	\begin{array}{lll}
		\displaystyle	\hspace{-0.2cm}C_{r,j}^{\text{MAX}} , &\textrm{if}~ \rho_j(k) \leq \alpha_{r,j} \rho_j^J,  \\[6pt]
		\displaystyle	\hspace{-0.2cm}\frac{C_{r,j}^{\text{MAX}}}{1-\alpha_{r,j}}\Bigg(1-\frac{\rho_j(k)}{\rho_j^J}\Bigg), & \textrm{otherwise},   
	\end{array} 
	\right. 
\end{equation}

\noindent where $C_{r,j}^{\text{MAX}}$ is the maximum inter-boundary capacity and $\alpha_{r,j} \rho_j^J$ is the point where the inter-boundary capacity starts to decrease with $0<\alpha_{r,j}<1,r\in\mathcal{R},j\in\mathcal{J}_r^-$, similar to \cite{Sirmatel2017}.

\subsubsection{Class-specific Traffic Dynamics}
\label{sec:couplinrrrg}

To represent class-specific dynamics we adapt the accumulation-based regional model. Let set $\mathcal{Y} = \{S,W\}$ contain the indices corresponding to the two classes of commuters. Specifically, $q_{o,r,s}^S(\rho_r(k))$ and $q_{o,r,d}^W(\rho_r(k))$ (veh/h) denote the intended transfer flow of commuters, originating from $o\in \mathcal{O}$ that exit from region $r\in \mathcal{R}$ at time-step $k$ heading to school $s\in\mathcal{S}_b$, $b\in\mathcal{B}$ and at destination $d\in\mathcal{D}$ for commuters of class $y\in\mathcal{Y}=\{S,W\}$, respectively. Furthermore, let variables $\rho_{o,r,s}^{S}(k)$ and $\rho_{o,r,d}^{W}(k)$ (veh/km) indicate the density that is in region $r\in\mathcal{R}$ at time step $k$ originating from $o\in\mathcal{O}$ heading to school $s\in\mathcal{S}_b, b\in\mathcal{B}$ and to the works located at destination $d\in\mathcal{D}$ for commuters of class $y\in\mathcal{Y}$, respectively. Then we can define the following quantities
\begin{align}
	\label{rrvtvv}
	\rho_r(k) &= \sum_{y\in\mathcal{Y}} \rho_r^y(k),\\
	%\textcolor{blue}{q_r(\rho_r(k))} &= \textcolor{blue}{q_r^{W}(\rho_r(k)) + q_r^H(\rho_r(k)) = \sum_{y\in\mathcal{Y}} q_r^y(\rho_r(k))},\\
	\label{ffsaagautg}
	\rho_{r}^{S}(k) &= \sum_{o\in \mathcal{O}}\sum_{b\in\mathcal{B}}\sum_{s\in\mathcal{S}_b}
	\rho_{o,r,s}^{S}(k),\\
	\label{fhshshuj}
	\rho_{r}^{W}(k) &= \sum_{o\in \mathcal{O}}\sum_{d\in\mathcal{D}} \rho_{o,r,d}^{W}(k),
\end{align}

\noindent where $\rho_{r}^{y}(k)$ (veh/km) indicates the density with respect to class $y\in\mathcal{Y}$ that is in region $r\in\mathcal{R}$ at time step $k$. It is true that variable $q_{o,r,s}^S(\rho_r(k))$ is given by 
\begin{align}
	\label{inssss}
	q_{o,r,s}^{S}(\rho_r(k)) &= \frac{\rho_{o,r,s}^{S}(k)}{\rho_r(k)}q_r(\rho_r(k)) = u_r(\rho_r(k))\rho_{o,r,s}^{S}(k)\frac{L_r}{l_r}.
\end{align} 

\noindent Similarly, variable $q_{o,r,d}^W(k)$ is obtained from
\begin{align}
	\label{naruto}
	q_{o,r,d}^{W}(\rho_r(k)) &= \frac{\rho_{o,r,d}^{W}(k)}{\rho_r(k)}q_r(\rho_r(k)) = u_r(\rho_r(k))\rho_{o,r,d}^{W}(k)\frac{L_r}{l_r}.
\end{align} 

\noindent Let also $q_{o,r,j,s}^{S}(\rho_r(k))$ and $\tilde{q}_{o,r,j,s}^{S}(\rho_r(k), \rho_j(k))$ (veh/h) denote the \textit{intended} and the \textit{actual} flow of commuters of class S that transit from region $r\in \mathcal{R}$ at time-step $k$ originating from $o\in \mathcal{O}$ through its neighbouring region $j\in\mathcal{J}_r^-$ and eventually heading to school $s\in\mathcal{S}_b, b\in\mathcal{B}$, respectively. Similarly, quantities $q_{o,r,j,d}^{W}(\rho_r(k))$ and $\tilde{q}_{o,r,j,d}^{W}(\rho_r(k), \rho_j(k))$ (veh/h) are defined accordingly for class W concerning destination $d\in\mathcal{D}$ as
\begin{align}
	\label{queen}
	q_{o,r,j,s}^S(\rho_r(k)) &= \theta_{r,j,b}(k)q_{o,r,s}^S(\rho_r(k)),\\
	\label{queen1}
	q_{o,r,j,d}^W(\rho_r(k)) &= \theta_{r,j,d}(k)q_{o,r,d}^W(\rho_r(k)),\\
	\label{tucer1}
	q_{o,r,s}^{S}(\rho_r(k)) &= \sum_{j\in\mathcal{J}_r^+}q_{o,r,j,s}^{S}(\rho_r(k)),\\
	q_{o,r,d}^{W}(\rho_r(k)) &= \sum_{j\in\mathcal{J}_r^+}q_{o,r,j,d}^{W}(\rho_r(k)).
\end{align}

\noindent In Eqs. \eqref{queen} - \eqref{queen1} the terms $\theta_{r,j,b}(k), \theta_{r,j,d}(k)\in[0,1]$ denote the ratio of vehicles that transit from region $r\in\mathcal{R}$ at each time-step $k$ through the neighbouring region $j\in\mathcal{J}_r^+$ towards the corresponding school $s$ located in $b\in\mathcal{B}$ or at destination $d\in\mathcal{D}$, respectively\footnote{In this work, the split ratio is considered a known regional parameter, which can be empirically estimated.}, such that $\sum_{j\in\mathcal{J}_r^+}\theta_{r,j,d}(k) = \sum_{j\in\mathcal{J}_r^+}\theta_{r,j,b}(k)=1, b\in\mathcal{B}, d\in\mathcal{D}, k\in\mathcal{K}$.\par

Let $q_{r,j}(\rho_r(k))$ (veh/h) denote the \textit{intended flow} of vehicles (in regional level) that transit from region $r\in\mathcal{R}$ at time step $k\in\mathcal{K}$ towards the directly accessible neighbouring region $j\in\mathcal{J}_r^+$. It is true that
\begin{align}
	\label{ffsarrraga}
	q_{r,j}(\rho_r(k)) = \sum_{o\in\mathcal{O}}\sum_{b\in\mathcal{B}}\sum_{s\in\mathcal{S}_b} q_{o,r,j,s}^S(\rho_r(k)) + \sum_{o\in\mathcal{O}}\sum_{d\in\mathcal{D}} q_{o,r,j,d}^W(\rho_r(k)).
\end{align}

\noindent From Eqs. \eqref{poeee} and \eqref{ffsarrraga} it is also true that
\begin{align}
	\label{interesting}
	q_r(\rho_r(k)) = \sum_{j\in\mathcal{J}_r^-}q_{r,j}(\rho_r(k)).
\end{align}

\noindent The \textit{actual flows}, $\tilde{q}_{o,r,j,s}^{S}(\rho_r(k), \rho_j(k))$ and $\tilde{q}_{o,r,j,d}^{W}(\rho_r(k), \rho_j(k))$ for commuters of class S and W, respectively, are defined as 
\begin{align}
	\label{secondtt}
	\tilde{q}_{o,r,j,s}^{S}(\rho_r(k), \rho_j(k)) = \min\Bigg(q_{o,r,j,s}^{S}(\rho_r(k)),\frac{C_{r,j}(\rho_j(k))}{q_{r,j}(\rho_r(k))}q_{o,r,j,s}^S(\rho_r(k))\Bigg),\\
	\label{secondtt1}
	\tilde{q}_{o,r,j,d}^{W}(\rho_r(k), \rho_j(k)) = \min\Bigg(q_{o,r,j,d}^{W}(\rho_r(k)),\frac{C_{r,j}(\rho_j(k))}{q_{r,j}(\rho_r(k))}q_{o,r,j,d}^W(\rho_r(k))\Bigg),
\end{align}

\noindent extending the modeling framework of \cite{menelaou2021joint} to multi-class-specific traffic dynamics. Next, the traffic dynamics with respect to commuters of Class S are given by
\begin{align}
	\label{arxontas1}
	\rho_{o,r,s}^{S}(k+1) = \rho_{o,r,s}^{S}(k) +  \frac{1}{L_r}d_{r,s}^{S}(k;\boldsymbol{\xi}) + \frac{T}{L_r}\Bigg(\sum_{j\in \mathcal{J}_r^+}\tilde{q}_{o,j,r,s}^{S}(\rho_j(k), \rho_r(k)) -\sum_{j\in \mathcal{J}_r^-}\tilde{q}_{o,r,j,s}^{S}(\rho_r(k), \rho_j(k))\Bigg),
\end{align}

\noindent where $d_{r,s}^{S}(k;\boldsymbol{\xi})\neq 0, r\in\mathcal{O}$. In Eq. \eqref{arxontas1}, $d_{r,s}^{S}(k;\boldsymbol{\xi})$ denotes the shifted demand of class S stemming from region $r$ that passes from school $s$ at time-step $k$, while the first and second summation terms denote the corresponding inflow and outflow, respectively. Eq. \eqref{arxontas1} is crucial for our analysis as it links the school-related demand term, $d_{r,s}^{S}(k;\boldsymbol{\xi})$ with the traffic density dynamics of commuters of class S.

\begin{figure}[t]
	\centering
	\hspace{-1.9cm}\includegraphics[width=9cm]{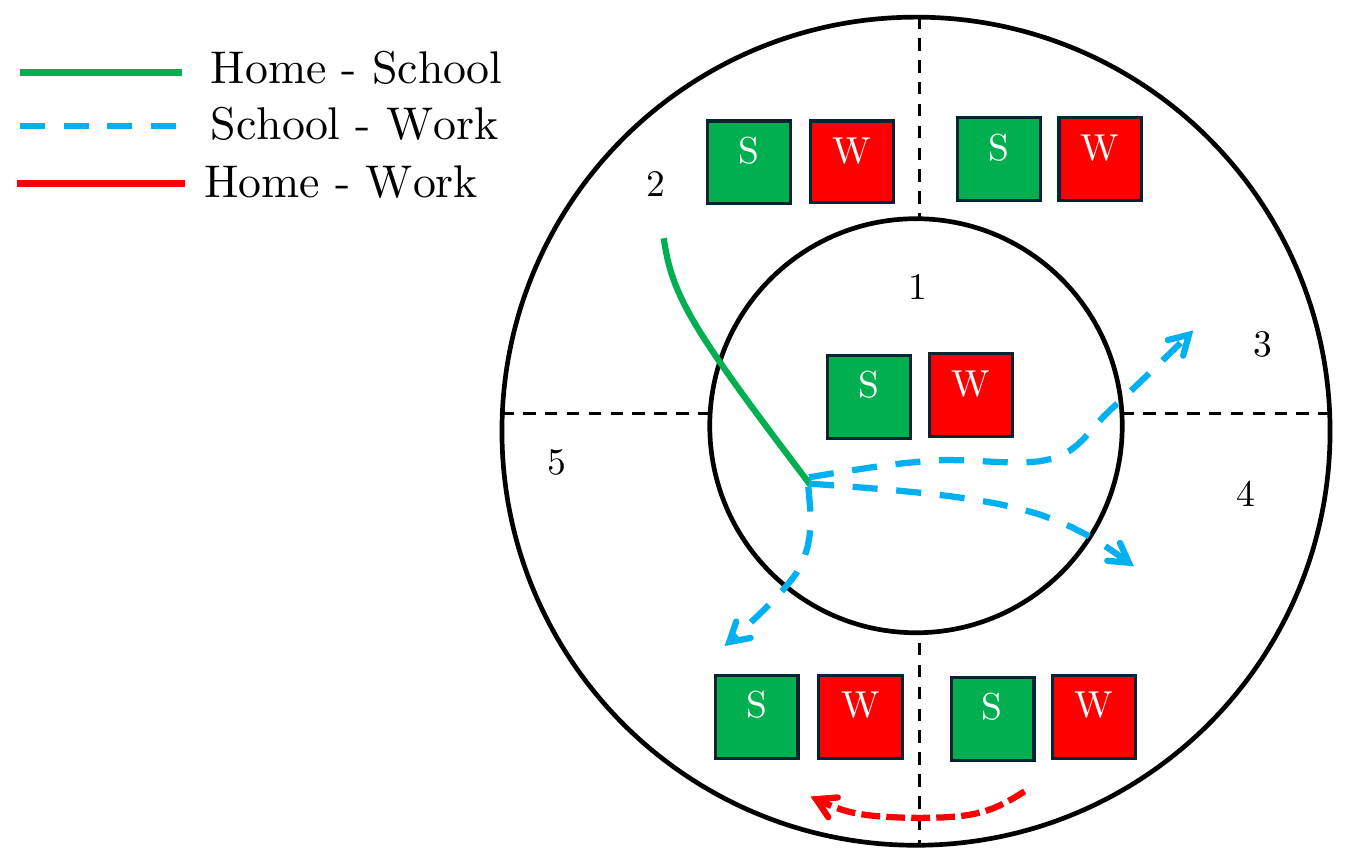} 
	\caption{The schematic of a multi-region urban network and some path examples for the two considered classes of commuters. The network consists of 5 regions, where every region may contain both, schools and workplaces. For the dual-destination commuters, a demand enters from region $r=2$ heading to a school located in region $r=1$ (represented with the green solid line). When this demand arrives at this school, then a \textit{new} demand is generated, representing the vehicles that after having arrived at this school, is subsequently dispersed to different workplaces or homes located in regions $r=\{3,4,5\}$ (this transition is shown with the blue dashed lines). For the single-destination commuters, an aggregate demand enters from region $r=4$, heading directly to workplaces located in region $r=5$ (represented with the red rounded line).}
	\label{fig:distea1}
\end{figure}

As an example, the movement of the two classes of commuters in an urban network is represented via Fig. \ref{fig:distea1}. Formally, we introduce demand $d_{b,d}^{S\rightarrow W}(k)$ (veh) to link the cumulative flow that arrives at all schools located in region $b\in\mathcal{B}$ from different origins with their final destination $d\in\mathcal{D}$. Hence,
\begin{equation}
	\label{strong2}
	d_{b,d}^{S\rightarrow W}(k) = T\sum_{s\in\mathcal{S}_b} \sum_{o\in\mathcal{O}} A_{o,s,d} \cdot \tilde{q}_{o,b,b,s}^S(\rho_b(k), \rho_b(k)),
\end{equation}

\noindent where $A_{o,s,d}$ indicates the ratio of vehicles originating from $o\in\mathcal{O}$ that travel from school $s\in\mathcal{S}_b$ to destination $d\in\mathcal{D}$, such that $\sum_{d\in\mathcal{D}}A_{o,s,d} = 1$ and $\tilde{q}_{o,b,b,s}^S(\rho_b(k), \rho_b(k))$ is derived from Eq. \eqref{secondtt} by putting $r=j=b$.\par

Without loss of generality, we assume that the demand term  $d_{b,d}^{S\rightarrow W}(k)$ is assimilated in the traffic dynamics with respect to commuters of class W as follows
\begin{align}
	\label{arxontase}
	\rho_{o,r,d}^{W}(k+1) = \rho_{o,r,d}^{W}(k) + \frac{1}{L_r}\Big(d_{r,d}^{W}(k;\boldsymbol{\xi}) + 
	d_{r,d}^{S\rightarrow W}(k)\Big) +\frac{T}{L_r}\Bigg(\sum_{j\in \mathcal{J}_r^+}\tilde{q}_{o,j,r,d}^{W}(\rho_j(k), \rho_r(k)) -\sum_{j\in \mathcal{J}_r^-}\tilde{q}_{o,r,j,d}^{W}(\rho_r(k), \rho_j(k))\Bigg),
\end{align}

\noindent where $d_{r,d}^{W}(k;\boldsymbol{\xi})\neq 0, ~\text{if}~ r\in\mathcal{O}, d_{r,d}^{S\rightarrow W}(k)\neq 0,~\text{if}~ r\in\mathcal{B}$. In Eq. \eqref{arxontase}, the first demand term, $d_{r,d}^{W}(k;\boldsymbol{\xi})$ denotes the aggregate demand of vehicles belonging to class W that enter to region $r$ at time-step $k$ heading to the works located in destination $d\in\mathcal{D}$, while the second demand term $d_{r,d}^{S\rightarrow W}(k)$ is defined in Eq. \eqref{strong2} with $r\equiv b$. 
We remind the reader that the work-related demand term $d_{r,d}^{W}(k;\boldsymbol{\xi})$ consists of the fixed and the flexible schedule component terms, $d_{r,d}^{W,fix}(k;\boldsymbol{\xi}), d_{r,d}^{W,flex}(k;\boldsymbol{\xi})$ as defined in Eqs. \eqref{eq:fixednew}, \eqref{eq:fixednew2}, respectively. The first and second summation terms denote the corresponding inflow and outflow to and from region $r\in\mathcal{R}$, respectively.\par 

We emphasize that the two demand terms appearing in Eq. \eqref{arxontase} play conceptually different roles. The work-related demand $d_{r,d}^{W}(k;\boldsymbol{\xi})$ represents the exogenous inflow of Class~W commuters entering the network and depends explicitly on the Upper-Level decision vector $\boldsymbol{\xi}$ through travelers' departure time choices. On the contrary, the term $d_{r,d}^{S \rightarrow W}(k)$ is not explicitly parameterized by $\boldsymbol{\xi}$. Instead, it represents an endogenous transfer flow generated by commuters of Class~S after arriving at region $r\in\mathcal{B}$, heading subsequently to their destinations. Nevertheless, the term $d_{r,d}^{S \rightarrow W}(k)$ is indirectly influenced by $\boldsymbol{\xi}$, since the school schedules determine when commuters of Class S arrive at school and subsequently later head to their destinations. Thus, the Upper-Level decision vector $\boldsymbol{\xi}$ essentially affects the travel behavior of both classes of commuters through explicit and implicit demand mechanisms.\par 

\subsubsection{Problem Formulation}
\label{sec:problem_formulation}

To formally describe the Upper-Level problem, we first define the objective cost function, $J_{\text{TTS}}$ (veh h). Toward this direction, we introduce variables $S_a(k)$ and $S_b(k)$, which denote the cumulative number of vehicles that entered the network and successfully exited from their destination $d$, respectively, such that 
\begin{align}
	\label{cumulative_in}
	&S_a(k+1) = S_a(k) + \sum_{o\in \mathcal{O}}\sum_{b\in\mathcal{B}}\sum_{s\in\mathcal{S}_b} d^S_{o,s}(k;\boldsymbol{\xi}) + \sum_{o\in\mathcal{O}}\sum_{d\in\mathcal{D}} d_{o,d}^{W}(k;\boldsymbol{\xi}) + \sum_{b\in \mathcal{B}}\sum_{d\in\mathcal{D}}d_{b,d}^{S\rightarrow W}(k), S_a(0) = 0, k\in\mathcal{K}, \\
	\label{final}
	&S_b(k+1) = S_b(k) + T\Bigg(\sum_{o\in \mathcal{O}}\sum_{b\in\mathcal{B}}\sum_{s\in\mathcal{S}_b}\tilde{q}_{o,b,b,s}^S(k)
	+\sum_{o\in \mathcal{O}}\sum_{d\in\mathcal{D}}\tilde{q}_{o,d,d,d}^W(k) \Bigg),S_b(0) = 0, k\in\mathcal{K},
\end{align}

\noindent where variables $\tilde{q}_{o,b,b,s}^S(k)$ and  $\tilde{q}_{o,d,d,d}^W(k)$ determine the total number of vehicles that arrived at school $s$ in the intermediate destination $b\in\mathcal{B}$ and the final destination $d\in\mathcal{D}$ at time step $k$, respectively. Hence,
\begin{align}
	\label{fixed}
	J_{\textrm{TTS}} = T\sum_{k\in\mathcal{K}}\Big(S_a(k) - S_b(k)\Big).
	%= T_s\cdot\sum_{k}\Bigg(\sum_{o}\sum_{z}^{}\sum_{v}^{}\sum_{d} d_{o,z,v,d}^{sw}(k,\tau_{z}^*) \\
	%&+ \sum_{o}\sum_{v}^{}\sum_{d} d_{o,v,d}^{w}(k,t_{v}^*) - T_s\sum_{d\in D}\tilde{q}_{d,d,d}^{sw}(k) - T_s \sum_{d\in D}\tilde{q}_{d,d,d}^{w}(k) \Bigg)
\end{align}

\noindent The second objective, $J_{\text{STC}}$ (min) is defined as
\begin{align}
	\label{QUO1}
	J_{\text{STC}} = \displaystyle	 \sum_{b\in\mathcal{B}}\sum_{s\in\mathcal{S}_b}\sum_{m\in\mathcal{M}}\xi_{m,s} |m| V,
\end{align}

\noindent where $\xi_{m,s}|m|V$ is equal to $|m|V$ if $\xi_{m,s}$ is equal to one and zero otherwise. This objective captures the total amount of shifting applied to school start times. On the one hand, shifting the start time of certain schools to a different time instant might be advantageous in terms of total time spent. On the other hand, such a shift would result in greater deviations from the originally scheduled school start times, which may be undesirable from an operational standpoint. In mathematical programming terms, the Upper-Level optimization problem can be stated as
\begin{subequations}
	\label{rfwgfa}
	\begin{align}
		\label{trr}
		(P_0) \quad \underset{\boldsymbol{\xi}\in\Xi_U}{\text{Minimize}} &~ \displaystyle \Big\{\hspace{-0.1cm}J_{\text{TTS}}, J_{\text{STC}}\hspace{-0.1cm}\Big\} \\
		\textrm{Subject To:}&~~\textrm{Traffic Dynamics}~ \eqref{mfd} -  \eqref{QUO1}, \nonumber \\
		\label{bin}
		&~~\sum_{m\in\mathcal{M}}\xi_{m,s}=1, \quad \forall s\in\mathcal{S}_b, \forall b\in\mathcal{B}, \\
		\label{jam}
		& ~~ 0\leq \rho_r(k) \leq \rho_r^J, \quad \forall r\in \mathcal{R},~\forall k\in\mathcal{K},\\
		\label{binary}
		&~~\xi_{m,s}\in\{0,1\}, \quad \forall m\in\mathcal{M},~\forall s\in\mathcal{S}_b, \forall b\in\mathcal{B},\\
		\label{init}
		\text{Initialization:}&~~ \rho_r(0)=\bar{\rho}_r,~\forall r\in\mathcal{R},\nonumber\\
		&~~\rho_r^y(0)=\bar{\rho}_r^y,~ \forall r\in\mathcal{R},\forall y\in\mathcal{Y},\nonumber\\
		&~~\rho_{o,r,s}^S(0)=\bar{\rho}_{o,r,s}^S,~\forall o\in\mathcal{O}, \forall r\in\mathcal{R},~\forall s\in\mathcal{S}_b, \forall b\in\mathcal{B},\nonumber\\ &~~\rho_{o,r,d}^W(0)=\bar{\rho}_{o,r,d}^W,~\forall o\in\mathcal{O}, \forall r\in\mathcal{R},~\forall  d\in\mathcal{D},\\
		\label{input}
		\text{Input:} 
		&~~ d_{o,s}^S(k;\boldsymbol{\xi}), 
		\quad \forall o\in\mathcal{O},~\forall s\in\mathcal{S}_b,~\forall b\in\mathcal{B},~\forall k\in\mathcal{K}, \nonumber\\
		&~~ d_{o,d}^W(k;\boldsymbol{\xi}), 
		\quad \forall o\in\mathcal{O},~\forall d\in\mathcal{D},~\forall k\in\mathcal{K}, \nonumber\\
		&~~\theta_{r,j,d}(k)\in[0,1],
		\quad \forall r\in\mathcal{R},~\forall j\in\mathcal{J}_r^-,~\forall k\in\mathcal{K},~\forall d\in\mathcal{D},\nonumber\\
		&~~\theta_{r,j,b}(k)\in[0,1],
		\quad \forall r\in\mathcal{R},~\forall j\in\mathcal{J}_r^-,~\forall k\in\mathcal{K},~\forall b\in\mathcal{B}.
	\end{align} 
\end{subequations}

\noindent Problem \eqref{rfwgfa} expresses the derived formulation, where constraints \eqref{mfd} - \eqref{QUO1} model the traffic dynamics according to a piecewise linear MFD, while constraint \eqref{bin} ensures that only one shifting index $m\in\mathcal{M}$ related to the shifted school start time can be assigned to each school $s\in\mathcal{S}_b, b\in\mathcal{B}$. Constraint \eqref{jam} sustains the density of each region within its physical limits, while constraint \eqref{binary} states that variable $\xi_{m,s}, \forall m\in\mathcal{M}, \forall s\in\mathcal{S}_b, \forall b\in\mathcal{B}$ can only take binary values 0 or 1. Next, constraint \eqref{init} represents the initial state of the network.\par 

For any feasible school start time configuration $\boldsymbol{\xi}\in\Xi_U$, the Lower-Level problem, defined in Section \ref{sec:lower_level} determines the corresponding demand profiles, i) the school-related demand $d_{o,s}^S(k;\boldsymbol{\xi})$ and ii) the work-related demand $d_{o,d}^W(k;\boldsymbol{\xi})$. These equilibrium demand profiles are treated as exogenous inputs to Problem~$P_0$ (see constraint \eqref{input}) and are used to evaluate the macroscopic traffic dynamics, the resulting density trajectories, and the corresponding values of the Upper-Level objective functions. The dependence of the demand terms on $\boldsymbol{\xi}$ reflects the behavioral response
of commuters to the implemented school start time configuration. The output of the optimization problem $P_0$ gives us the values for the binary variables $\xi_{m,s}, \forall m\in\mathcal{M}, \forall s\in\mathcal{S}_b, \forall b\in\mathcal{B}$, meaning that the start time of each school would be $\tilde{\tau}_s = \tau_s + \xi_{m,s}mV$. \par

Given the demand vectors $\mathbf{d}_{o,s}^S(\boldsymbol{\xi})\in\mathbb{R}^{K\times 1}$, $\forall o\in\mathcal{O},\forall s\in\mathcal{S}_b,\forall b\in\mathcal{B}$ and $\mathbf{d}_{o,d}^W(\boldsymbol{\xi})\in\mathbb{R}^{K\times 1}$, $\forall o\in\mathcal{O},\forall d\in\mathcal{D}$, the Upper-Level problem $P_0$ represents a Bi-Objective Mixed-Integer Nonlinear Program (MINLP) due to the existence of the nonlinear MFD function shown in \eqref{mfd}, the bilinear terms in \eqref{inssss} and \eqref{naruto}, the nonlinear functions in \eqref{upoptos}, \eqref{secondtt} and \eqref{secondtt1}, the binary variables in \eqref{binary} and the incorporation of two objective metrics of interest shown in \eqref{trr}. Efficiently solving problem $P_0$ with standard nonlinear solvers can be quite challenging since nonconvex and nonlinear constraints arise with respect to the traffic dynamics. To add to that, the presence of the binary variables along with the minimization with respect to two objective metrics, i.e., $J_{\text{TTS}}, J_{\text{STC}}$ in problem $P_0$ make the problem combinatorial and hence more difficult to solve. As a result, $P_0$ constitutes a challenging optimization problem, even when the demand vectors $\mathbf{d}_{o,s}^S(\boldsymbol{\xi})$ and $\mathbf{d}_{o,d}^W(\boldsymbol{\xi})$ are known.

\subsection{Lower-Level problem}
\label{sec:lower_level}

The Upper-Level problem selects the school start time decision vector $\boldsymbol{\xi}$. However, the impact of these decisions depends critically on how commuters respond to the resulting traffic conditions. In particular, changes in school start times affect congestion patterns across the network, which in turn influence the departure time choices of both school- and work-bound commuters. Capturing this endogenous behavioral response is therefore essential for evaluating the effectiveness of these Upper-Level decisions.\par

The Lower-Level problem is designed to model this response by computing the equilibrium allocation of departure times across the analysis horizon for all OD pairs. At the level of spatial and temporal aggregation considered in this work, departure-time choice constitutes the primary behavioral response to changes in system conditions, while route choice has a comparatively limited impact on the evolution of system-wide regional congestion patterns and hence is beyond the scope of this work \citep{ARNOTT1990111, small1982scheduling}.\par %As a result, the Lower-Level problem aims to distribute demand across time, rather than selecting alternative routes.\par

%Route choice is implicitly determined by the acyclic regional graph and the associated destination-specific split ratios introduced in Section~\ref{sec:path_identification}. This assumption is consistent with macroscopic MFD-based traffic modeling approaches, where traffic dynamics are described at an aggregate regional level and the detailed structure of individual routes is not explicitly represented \citep{GEROLIMINIS2008759, doi:10.3141/2124-12, gu2018big}. Furthermore, at the level of spatial and temporal aggregation considered in this work, departure-time choice constitutes the primary behavioral response to changes in system conditions, while route switching has a comparatively limited impact on the evolution of regional congestion patterns \citep{ARNOTT1990111, small1982scheduling}. Fixing the split ratios therefore allows the model to capture the dominant congestion formation mechanisms induced by school start time policies, while preserving computational tractability. Accordingly, all commuters sharing the same OD pair follow the regional transitions induced by the split-ratio structure.
%The only behavioral decision considered at the Lower-Level is the departure time choice. Thus, the Lower-Level equilibrium allocates fixed OD demand across departure times, rather than across alternative routes. 

For a given Upper-Level decision vector $\boldsymbol{\xi}\in\Xi_U$, the Lower-Level problem determines an equilibrium-based demand profile denoted by $\hat{\mathbf d}(\boldsymbol{\xi})$. Specifically, $\hat{\mathbf d}(\boldsymbol{\xi})$ comprises the collection of the equilibrium demand vectors for the two classes of commuters over the entire analysis horizon, namely the school-related demand $\hat{\mathbf d}_{o,s}^{S}(\boldsymbol{\xi})$ and the work-related demand $\hat{\mathbf d}_{o,d}^{W}(\boldsymbol{\xi})$, $\forall o\in\mathcal{O}$, $\forall s\in\mathcal{S}_b$, $\forall b\in\mathcal{B}, \forall d\in\mathcal{D}$. Throughout the remainder of this section, $\boldsymbol{\xi}$ is treated as fixed and its dependence is omitted from the demand variables whenever this does not create ambiguity.

\subsubsection{Identification and Characterization of Paths: Split-Ratio Representation}
\label{sec:path_identification}

In MFD-based dynamic traffic assignment frameworks, a clear distinction is made between \emph{trips} and \emph{paths}, reflecting the different levels of network resolution. As discussed in \cite{BatistaSeppecherLeclercq2021a}, trips are defined as ordered sequences of links connecting origin and destination nodes in the urban network, whereas \emph{paths} correspond to ordered sequences of regions connecting an origin region to a destination region in the aggregated (regional) network. This distinction is fundamental in regional-scale modeling, where individual vehicle trajectories are not explicitly represented and traffic dynamics are governed by region-level states.

In the existing literature, regional paths are typically identified by aggregating urban-scale trips according to the sequence of regions they traverse, thereby constructing origin--destination connections at the regional level \citep{BatistaGeroliminisLeclercq2019}. Such paths provide a macroscopic description of how demand propagates across regions while remaining consistent with MFD-based traffic dynamics \citep{BatistaLeclercq2020}. 

In the present work, instead of explicitly identifying regional paths from disaggregated trip data, we adopt an equivalent but implicit representation, based on destination-specific split ratios, fully consistent with the macroscopic traffic model introduced in Section~\ref{sec:tbyjhgt}. Specifically, the regional network is represented as a directed graph whose nodes correspond to regions and whose arcs represent direct accessibility between neighboring regions. For each origin region $o$ and destination $d$, the split ratio $\theta_{o,j,d}$ incorporated in Eqs. \eqref{queen}, \eqref{queen1}, represents the fraction of vehicles transiting from origin $o\in\mathcal{O}$ through the adjacent region $j\in\mathcal{J}_r^+$ heading to destination $d\in\mathcal{D}$, satisfying the condition $\sum_{j \in J_r^+} \theta_{o,j,d} = 1, \forall o\in\mathcal{O},\; \forall d\in\mathcal{D}$. The split ratios therefore define the routing structure of the regional network. 

We assume that for each origin-destination pair $(o,d)$, the region-to-region transitions associated with nonzero split ratios $\theta_{o,j,d}$ define an \emph{acyclic directed graph}. This assumption excludes circular movements at the regional level and ensures that any vehicle departing from region $r$ reaches its destination region $d$ in a finite number of regional transition steps. This assumption is consistent with macroscopic MFD-based traffic modeling approaches, where traffic dynamics are described at an aggregate regional level and the detailed structure of individual routes is not explicitly represented \citep{GEROLIMINIS2008759, doi:10.3141/2124-12, gu2018big}.\par 

Under this assumption,  routing is fully determined by the split-ratio structure. Each admissible sequence of regions defined by the acyclic graph convention corresponds to a feasible regional path connecting an origin to a destination region. Consequently, routing is fixed at the Lower-Level problem and commuters do not switch routes. This modeling choice ensures consistency with the aggregated MFD-based representation and allows us to focus exclusively on departure-time choice decisions.\par

\begin{figure}[t]
	\centering 
	\includegraphics[width=5.5cm]{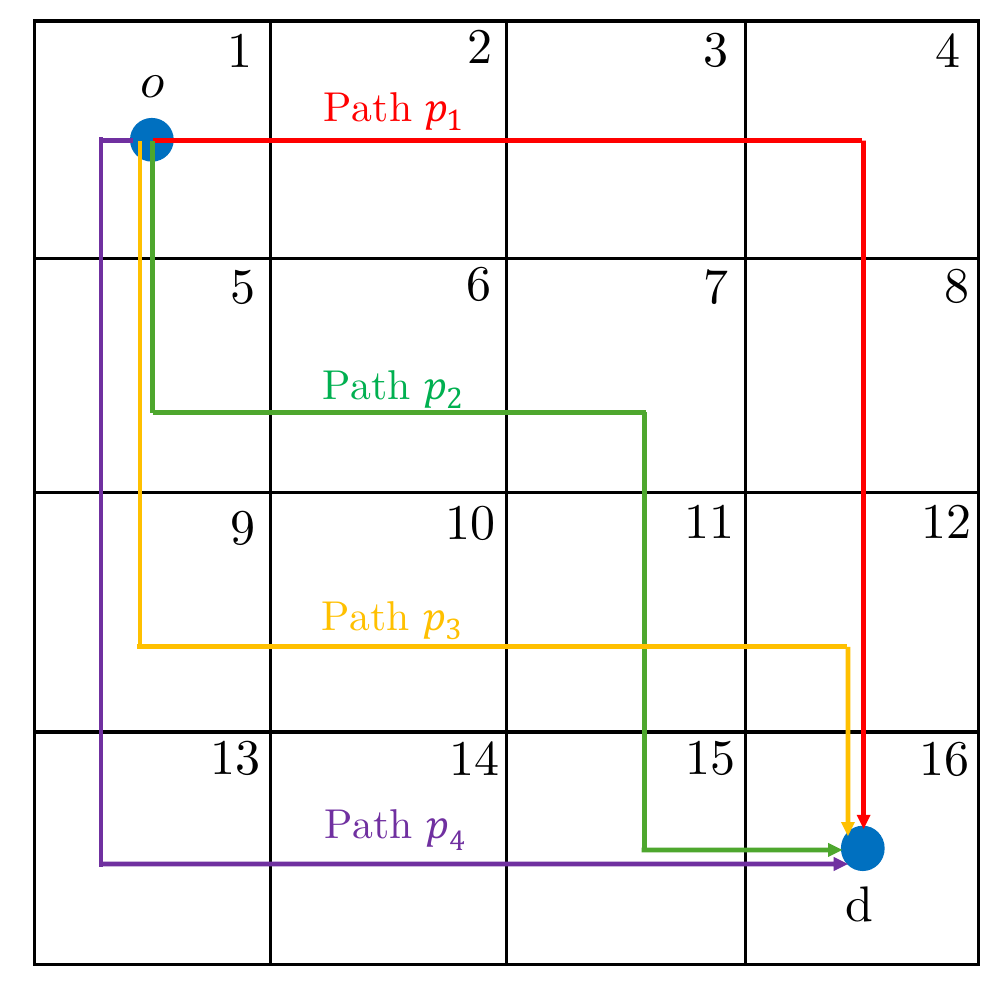} 
	\caption{Illustrative regional path selection process from origin $o=1$ to destination region $d=16$ consisting of 4 indicative paths. 
		Each arrow represents a regional transition with a destination-specific split ratio $\theta_{r,j,d}$ with respect to region $r$, the directly accessible neighbouring region $j$ and the destination region $d$. 
		The probability $\varphi_{p_1,1,16}$ of selecting path $p_1$ is obtained as the product of the split ratios along the path, namely
		$\varphi_{p_1,1,16} = \theta_{1,2,16} \cdot \theta_{2,3,16} \cdot \theta_{3,4,16} \cdot \theta_{4,8,16} \cdot \theta_{8,12,16} \cdot \theta_{12,16,16}$. Similarly for path $p_2$, we obtain the probability $\varphi_{p_2,1,16} = \theta_{1,5,16} \cdot \theta_{5,6,16} \cdot \theta_{6,7,16} \cdot \theta_{7,11,16} \cdot \theta_{11,15,16} \cdot \theta_{15,16,16}$}
	\label{fig:regional_path_probability}
\end{figure}

Fig.~\ref{fig:regional_path_probability} schematically illustrates an example of how regional transitions define feasible paths in a network. Each path connecting an origin to a destination region corresponds to a sequence of transitions, with the probability of selecting that path determined by the multiplicative combination of the split ratios along the sequence.

\subsubsection{Travel Times on regional networks}

Let $\mathcal{Q}$ denote the set of all origin-destination (OD) pairs considered in the network. Given the acyclic regional graph structure and the induced set of feasible regional paths $\mathcal{P}_{o,d}$ for each origin-destination pair $(o,d)\in\mathcal{Q}$, we first define
the travel time associated with a \emph{fixed regional path} and then derive the expected OD-level travel time. Consider a path $p \in \mathcal{P}_{o,d}$, defined as an ordered sequence of regions connecting an origin region with a destination region. Following the regional MFD-based formulation in \cite{BATISTA2025104980}, the travel time $TT_{p,o,d}(k)$ (measured in hours) of path $p$ of $(o,d)$ pair that departs at time-step $k$ is defined as
\begin{equation}
	\label{trave_time_equation}
	TT_{p,o,d}(k)
	=
	\sum_{r\in\mathcal{R}}
	\frac{L_{r,p}}{u_r(\rho_r(k))} \, \delta_{r,p},
	\qquad
	\forall p\in\mathcal{P}_{o,d}, \forall (o,d)\in\mathcal{Q},
\end{equation}
where $L_{r,p}$ denotes the travel distance associated with region $r$ along path $p$, $u_r(\rho_r(k))$ is the MFD-based speed in region $r$ at time-step $k$, and $\delta_{r,p}$ is a binary variable that is equal to 1 if path $p$ traverses region $r$, and 0 otherwise.
This expression yields the deterministic travel time of path $p$, conditional on the current regional traffic states. For the expression shown in Eq. \eqref{trave_time_equation}, we adopt a \emph{quasi-static assumption}, in which the regional densities $\rho_r(k)$ are assumed to remain fixed during the travel of path $p$. This simplification allows us to compute deterministic travel times based on the current traffic state. A more detailed discussion with respect to this assumption and its implications is provided in Section~\ref{sec:formulation_ddue}.\par

We reiterate that routing decisions are implicitly determined by the destination-specific split ratios that govern the regional transitions. Consequently, the probability that a trip traveling between origin $o\in\mathcal{O}$ and destination $d\in\mathcal{D}$ follows a particular regional path $p$ can be obtained as the product of the transition probabilities along that path. Let path $p$ correspond to the ordered sequence of regional transitions
\[
p=(r_0,r_1,r_2,\ldots,r_n),
\]
where $r_0=o$ and $r_n=d$. With the aid of the split ratio notion, the probability that a trip between $o$ and $d$ follows path $p$ is therefore expressed as
\begin{equation}
	\label{eq:path_probability}
	\varphi_{p,o,d}
	=
	\prod_{i=0}^{n-1}\theta_{r_i,r_{i+1},d},
	\qquad
	\forall p\in\mathcal{P}_{o,d},
\end{equation}

\noindent such that
\begin{equation*}
	\sum_{p \in \mathcal{P}_{o,d}} \varphi_{p,o,d} = 1, \qquad \forall (o,d) \in \mathcal{Q}.
\end{equation*}

\noindent These path probabilities define how the demand associated with pair $(o,d)$ is distributed across the feasible regional paths induced by the split-ratio structure. Accordingly, the expected travel time, denoted by $\overline{TT}_{o,d}(k)$ (measured in hours) at time-step $k$, is obtained as the probability-weighted average of the travel times associated with all feasible regional paths connecting the pair $(o,d)$:
\begin{equation}
	\label{eq:travel_time_expression}
	\overline{TT}_{o,d}(k)
	=
	\sum_{p \in \mathcal{P}_{o,d}}
	\varphi_{p,o,d}\,TT_{p,o,d}(k),
	\qquad
	\forall (o,d)\in\mathcal{Q}.
\end{equation}

%\textcolor{blue}{\textbf{16/03/2025 say that route choice is beyond the scope of this work. Justify!}}\\

\subsubsection{Departure Time Choice Modeling}
\label{sec:depart_complete}

This section defines the departure time choice model, which captures the behavioral response of commuters to school start times, fixed work hours, and flexible work hours.
We develop a modified variant of the classic bottleneck model for the morning commute problem,~\citep{vickrey1969congestion}, to account for the co-existence of different classes of commuters in an urban network considering regional MFD-based traffic dynamics. 

Let $t_k = t_0 + k T$ denote the departure time (in hours) corresponding to discrete time-step $k$, where $t_0$ is the starting time of the analysis horizon coinciding with $k=0$ and $T$ is the duration of the discrete time-step $k$. Desired arrival times, expressed in hours, are denoted by $\tilde{\tau}_s$ for schools, $\bar{t}_d$ for fixed work, and $\bar{t}_d^{flex}$ for flexible work. For schools, $\tilde{\tau}_s$ depends on the Upper-Level decision vector $\boldsymbol{\xi}$, i.e., $\tilde{\tau}_s = \tau_s + \xi_{m,s} m V$, 
where $\tau_s$ is the initial school start time and $\xi_{m,s}, m\in\mathcal{M}, s\in\mathcal{S}_b, b\in\mathcal{B}$ denotes the binary variables responsible for the school start time shift.

Each OD pair is associated with a perceived travel cost when departing at $t_k$. We adapt the $\alpha$-$\beta$-$\gamma$ preference model of \cite{vickrey1969congestion} (generalized by \cite{ARNOTT1990111}) to the OD-pair level to capture aggregate departure time preferences. The perceived cost function for each class of commuters at the regional OD-pair level is expressed as
\begin{align}
	\label{cost1}
	&P_{o,s,b}^{S}(k) = \alpha^S \cdot \overline{TT}_{o,b}(k)
	+ \beta^S \cdot \max\big(0, \tilde{\tau}_s - t_k - \overline{TT}_{o,b}(k) \big)
	+ \gamma^S \cdot \max\big(0, t_k + \overline{TT}_{o,b}(k) - \tilde{\tau}_s \big), \\
	\label{cost2}
	&P_{o,d}^{W}(k) = \alpha^W \cdot \overline{TT}_{o,d}(k)
	+ \beta^W \cdot \max\big(0, \bar{t}_d - t_k - \overline{TT}_{o,d}(k) \big)
	+ \gamma^W \cdot \max\big(0, t_k + \overline{TT}_{o,d}(k) - \bar{t}_d \big), \\
	\label{cost3}
	&P_{o,d}^{W,flex}(k) = \alpha^{W,flex} \cdot \overline{TT}_{o,d}(k)
	+ \beta^{W,flex} \cdot \max\big(0, \bar{t}_d^{flex} - t_k - \overline{TT}_{o,d}(k) \big)
	+ \gamma^{W,flex} \cdot \max\big(0, t_k + \overline{TT}_{o,d}(k) - \bar{t}_d^{flex} \big),
\end{align}

\noindent where $P_{o,s,b}^{S}(k), P_{o,d}^{W}(k)$, and $P_{o,d}^{W,flex}(k)\in\mathbb{R}^{+}$ are expressed in monetary units (\euro) with respect to Class~S or Class~W (with fixed and flexible work hours), when departing at discrete time-step $k$ toward school $s\in\mathcal{S}_b, b\in\mathcal{B}$ or destination $d\in\mathcal{D}$, respectively. Here, $\alpha^S, \beta^S, \gamma^S$ and $\alpha^W, \beta^W, \gamma^W, \alpha^{W,flex}, \beta^{W,flex}, \gamma^{W,flex}$ denote the shadow values (\euro/hour) of travel time, and this of early and late arrival for commuters of class S and W (with fixed and flexible work hours), respectively. 

The perceived cost functions defined in Eqs.~\eqref{cost1}–\eqref{cost3} depend on congestion through the OD-level travel times $\overline{TT}_{o,b}(k)$ and $\overline{TT}_{o,d}(k)$, which are in turn determined by the regional density trajectories 
\(\boldsymbol{\rho} = \{\rho_r(k)\}_{r \in \mathcal{R}, k \in \mathcal{K}}\), produced by the macroscopic traffic model described in Section~\ref{sec:tbyjhgt}, as well as on the imposed school and work start times, $\tilde{\tau}_s, s\in\mathcal{S}_b, b\in\mathcal{B}, \bar{t}_d, \bar{t}_d^{flex}, d\in\mathcal{D}$, respectively. Specifically, for an OD pair of commuters departing from origin \(o\) to destination $d$ at time-step \(k\), the expected travel time along a regional path $p$ is computed as the sum of travel times across the regions traversed by that path, accounting for the speed–density relationship of each region. For instance, the travel time through region \(r\) at time-step \(k\) is given by $\frac{L_{r,p}}{u_r(\rho_r(k))}$,
where \(u_r(\rho_r(k))\) is the average speed obtained from Eq.~\eqref{speed}. The OD-level travel time \(\overline{TT}_{o,d}(k)\) is then obtained by aggregating these regional travel times along the feasible path $p$ connecting \(o\in\mathcal{O}\) to \(d\in\mathcal{D}\) (or from \(o\) to school \(s\) located in region $b\in\mathcal{B}$). Consequently, any changes in the regional densities \(\rho_r(k)\), induced by school start time shifts or departure-time choices, directly affect the travel times and thus the perceived costs faced by commuters.\par 

In the Lower-Level problem, commuters choose their departure times by allocating demand across the discrete time steps $k\in\mathcal{K}$. These choices are represented by the equilibrium-based demand vector $\hat{\mathbf d}$. Based on this demand vector and the Upper-Level decision vector $\boldsymbol{\xi}$, the macroscopic traffic model determines the resulting evolution of regional traffic densities through the discrete-time mapping
\begin{equation}
	\label{mapping}
	\boldsymbol{\rho}(k+1)
	=
	\Upsilon\!\left(
	\boldsymbol{\rho}(k),
	\hat{\mathbf d}(k)
	\right),
	\quad \forall k\in\mathcal{K},
\end{equation}

\noindent where $\Upsilon(\cdot)$ captures the aggregate macroscopic traffic dynamics governing the propagation of congestion across regions, based on Eqs.~\eqref{eq:fixednew}–\eqref{arxontase}.

%Accordingly, the lower-level decision space is defined as the set of all non-negative demand vectors and density states such that
%\[
%\Xi_L
%=
%\Bigl\{
%\hat{\mathbf d}, \boldsymbol{\rho}\geq 0~\text{:}~ \text{the mapping \eqref{mapping} is satisfied}~ \forall k\in\mathcal{K}
%\Bigr\}.
%\]

\subsubsection{Formulation of the Deterministic Dynamic Multi-Class User Equilibrium}
\label{sec:formulation_ddue}

We consider a deterministic dynamic user equilibrium formulation with instantaneous information, consistent with Wardrop’s first principle extended to a time-dependent setting \citep{wardrop1952}. Commuters of each class are assumed to be perfectly rational and to observe the prevailing network conditions at each departure interval. However, they do not possess perfect forecasts of the future congestion evolution. All commuters sharing the same OD pair exhibit identical deterministic choice behavior.

Equilibrium is attained when no commuter can reduce their perceived travel cost by unilaterally shifting to another departure interval, given the prevailing network state. The formulation we adopt in this work is consistent with the \emph{instantaneous dynamic user equilibrium} (IDUE) paradigm used in the dynamic traffic assignment literature \citep{ma2018link}. Instantaneous equilibrium approaches assume that travelers base their decisions on prevailing traffic conditions rather than on fully anticipatory predictions of future congestion. Such formulations have been widely used in dynamic traffic assignment and perimeter control problems \citep{GUO202087, PAPAGEORGIOU1990471,  BAN2012360}. In large-scale, region-based networks, accurately forecasting future region-level congestion over the entire trip duration becomes particularly challenging due to nonlinear MFD-based dynamics and endogenous demand interactions. Predicting the full future evolution of regional densities would require anticipating the collective departure decisions of all commuters and the nonlinear propagation of congestion through the MFD dynamics. Extending the IDUE rationale to account for departure-time choice, therefore provides a behaviorally consistent and computationally tractable modeling assumption.\par

Given the induced discrete-time density mapping defined in Eq. \eqref{mapping} acting for a demand vector $\hat{\mathbf d}$, the perceived generalized cost operator is defined as
\[
\mathbf{P}(\hat{\mathbf d})
=
\mathbf{P}\big(
\hat{\mathbf d},
\boldsymbol{\rho}(\hat{\mathbf d})
\big),
\]

where $\mathbf{P}(\hat{\mathbf d})$ denotes the collection of perceived costs evaluated at a discrete time step $k$, i.e.,
\[
\mathbf{P}(\hat{\mathbf d})
=
\left\{
P_{o,s,b}^{S}(k),\;
P_{o,d}^{W}(k),\;
P_{o,d}^{W,\text{flex}}(k)
\right\}_{\forall (o,s,b),\,(o,d)}.
\]

The perceived costs per OD pair for the aforementioned classes of commuters are specified in 
Eqs.~\eqref{cost1}--\eqref{cost3}. Thus, the perceived costs depend endogenously on the departure times of commuters through the nonlinear density mapping (see Eq. \eqref{mapping}).

\medskip

\noindent
\paragraph{Feasible demand set}

\noindent Let parameters 
$D^{S}_{o,s}$,
$D^{W}_{o,d}$,
and
$D^{W,flex}_{o,d}$
denote the fixed total travel demand over the entire planning horizon. These quantities are exogenous and represent aggregate daily travel volumes. Specifically,

\begin{itemize}
	\item $D^{S}_{o,s}$ denotes the total number of school-related trips originating from origin $o\in\mathcal{O}$ and destined to school schedule $s\in\mathcal{S}_b$, where the school is located in region $b\in\mathcal{B}$.
	
	\item $D^{W}_{o,d}$ denotes the total number of work trips from origin $o$ to destination $d\in\mathcal{D}$ for commuters with fixed work start times.
	
	\item $D^{W,flex}_{o,d}$ denotes the total number of work trips from origin $o$ to destination $d\in\mathcal{D}$ for commuters with flexible work start times.
\end{itemize}

These volumes of vehicles are conserved across different departure times; the Lower-Level equilibrium redistributes demand across time but does not alter total volumes. Demand conservation is therefore enforced for each class of commuters:
\begin{align}
	\label{conservation1}
	& \sum_{k\in\mathcal{K}}
	\hat d^{S}_{o,s}(k)
	=
	D^{S}_{o,s},
	&& \forall o\in\mathcal{O},\;
	\forall b\in\mathcal{B},\;
	\forall s\in\mathcal{S}_b,
	\\[4pt]
	\label{conservation2}
	& \sum_{k\in\mathcal{K}}
	\hat d^{W}_{o,d}(k)
	=
	D^{W}_{o,d},
	&& \forall o\in\mathcal{O},\;
	\forall d\in\mathcal{D},
	\\[4pt]
	\label{conservation3}
	& \sum_{k\in\mathcal{K}}
	\hat d^{W,flex}_{o,d}(k)
	=
	D^{W,flex}_{o,d},
	&& \forall o\in\mathcal{O},\;
	\forall d\in\mathcal{D}.
\end{align}

\noindent For an Upper-Level decision vector $\boldsymbol{\xi}\in\Xi_U$, let the feasible set of demand profiles for the Lower-Level problem be defined as
\[
\Xi_L
=
\left\{
\hat{\mathbf d} \ge 0 :
\text{demand conservation constraints \eqref{conservation1} - \eqref{conservation3} hold}
\right\}.
\]

\noindent The above behavioral assumptions can be formally expressed as a deterministic dynamic user equilibrium condition over departure-time choices. This condition can be equivalently formulated as a finite-dimensional variational inequality (VI) problem over the feasible set of demand vectors \citep{Friesz1993,friesz2011,han2013}.

%\noindent Formally, for each OD pair and commuter class, the DUE can be expressed via the following condition
%
%\begin{equation}
%	\label{equilibrium_condition}
%	\hat d^{c,*}_{o,z}(k) > 0
%	\;\Rightarrow\;
%	P^{c}_{o,z}\big(k;\hat{\mathbf d}^*\big)
%	=
%	\min_{k'\in\mathcal{K}}
%	P^{c}_{o,z}\big(k';\hat{\mathbf d}^*\big), \forall o\in\mathcal{O}
%\end{equation}
%
%\noindent where $c\in\{S,W,\{W,\mathrm{flex}\}\}$ denotes the commuter class and $z\in\{s,d\}$ the index for school $s$ and destination $d$, respectively. The Wardrop equilibrium condition shown in \eqref{equilibrium_condition} states that positive demand values
%occur only at departure times that minimize perceived generalized cost.
%This complementarity-type condition can be equivalently expressed as a
%finite-dimensional variational inequality (VI) problem over the feasible set of demand vectors \citep{Friesz1993,friesz2011,han2013}.

\paragraph{Variational inequality formulation}

The deterministic dynamic user equilibrium consists in finding 
$\mathbf d^* \in \Xi_L$ for a fixed Upper-Level decision vector $\boldsymbol{\xi}$ such that
\begin{align}
	\label{eq:VI_DUE_structured}
	& \sum_{o\in\mathcal{O}}
	\sum_{b\in\mathcal{B}}
	\sum_{s\in\mathcal{S}_b}
	\sum_{k\in\mathcal{K}}
	P^{S}_{o,s,b}\big(k;\boldsymbol{\rho}(\mathbf d^*)\big)
	\Big(
	\hat d^{S}_{o,s}(k)
	-
	d^{S,*}_{o,s}(k)
	\Big)
	\nonumber \\
	& +
	\sum_{o\in\mathcal{O}}
	\sum_{d\in\mathcal{D}}
	\sum_{k\in\mathcal{K}}
	P^{W}_{o,d}\big(k;\boldsymbol{\rho}(\mathbf d^*)\big)
	\Big(
	\hat d^{W}_{o,d}(k)
	-
	d^{W,*}_{o,d}(k)
	\Big)
	\nonumber \\
	& +
	\sum_{o\in\mathcal{O}}
	\sum_{d\in\mathcal{D}}
	\sum_{k\in\mathcal{K}}
	P^{W,flex}_{o,d}
	\big(k;\boldsymbol{\rho}(\mathbf d^*)\big)
	\Big(
	\hat d^{W,flex}_{o,d}(k)
	-
	d^{W,flex,*}_{o,d}(k)
	\Big)
	\ge 0,
	\quad
	\forall \hat{\mathbf d}\in\Xi_L,
\end{align}

\noindent where the traffic density evolution $\boldsymbol{\rho}$ is captured through the mapping shown in Eq. \eqref{mapping}. Let $\mathbf{F}: \Xi_L \rightarrow \mathbb{R}^{n}$ be the mapping whose components are the perceived costs $P^S_{o,s,b}(k;\boldsymbol{\rho}(\mathbf{d}))$, $P^W_{o,d}(k;\boldsymbol{\rho}(\mathbf{d}))$, and $P^{W,\mathrm{flex}}_{o,d}(k;\boldsymbol{\rho}(\mathbf{d}))$, indexed over all classes, OD pairs, and time steps $k \in \mathcal{K}$, where
\begin{equation*}
	n := |\mathcal{O}|\left(\sum_{b\in\mathcal{B}}|\mathcal{S}_b| + 2|\mathcal{D}|\right)|\mathcal{K}|.
\end{equation*} 

\noindent Then the variational inequality~\eqref{eq:VI_DUE_structured} can be compactly written as
\begin{equation*}
	\big\langle 
	\mathbf{F}(\mathbf d^*),
	\hat{\mathbf d} - \mathbf d^*
	\big\rangle
	\ge 0,
	\quad \forall \hat{\mathbf d} \in \Xi_L,
\end{equation*}
where $\langle \cdot,\cdot \rangle$ denotes the standard inner product.\par

The variational inequality implies that any feasible deviation from the equilibrium departure demand vector $\mathbf d^*$ cannot reduce the total perceived travel cost. Hence, no reallocation of departures across time intervals yields a lower aggregate cost, which is consistent with Wardrop’s equilibrium principle.\par

\subsection{Bilevel Optimization Problem}
\label{sec:bilevel_problem}

The overall problem exhibits a leader-follower structure 
\citep{RAMEZANI20151}. In the Upper-Level problem, the school start times are selected, while explicitly anticipating the equilibrium response of commuters. In the Lower-Level problem, the deterministic dynamic user equilibrium is defined by the variational inequality~\eqref{eq:VI_DUE_structured}. Formally, the bilevel optimization problem can be expressed as
\begin{subequations}
	\label{eq:bilevel_VI}
	\begin{align}
		\left(\BilevelOpt\right)\quad
		\underset{\boldsymbol{\xi}\in\Xi_U}{\text{Minimize}}~
		& \displaystyle \Big\{ J_{\mathrm{TTS}}(\boldsymbol{\xi}, \hat{\mathbf d}), 
		J_{\mathrm{STC}}(\boldsymbol{\xi}) \Big\} 
		\\[1mm]
		\text{s.t.} \quad
		& \text{Constraints from Upper-Level problem: } 
		\eqref{mfd} - \eqref{QUO1}, ~
		\eqref{bin} - \eqref{binary}, 
		\nonumber
		\\[1mm]
		& \text{(DUE-VI)} \quad
		\big\langle 
		\mathbf{F}(\mathbf d^*; \boldsymbol{\xi}),
		\hat{\mathbf d} - \mathbf d^*
		\big\rangle
		\ge 0,
		\quad \forall \hat{\mathbf d} \in \Xi_L.
		\label{eq:lower_level_VI}
	\end{align}
\end{subequations}

\noindent It is worth reiterating that

\begin{itemize}
	\item $J_{\mathrm{TTS}}(\boldsymbol{\xi}, \hat{\mathbf d})$ is the TTS value, which depends on the Upper-Level decision vector $\boldsymbol{\xi} \in \Xi_U$ and the Lower-Level demand vector $\hat{\mathbf d}\in\Xi_L$.  
	\item $J_{\mathrm{STC}}(\boldsymbol{\xi})$ is the STC value.  
	\item $\Xi_L$ is the feasible set of demand profiles for all commuter classes, satisfying demand conservation constraints for each OD pair (Eqs.~\eqref{conservation1}-\eqref{conservation3} above).  
	\item $\mathbf{F}(\mathbf d^*; \boldsymbol{\xi})$ is the generalized cost mapping, which returns the vector of perceived costs for all commuters given the demand vector $\mathbf d^*$ obtained from the Lower-Level and the Upper-Level decision vector $\boldsymbol{\xi}$. 
\end{itemize}

\medskip

In its current form, Problem~$\BilevelOpt$ constitutes a Bi-Objective \emph{Mathematical Program with Equilibrium Constraints} (MPEC), where the equilibrium constraint is represented by the Lower-Level problem variational inequality~\eqref{eq:lower_level_VI}. In the proposed solution methodology (see Section \ref{sec:solution}), the two objectives are handled together using the $\epsilon$-constraint scalarization technique \citep{Mavrotas2009}. For completeness and consistency between the problem statement and the solution methodology, we next present the equivalent single-objective bilevel reformulation.

\subsection{Single-objective $\epsilon$-constraint bilevel reformulation}
\label{sec:bilevel_eps_constraint}

From a planning perspective, large changes in school start times are not considered acceptable in practice. Consequently, we assume that the set of admissible solutions is implicitly restricted by a threshold $\epsilon$ (measured in minutes), representing the overall permissible school start time change. In other words, only school scheduling policies that satisfy the condition $J_{\mathrm{STC}}(\boldsymbol{\xi}) \leq \epsilon$ are regarded as 
acceptable. \par

To handle the bi-objective nature of Problem~$\BilevelOpt$, we adopt the $\epsilon$-constraint technique, in which the STC metric is converted into a hard feasibility constraint, while the TTS metric remains the single optimization objective. Then for a fixed value of $\epsilon$, the planning problem can be reformulated as the following bilevel MPEC:
\begin{subequations}
	\label{eq:bilevel_eps}
	\begin{align}
		(\BilevelOptwithepsilon)\quad
		\underset{\boldsymbol{\xi}\in\Xi_U}{\text{Minimize}}~~
		& J_{\mathrm{TTS}}(\boldsymbol{\xi}, \hat{\mathbf d})
		\\[1mm]
		\text{s.t.}\quad
		& J_{\mathrm{STC}}(\boldsymbol{\xi}) \le \epsilon,
		\label{eq:STC_constraint_bilevel}
		\\
		& \text{Constraints from Upper-Level problem: } 
		\eqref{mfd} - \eqref{QUO1},~
		\eqref{bin} - \eqref{binary},
		\nonumber
		\\[1mm]
		& \text{(DUE-VI)} \quad
		\big\langle 
		\mathbf{F}(\mathbf d^*; \boldsymbol{\xi}),
		\hat{\mathbf d} - \mathbf d^*
		\big\rangle
		\ge 0,
		\quad \forall \hat{\mathbf d} \in \Xi_L.\nonumber
		\label{eq:lower_level_VI_eps}
	\end{align}
\end{subequations}

\noindent Problem~$\BilevelOptwithepsilon$ represents the practical planning problem considered in this paper. For a given $\epsilon$, the decision-maker seeks the school start time vector $\boldsymbol{\xi}$ that minimizes network congestion, expressed through the TTS metric, accounting at the same time for the dynamic user equilibrium conditions. 

Note that the constraint $J_{\mathrm{STC}}(\boldsymbol{\xi}) \le \epsilon$ is imposed at the level of the bilevel problem, but it only involves the Upper-Level decision variables $\boldsymbol{\xi}$ (school start times). Consequently, after the bilevel problem is reformulated into a single-level optimization problem, the $\epsilon$-constraint naturally appears as a constraint of the Upper-Level optimization problem.

\section{Solution Approach}
\label{sec:solution}

\begin{figure}[t]
	\centering
	\includegraphics[width=15cm]{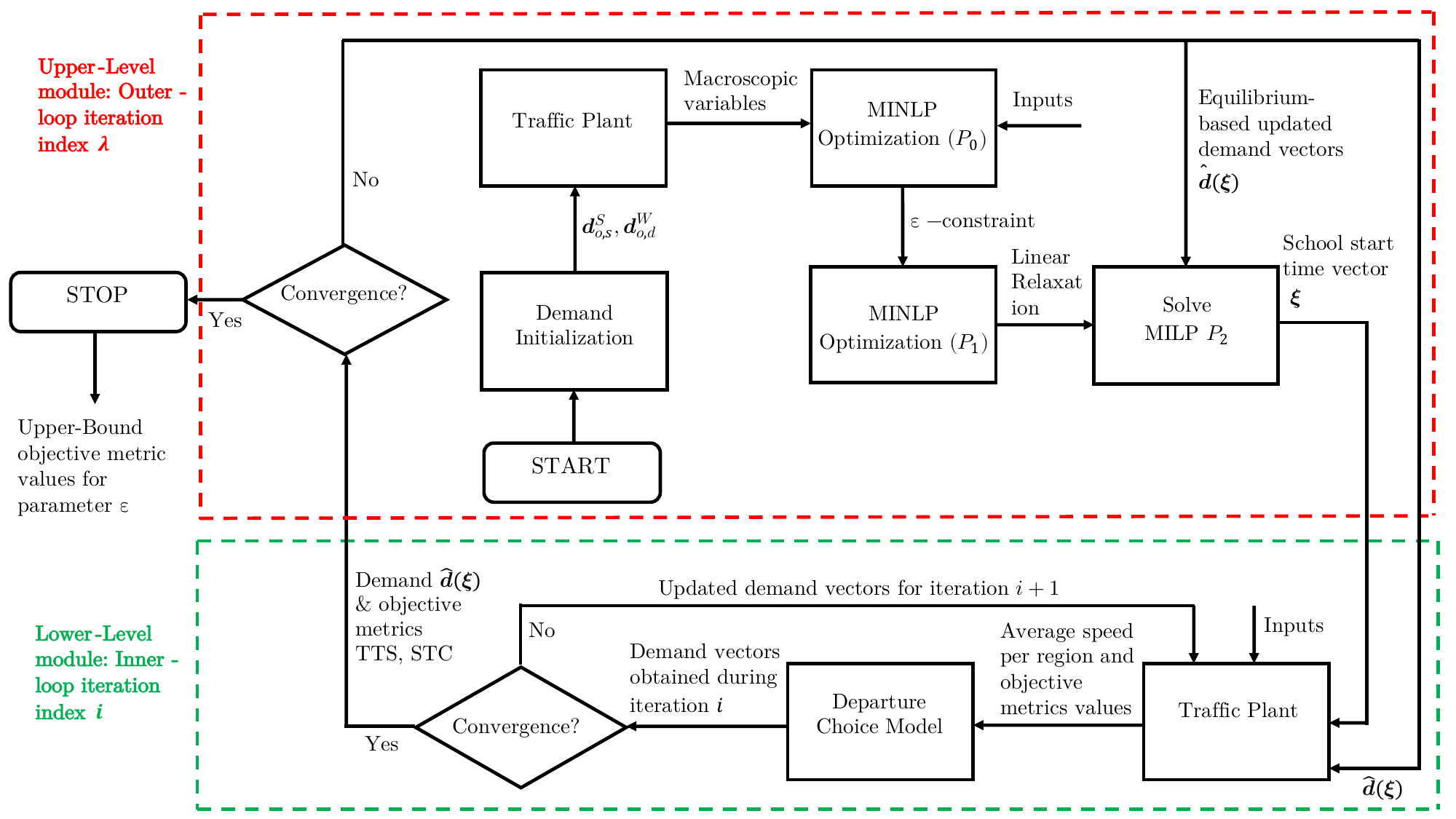} 
	\caption{Flowchart of the proposed solution approach for school start time selection for a given value of the overall permissible school start time change $\epsilon$. The flowchart consists of an Upper-Level module, which determines school start time decisions, and a Lower-Level module, which captures equilibrium-based demand responses under the derived school start times. The interaction between the two modules enables the evaluation of network-wide impacts of alternative school schedule configurations.}
	\label{fig:flowchart}
\end{figure}

One common approach to deal with bilevel optimization problems is to reformulate them as single-level optimization problems through Karush--Kuhn--Tucker (KKT) conditions \citep{dempe2002foundations, bard2013practical}. However, even under such a reformulation, the resulting problem remains computationally intractable in our setting due to the nonlinear and nonconvex macroscopic traffic dynamics induced by the MFD-based regional model.

To overcome this difficulty, we adopt an alternating solution framework that decomposes the bilevel problem \BilevelOptwithepsilon\ into two interacting modules, as illustrated in Fig.~\ref{fig:flowchart}. The Upper-Level module (red dashed box) determines the school start time decision vector $\boldsymbol{\xi}$ for a given demand vector $\hat{\mathbf{d}}$, while the Lower-Level module (green dashed box) computes equilibrium demand responses under the selected schedule, capturing commuters' departure time adjustments. Within Fig.~\ref{fig:flowchart}, each module is depicted as a structured pipeline comprising multiple building blocks, traffic state evaluation through the MFD-based traffic plant, demand updates, and convergence checks. Accordingly, the Upper-Level and Lower-Level modules should not be interpreted as being identical to the optimization problem $P_0$ and the Variational Inequality problem \eqref{eq:VI_DUE_structured}, respectively, but rather as higher-level algorithmic constructs that embed these problems within an iterative solution procedure. Both modules share common inputs, including school and workplace start times and regional macroscopic parameters, while the Lower-Level module additionally requires behavioral parameters related to earliness, lateness, and travel-time costs. Next, we provide a high-level description of the Upper-Level and Lower-Level modules before delving into the details of each module.

The Upper-Level module operates in two stages. First, the bi-objective problem $P_0$ is reformulated into a single-objective formulation through the introduction of a hard constraint, yielding problem $P_1$, while the demand is initialized under free-flow conditions (Section~\ref{sec:demand_initialization}). Second, problem $P_1$ is relaxed to obtain a mixed-integer linear programming problem, referred to as $P_2$, which provides a lower-bound to the original MINLP $P_1$. When the obtained lower-bound also satisfies the original MINLP constraints, then this solution is deemed feasible. This lower-bound contains the school start time vector $\boldsymbol{\xi}$ together with the lower-bound values associated with the TTS and STC metrics, respectively. 

The resulting vector $\boldsymbol{\xi}$ is subsequently passed as input to the Lower-Level module. This module is based on the variational inequality \eqref{eq:VI_DUE_structured}, which characterizes the user-equilibrium conditions governing commuters’ departure time choices. Due to the time-dependent nature of macroscopic traffic dynamics, the interaction of multiple classes of commuters, and the network-wide coupling induced by the macroscopic traffic model, solving the Variational Inequality \eqref{eq:VI_DUE_structured} exactly is computationally challenging. As a result, most studies in the literature rely on numerical, simulation-based, or heuristic procedures to approximate equilibrium conditions rather than deriving optimal solutions \citep{Han2014,GUO202087}. Hence, in the present work, the Lower-Level module employs a numerical scheme that aims to approximate the deterministic dynamic user equilibrium. Specifically, the Traffic Plant is first simulated to obtain region-level traffic states, including average speeds and the corresponding values of the objective metrics, based on the demand initialization performed in the Upper-Level module. These outputs are then provided to the departure-time choice model (defined in Section \ref{sec:depart_complete}), which produces updated demand vectors for each class of commuters. This procedure is performed iteratively until approximate user-equilibrium conditions are reached, yielding the demand vector $\hat{\mathbf d}(\boldsymbol{\xi})$.

The interaction between the two modules follows a sequential and iterative information exchange. At each outer-loop iteration, the Upper-Level module produces the school start time vector $\boldsymbol{\xi}\in\Xi_U$, while the Lower-Level module returns the corresponding equilibrium demand response $\hat{\mathbf d}(\boldsymbol{\xi})$. The updated demand vectors are then fed back to the Upper-Level module (see red dashed line box in Fig. \ref{fig:flowchart}), and the same procedure is pursued until convergence, yielding a feasible upper-bound to the Upper-Level problem $P_1$. Finally, the entire procedure is repeated for multiple values of the permissible overall school start time change $\epsilon$ to construct the Pareto front describing the trade-off between TTS and STC.

\subsection{Demand Initialization}
\label{sec:demand_initialization}

We initialize the demand vectors for each class of commuters, i.e., $d_{o,s}^S(k)$ and $d_{o,d}^W(k)$, $\forall k \in \mathcal{K}$, under free-flow conditions. This procedure yields the baseline demand input to the Upper-Level module.

For \emph{school commuters} (class $S$), the demand profile is constructed based on the initial school start time $\tau_s$ and the shortest-path travel time from origin region $o \in \mathcal{O}$ to destination region $b \in \mathcal{B}$, where school $s \in \mathcal{S}_b$ is located. In particular, under free-flow conditions, the preferred departure time from origin $o$ is taken to be centered around $\tau_s - t_{o,b}^{\text{SP}}$, where $t_{o,b}^{\text{SP}}$ denotes the shortest-path travel time from $o$ to $b$.

The initialized demand vector is therefore constructed so that the majority of departures are concentrated in a neighborhood around the time spectrum $\tau_s - t_{o,b}^{\text{SP}}$, while allowing for temporal dispersion, i.e., departures may also occur in earlier or later discrete time steps $k \in \mathcal{K}$.

An analogous initialization procedure is applied to commuters of class $W$, for both fixed and flexible work schedules.

\subsection{Lower-Level Module}
\label{sec:lower_level_controller}

The Lower-Level module (green dashed line box in Fig.~\ref{fig:flowchart}) represents the inner-loop component of the proposed solution approach and is responsible for modeling commuters’ behavioral responses to a given school start time vector $\boldsymbol{\xi}$. Its role is to approximate deterministic dynamic user equilibrium conditions in the network.\par

\paragraph{Relation to the Variational Inequality formulation}

The procedure described in this section is designed to compute an approximate solution of the variational inequality \eqref{eq:VI_DUE_structured}, which formalizes the deterministic dynamic user equilibrium with departure-time
choice. %Recall that the VI expresses Wardrop’s equilibrium principle in a time-dependent setting: for each OD pair and commuter class, positive demand is assigned only to departure times that minimize the perceived generalized cost. This equilibrium condition can be equivalently expressed as a fixed-point problem \citep{sheffi1985urban}. \par
%Let $H(\mathbf{\hat d})$ denote the departure-time update 
%mapping that given the current network conditions induced by a demand vector $\mathbf{\hat d}$, reallocates demand toward cost-minimizing departure-time intervals. At equilibrium, the demand 
%vector must satisfy the fixed-point condition
%\[
%\mathbf d^* = H(\mathbf d^*),
%\]
%meaning that no OD pair of commuters can further reduce their perceived cost by shifting their departure time.

%The algorithm presented next is constructed precisely to approximate this fixed point. 
Starting from an initial demand vector derived under free-flow conditions, a general methodology for solving the variational inequality \eqref{eq:VI_DUE_structured} proceeds iteratively through three main steps:

\begin{enumerate}[leftmargin=1.5cm]
	\item[\textbf{Step 1}:] Simulate network conditions under the current demand.
	\item[\textbf{Step 2}:] Compute best-response departure-time adjustments that reduce perceived
	costs over a range of candidate departure times.
	\item[\textbf{Step 3}:] Update the demand vector toward this best-response solution using the Method of Successive Averages (MSA).
\end{enumerate}

This general methodology is widely used for solving dynamic traffic assignment problems and for computing equilibria, formulated as variational inequalities \citep{sheffi1985urban,friesz2011,Yildirimoglu2021}. When the iterative process converges, the resulting demand vector satisfies the Wardrop's equilibrium condition and therefore constitutes an approximate solution of the variational inequality \eqref{eq:VI_DUE_structured}.

Building on this three-step methodology, Algorithm~\ref{alg:user_equilibrium}, which constitutes the numerical core of the Lower-Level module, computes an approximate solution to the dynamic multi-class user equilibrium problem, equivalently expressed by the VI \eqref{eq:VI_DUE_structured}. The algorithm takes as input the traffic network parameters, external demand profiles, departure-time choice parameters, and simulation parameters including the permissible overall school start time change $\epsilon$, the initial school start time vector $\boldmath{\tau}^{\text{INIT}}$$(\boldsymbol{\xi}, \epsilon)$ (which is a function of the Upper-Level decision vector $\boldsymbol{\xi}$), the work start time vectors $\boldsymbol{\bar{t}}_d$, $\boldsymbol{\bar{t}}_d^{flex}$, the convergence tolerance $\delta_1$, and the search parameters $\varPhi$ and $\Delta\phi$ (Line~1). The algorithm begins by initializing the iteration counter $i=0$, the convergence indicator $\text{Gap}$, and the initial demand vectors $\mathbf{\hat{d}}^{i=0}$ based on free-flow departure times (Line~2). Then Algorithm \ref{alg:user_equilibrium} presents a tailored implementation of the general three-step methodology for approximately solving the VI problem. In particular:

\begin{algorithm}[t]
	\algsetup{linenosize=\tiny}
	\footnotesize
	\begin{algorithmic}[1]
		\STATE \textbf{Input:}
		\begin{itemize}[leftmargin=*]
			\item Traffic Network Parameters: $u_r^f, \rho_r^C, \rho_r^J, q_r^C, g_r^C, L_r, l_r, r\in\mathcal{R}, C_{r,j}^{\text{MAX}}, \alpha_{r,j}, r\in\mathcal{R}, j\in\mathcal{J}_r^-$, $\theta_{r,j,d}(k), r\in\mathcal{R},  j\in\mathcal{J}_r^-,  k\in\mathcal{K}, d\in\mathcal{D}, \theta_{r,j,b}(k), r\in\mathcal{R},  j\in\mathcal{J}_r^-,  k\in\mathcal{K}, b\in\mathcal{B}$.
			\item External demand: $d_{o,s}^S(k),  \forall o\in\mathcal{O}, \forall b\in\mathcal{B}, \forall s\in\mathcal{S}_b$, $\forall k\in\mathcal{K}, d_{o,d}^W(k), \forall o\in\mathcal{O}, \forall d\in\mathcal{D}, \forall k\in\mathcal{K}, d_{o,s,m}^S(k), \forall o\in\mathcal{O},  \forall b\in\mathcal{B}, \forall s\in\mathcal{S}_b, \forall m\in\mathcal{M}, \forall k\in\mathcal{K}$.
			\item Departure-time choice parameters: $\alpha, \beta, \gamma$ for school, fixed work, flexible work
			\item Simulation parameters: $\boldsymbol{\xi}, \epsilon, \delta_1, K$, $\boldmath{\tau}^{\text{INIT}}(\boldsymbol{\xi}$$, \epsilon)$, $\boldsymbol{\bar{t}}_d, \boldsymbol{\bar{t}}_d^{flex}, \varPhi, \Delta\phi$
		\end{itemize}
		\STATE \textbf{Initialization:} $i=0$, $\text{Gap}=10000$, initialize $\mathbf{\hat{d}}^{i=0}$ using free-flow departure times.
		
		\WHILE{$\text{Gap}_i > \delta_1$}
		
		\STATE Simulate the traffic system dictated by Eqs. \eqref{mfd} - \eqref{QUO1} providing as input the school start time vector $\boldmath{\tau}^{\text{INIT}}(\boldsymbol{\xi}$$, \epsilon)$ and the work start time vectors $\boldsymbol{\bar{t}}_d, \boldsymbol{\bar{t}}_d^{flex}$, respectively, and compute average speeds $\mathbf{u}_r^i = \text{col}(u_r^i(\rho_r(k)))$ from Eq. \eqref{speed} and performance metrics $J_{\text{TTS}}^i, J_{\text{STC}}^i$ utilizing Eqs. \eqref{fixed} and \eqref{QUO1}, respectively for iteration $i$.
		
		\FOR{each commuter class $c \in \{S,W,W^{flex}\}$}
		
		\FOR{each region $r$ and destination $z\in\{s,d\}$}
		
		\STATE Compute perceived cost $P^c$ for the current departure time using Eqs.~\eqref{cost1}--\eqref{cost3}.
		
		\STATE Generate candidate departure times within $[-\varPhi,+\varPhi]$ minutes with increment $\Delta\phi$ and evaluate their perceived costs.
		
		\STATE Select the candidate departure time $\chi_{r,z}^{c,*}$ that minimizes the perceived cost.
		
		\STATE Update the departure time using the MSA scheme: $\chi_{r,z}^{c,i,\text{MSA}}=(1-\mu)\chi_{r,z}^{c,i-1}+\mu \chi_{r,z}^{c,*},
		\qquad \mu=1/i.$
		
		\STATE Update the demand vector $\mathbf{\hat{d}}_{r,z}^{c,i}$ using Eq.~\eqref{shift_departure_1}.
		
		\STATE Set $\chi_{r,z}^{c,i}=\chi_{r,z}^{c,i,\text{MSA}}$.
		
		\ENDFOR
		\ENDFOR
		
		\IF{$i>0$}
		\STATE $\text{Gap}_i = \big\|\mathbf{\hat{d}}^{\,i} - \mathbf{\hat{d}}^{\,i-1}\big\|$.
		\ENDIF
		
		\STATE $i=i+1$
		
		\ENDWHILE
		
		\STATE \textbf{Output:} User-equilibrium demand vectors $\mathbf{\hat{d}}^{\Iota}$ obtained in iteration $\Iota\equiv i$, and performance metrics $J_{\text{TTS}}^{\Iota}, J_{\text{STC}}^{\Iota}$.
		
	\end{algorithmic}
	\caption{\Approximate\ (\ApproximateAcronym) Algorithm}
	\label{alg:user_equilibrium}
\end{algorithm}

\begin{itemize}
	\item \textbf{Step 1 – Network simulation (Line 4):}
	At each outer-loop iteration, the algorithm simulates the traffic system using the current demand vectors and the given school and work start times (Line~4). This step produces the region-based speed profiles $\mathbf{u}_r^i$ and the associated system performance indicators required to evaluate commuters’ perceived costs.
	
	\item \textbf{Step 2 – Best-response departure-time update (Lines 5--12):}
	For each commuter class, origin region, and destination, the algorithm computes perceived costs (Line~7), evaluates candidate departure times within the search window (Line~8), selects the cost-minimizing candidate (Line~9), and applies the MSA to update departure times and the corresponding demand allocation (Lines~10--12).
	
	\item \textbf{Step 3 – Convergence check and iteration update (Lines 13--18):}
	The global change in the demand vector between successive iterations is computed to form the convergence indicator $\text{Gap}_i$ (Lines~15-17). The iteration counter is updated and the procedure repeats until $\text{Gap}_i \le \delta_1$.
\end{itemize}

We now provide a detailed mathematical description of the main components of Algorithm~\ref{alg:user_equilibrium}. 

\paragraph{Step 1: Network simulation}
At each iteration $i$, the traffic system is simulated using Eqs.~\eqref{mfd}--\eqref{QUO1} with the current demand vectors and the given school and work start times. This simulation yields the region-based average speed vectors $\mathbf{u}_r^i$, from which the system performance indicators $J_{\text{TTS}}^i$ and $J_{\text{STC}}^i$ can be obtained (Line~4). The vector $\mathbf{u}_r^i$ collects the average speeds at each discrete time step $k\in\mathcal{K}$ for region $r$, i.e.,
\[
\mathbf{u}_r^i = \text{col}(u_r^i(\rho_r(k))), \quad k \in \mathcal{K},
\]  
where $u_r^i(\rho_r(k))$ denotes the average speed obtained from Eq.~\eqref{speed} during iteration $i$.

\paragraph{Step 2: Best-response departure-time update}
The algorithm evaluates departure-time adjustments through an iterative procedure (Lines 5-11). For each class $c\in\{\text{School},\text{Work},\text{Work-flex}\}$ we process all origin regions $r$ and destinations $d$ (or schools $s$) associated with the corresponding OD pairs (Line~6). Within these loops, the perceived cost $P^c$ associated with the current departure time is computed using the cost functions defined in Eqs.~\eqref{cost1}--\eqref{cost3} (Line~7). Subsequently, candidate departure times are generated within the interval $[-\varPhi,+\varPhi]$ minutes around the current departure time using a discretization step $\Delta\phi$, and the corresponding perceived costs are evaluated for each candidate departure time (Line~8). Among the generated candidates, the departure time that minimizes the perceived cost is identified and denoted by $\chi_{r,z}^{c,*}$ (Line~9). To prevent oscillatory behavior caused by highly reactive user updates and non-monotonic congestion feedbacks, the Method of Successive Averages is applied \citep{UKKUSURI2009625,HAN201516}. Specifically, the updated departure time $\chi_{r,z}^{c,i,\text{MSA}}$ is computed as a linear combination of the previous departure time $\chi_{r,z}^{c,i-1}$ and the cost-minimizing candidate $\chi_{r,z}^{c,*}$ with step size $\mu=1/i$ (Line~10). The demand vector $\mathbf{\hat{d}}_{r,z}^{c,i}$ is then updated according to Eq.~\eqref{shift_departure_1}, where $\Delta \chi_{r,z}^c = |\chi_{r,z}^{c,i,\text{MSA}} - \chi_{r,z}^{c,i}|$ denotes the difference in time-steps between the current and the \textit{new} candidate departure time (Line 11). Next, the updated departure time is stored for use in the subsequent outer-loop iteration (Line 12).

\begin{equation}
	\label{shift_departure_1}
	\hat{d}_{r,z}^{c,i}(k)=
	\begin{cases}
		d_{r,z}^{c}(k-\Delta\chi_{r,z}^{c}), 
		& \text{if } \Delta\chi_{r,z}^{c}>0 \text{ and } \Delta\chi_{r,z}^{c}<k\le K, \\[6pt]
		
		0, 
		& \text{if } \Delta\chi_{r,z}^{c}>0 \text{ and } 1\le k\le\Delta\chi_{r,z}^{c}, \\[6pt]
		
		d_{r,z}^{c}(k), 
		& \text{if } \Delta\chi_{r,z}^{c}=0, \\[6pt]
		
		d_{r,z}^{c}(k-\Delta\chi_{r,z}^{c}), 
		& \text{if } \Delta\chi_{r,z}^{c}<0 \text{ and } 1\le k\le K+\Delta\chi_{r,z}^{c}, \\[6pt]
		
		0, 
		& \text{if } \Delta\chi_{r,z}^{c}<0 \text{ and } K+\Delta\chi_{r,z}^{c}<k\le K,
	\end{cases}
\end{equation}

\noindent where $z\in\{s,d\}$ denotes the destination index (school or workplace) and $c\in\{S,W,W^{flex}\}$ represents the commuter class.

\paragraph{Step 3: Stopping criterion based on the variational inequality}
Ideally, the convergence of the algorithm would be assessed using the residual of the variational inequality~\eqref{eq:VI_DUE_structured}. At
iteration $i$, this residual can be expressed as
\begin{equation*}
	\text{VI-Residual}_i = \max_{\hat{\mathbf d}^i \in \Xi_L} 
	\big\langle \mathbf{F}(\mathbf d^*), \hat{\mathbf d}^i - \mathbf d^* \big\rangle,
\end{equation*}

\noindent which represents the maximum possible improvement in perceived cost over all feasible deviations from the current demand vector $\mathbf{\hat{d}}^i$. However, since the true equilibrium $\mathbf d^*$ is unknown and the feasible demand set $\Xi_L$ can be very large, the exact residual cannot be computed directly.\par 

%To formalize the convergence metric, we define the \emph{global demand vector} at iteration $i$ as
%\begin{equation*}
%	\mathbf{\hat{d}}^i = \Big[ \mathbf{\hat{d}}_{r,z}^{c,i} \Big], 
%	\quad \forall r \in \mathcal{R}, \; \forall z \in \{s,d\}, \; \forall c \in \{S,W,W^{flex}\},
%\end{equation*}
%where $\mathbf{\hat{d}}_{r,z}^{c,i}$ collects the demand values for all discrete time steps $k \in \mathcal{K}$. In other words, $\mathbf{\hat{d}}^i$ concatenates the demand across all regions, destinations, commuter classes, and discrete time steps. 

A practical surrogate for the VI residual is the global change in demand between successive iterations: Formally, the convergence indicator is defined as
\begin{equation*}
	\text{Gap}_i = \big\|\mathbf{\hat{d}}^{\,i} - \mathbf{\hat{d}}^{\,i-1}\big\|,
\end{equation*}
where $\|\cdot\|$ denotes the Euclidean norm. The algorithm stops when $\text{Gap}_i \le \delta_1$, which ensures that no significant reallocation of departure times would further reduce the aggregate perceived cost of commuters, consistent with Wardrop’s equilibrium principle (Line 16). Next, the iteration counter is incremented (Line~18), and the algorithm returns to the beginning of the outer loop until the stopping condition $\text{Gap}_i\le\delta_1$ is satisfied (Line~3).\par

Upon convergence, the algorithm provides as output the approximate user-equilibrium demand vectors $\mathbf{\hat{d}}^{\Iota}$ obtained at the final outer-loop iteration $\Iota \equiv i$ of Algorithm \ref{alg:user_equilibrium}, together with the corresponding system performance indicators $J_{\text{TTS}}^{\Iota}$ and $J_{\text{STC}}^{\Iota}$ (Line~20).

\subsection{Upper-Level Module}
\label{sec:upper_level_controller}

The Upper-Level module corresponds to the red dashed box in Fig.~\ref{fig:flowchart} and represents the outer-loop component of the proposed solution approach. Its role is to determine the school start time decision vector $\boldsymbol{\xi}$, given an approximate equilibrium demand vector $\mathbf{\hat{d}}$. By systematically exploring alternative school start time configurations, the Upper-Level module evaluates their network-wide impacts through the induced changes in commuter demand patterns and resulting equilibrium conditions.

\subsubsection{Demand Projection Mechanism}
\label{sec:demand_projection}

\begin{figure}[t]
	\centering
	\includegraphics[width=10cm]{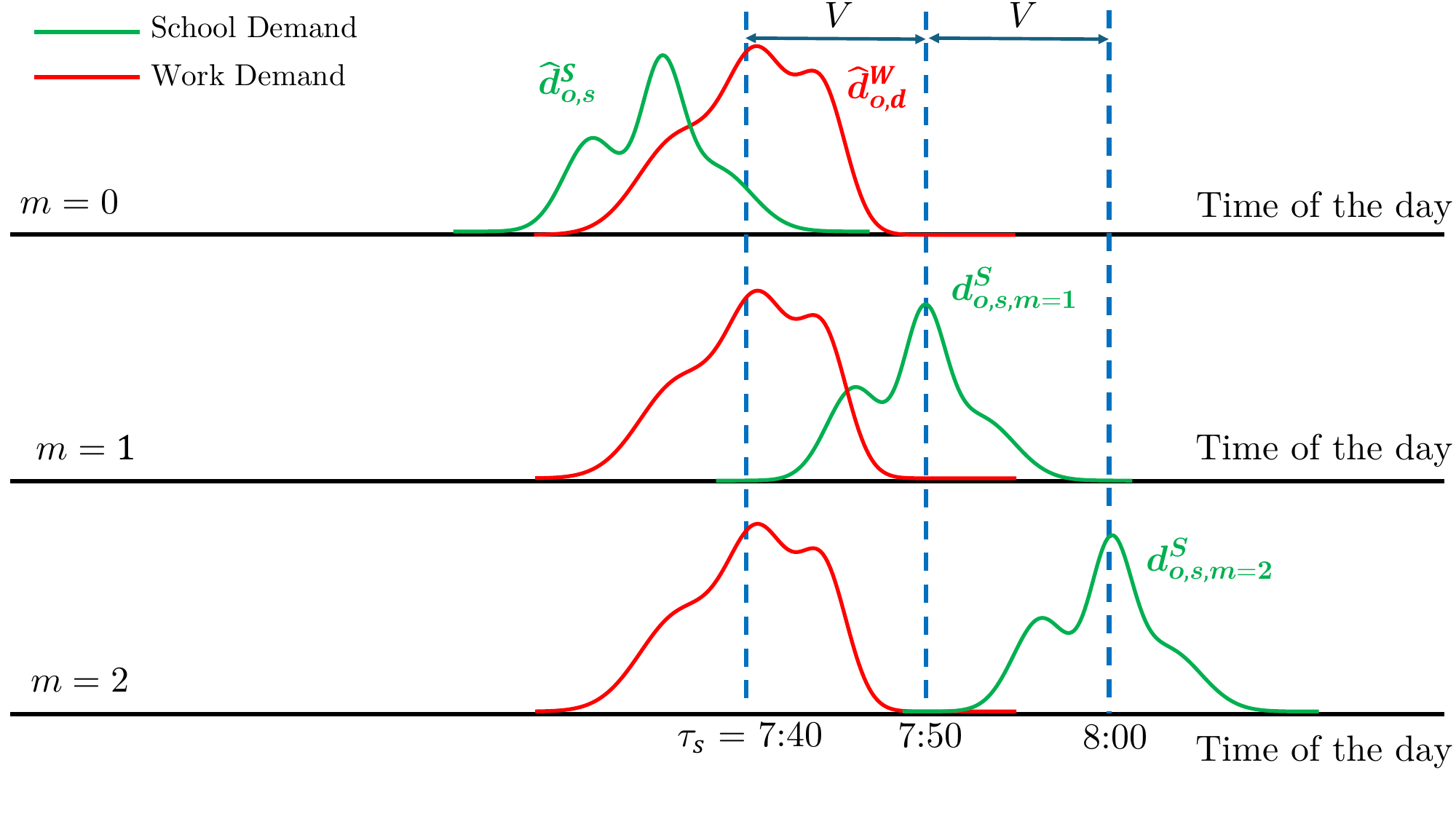} 
	\caption{Illustration of OD-level demand distributions for work-related commuters (class W) and school-related commuters (class S) under different realizations of the discrete shifting index $m$. Only one shifting index can be selected for each school. The case $m=0$ corresponds to the baseline demand distributions for each commuter class as obtained from Algorithm \ApproximateAcronym\ (incorporated within the Lower-Level module), i.e., $\hat{\mathbf{d}}_{o,s}^S\in\mathbb{R}^{K\times 1}, o\in\mathcal{O}, s\in\mathcal{S}_b, b\in\mathcal{B}$ for school-related demand and $\hat{\mathbf{d}}_{o,d}^W\in\mathbb{R}^{K\times 1}, o\in\mathcal{O}, d\in\mathcal{D}$ for work-related demand. The figure depicts the demand projection approximation step, in which school-related demand distributions are shifted according to the school start time decisions. While the illustration shows the forward shifting case ($m>0$), the proposed mathematical framework is equally applicable to backward shifts ($m<0$).}
	\label{fig:ggwheh_UE}
\end{figure}

Typically, a change in the Upper-Level school start time decision vector $\boldsymbol{\xi}$, implies a change in the demand distribution associated with commuters of class S. For small deviations of $\boldsymbol{\xi}$ it is reasonable to assume that the shape of school-related demand (class S) is not significantly altered; instead it is shifted earlier or later by approximately the same amount as the change in the school start time. In other words, the temporal profile of school-related travel demand is expected to remain largely unchanged. This assumption is also validated from previous studies in an MFD-based regional context, where small perturbations in departure times are translated into approximately uniform shifts of the departure distribution, rather than complete redistributions of trips \citep{LAMOTTE2018794}. Motivated by this fact, we introduce a demand projection mechanism that shifts the baseline equilibrium demand associated with commuters of class S, i.e., $\hat{\boldsymbol{d}}_{o,s}^{S}\in\mathbb{R}^{K\times 1}, o\in\mathcal{O}, s\in\mathcal{S}_b, b\in\mathcal{B}$ (derived from Algorithm \ref{alg:user_equilibrium}) forward or backward, according to the discrete shift in school start times encoded by $\boldsymbol{\xi}$, without changing the shape of demand over time.

This behaviour is assimilated within the Upper-Level optimization problem $P_0$. Let us remind the reader that the duration of each discrete shifting interval is equal to $V$ (min), where every interval comprises $h$ time-steps of duration $T$ such that $V = h T$. Let $d_{o,s,m}^{S}(k)$ (veh) denote the number of vehicles belonging to class S that enter the network at time-step $k+mh$ from origin $o\in\mathcal{O}$ towards school $s\in\mathcal{S}_b, b\in\mathcal{B}$. The corresponding demand vector for the entire morning commute period is defined as $\boldsymbol{d}_{o,s,m}^{S} =
\big[ d_{o,s,m}^{S}(1),\, d_{o,s,m}^{S}(2),\, \dots,\, d_{o,s,m}^{S}(K) \big]^{\top} \in \mathbb{R}^{K \times 1}$.\par

Formally, for each origin $o$, school $s$, time step $k$, and discrete shift index $m\in\mathcal{M}$, we define the shifted school demand as
\begin{equation}
	\label{shift_lhhg}
	d_{o,s}^{S}(k;\boldsymbol{\xi}) = \sum_{m\in\mathcal{M}} \xi_{m,s} \, d_{o,s,m}^{S}(k),
\end{equation}
where
\begin{equation}
	\label{mapping_new}
	d_{o,s,m}^{S}(k) =
	\begin{cases}
		\hat{d}^{S}_{o,s}(k - mh), & m>0,~ mh < k \le K,\\
		0, & m>0,~ 1 \le k \le mh,\\
		\hat{d}^{S}_{o,s}(k), & m=0,\\
		\hat{d}^{S}_{o,s}(k - mh), & m<0,~ 1 \le k \le K+mh,\\
		0, & m<0,~ K+mh < k \le K.
	\end{cases}
\end{equation}

\noindent Equations~\eqref{shift_lhhg}–\eqref{mapping_new} define how the school-related component of demand is projected under a given decision vector $\boldsymbol{\xi}$ by constructing and selecting temporally shifted versions of the baseline equilibrium demand $\hat{d}_{o,s}^{S}(k)$. In particular, Eq.~\eqref{mapping_new} constructs the demand distribution associated with a given discrete shift index $m$ by translating the baseline equilibrium demand $\hat{d}_{o,s}^{S}(k)$ obtained from Algorithm~\ref{alg:user_equilibrium}, while Eq.~\eqref{shift_lhhg} combines these shifted profiles through the binary decision variables $\xi_{m,s}$. However, the demand projection mechanism must explicitly specify how the demand distributions of all commuter classes entering the Upper-Level problem are affected by the decision vector $\boldsymbol{\xi}$. In the proposed framework, school start time decisions are assumed to affect only commuters of class S. Consequently, the demand associated with work-related commuters (class W) is assumed to remain independent of the discrete school start time shifts and therefore retains the equilibrium reading obtained from Algorithm~\ref{alg:user_equilibrium}. This is expressed as
\begin{equation}
	\label{projection_works}
	d_{o,d}^{W}(k;\boldsymbol{\xi}) \equiv \hat{d}_{o,d}^{W}(k),\quad k\in\mathcal{K}.
\end{equation}\par 

\noindent Taken together, Eqs.~\eqref{shift_lhhg}–\eqref{projection_works} jointly define the demand projection mechanism that maps the decision vector $\boldsymbol{\xi}$ to the complete time-dependent demand entering the Upper-Level optimization problem $P_0$. This projection mechanism constitutes a \emph{modeling approximation} introduced to represent the impact of discrete school start time changes within the Upper-Level optimization problem $P_0$.  This demand projection mechanism is illustrated in Fig.~\ref{fig:ggwheh_UE}.\par 

Because the proposed mechanism operates through temporal translations of the baseline school-related demand distributions, any modification of the school start time $\tilde{\tau}_s = \tau_s +\xi_{m,s}mV$ directly corresponds to a translation of the associated demand profile along the time axis. Therefore, controlling variations in school start times is equivalent to controlling how much the projected demand distributions are allowed to shift in time. This observation will be used later in Section \ref{sec:relaxation}.\par 

Overall, the demand terms $d_{o,s}^{S}(k;\boldsymbol{\xi})$ and $d_{o,d}^{W}(k;\boldsymbol{\xi})$, obtained through the demand projection mechanism and the solution of Algorithm~\ref{alg:user_equilibrium}, respectively, jointly constitute the demand input of the optimization problem $P_0$ and appear explicitly in Eq.~\eqref{input}.

\subsubsection{Reformulation of Upper-Level problem $P_0$ using the $\epsilon$-constraint technique}
\label{sec:reformulation}

Building upon the bilevel $\epsilon$-constraint formulation introduced in Section~\ref{sec:bilevel_eps_constraint}, we now derive the single-level optimization problem that is solved in practice. Recall that the $\epsilon$-constraint is part of the bilevel optimization problem $\BilevelOptwithepsilon$ and restricts the admissible school start times through the condition $J_{\text{STC}}(\boldsymbol{\xi}) \le \epsilon$. Since this constraint depends exclusively on the Upper-Level decision variables, it naturally appears as a constraint of the resulting Upper-Level optimization problem $P_0$ after the bilevel problem is reformulated. In other words, the $\epsilon$-constraint restricts the leader’s feasible region, not the follower’s equilibrium.

For a fixed value of $\epsilon$, the Upper-Level optimization problem is derived. The aim is to determine the shifted school start times that minimize the TTS metric while ensuring that the STC metric does not exceed a user-specified threshold~$\epsilon$.  The reformulated optimization problem is therefore expressed as
\begin{subequations}
	\label{rfwga}
	\begin{align}
		\label{trb}
		(P_1) \quad \underset{\xi_{m,s}, \forall m,s}{\text{Minimize}} &~ \displaystyle J_{\textrm{TTS}} = T\sum_{k\in\mathcal{K}}\Big(S_a(k) - S_b(k)\Big) \\
		\textrm{Subject To:}&~~\textrm{Demand Dynamics}~ \eqref{shift_lhhg} - \eqref{projection_works}, \nonumber \\
		&~~\textrm{Traffic Dynamics}~ \eqref{mfd} -  \eqref{QUO1}, \nonumber \\
		&~~\text{Constraints}~ \eqref{bin}, \eqref{jam}, \eqref{binary}, \nonumber \\
		\label{epsi}
		&~~  J_{\text{STC}} \leq  \epsilon,\\
		\textrm{Initialization:}&~~\eqref{init},\nonumber\\
		\textrm{Input:}&~~\eqref{input}\nonumber.
	\end{align} 
\end{subequations}

\noindent The resulting optimization problem $P_1$ is a Mixed-Integer Nonlinear Program (MINLP). The set of constraints labeled \textit{Demand Dynamics} (Eqs.~\eqref{shift_lhhg}–\eqref{projection_works}) correspond exactly to the demand projection mechanism introduced in Section~\ref{sec:demand_projection}, which maps the school start time decision vector $\boldsymbol{\xi}$ to the demand entering the network. 

\subsubsection{Linear Relaxation: Computation of lower-bound to problem $P_1$}
\label{sec:relaxation}

Problem $P_1$ poses a non-convex combinatorial optimization problem, since in addition to having binary variables $\xi_{m,s}$, the constraints \eqref{upoptos}, \eqref{inssss}, \eqref{naruto}, \eqref{secondtt} and \eqref{secondtt1} that are included in $P_1$ are nonlinear and non-convex. This section provides a tractable reformulation of $P_1$ that relaxes the nonlinear and non-convex constraints incorporated in problem $P_1$ with convex approximation constraints. To accomplish that, each non-convex constraint of the optimization problem $P_1$ is relaxed with linear constraints that lie in convex domains that are supersets of the corresponding non-convex domain of the original optimization problem $P_1$. 
As an outcome, the obtained solution from this formulation yields lower bounds to the original MINLP problem $P_1$. \par 
For ease of notation in the relaxation procedure pursued in this section, we omit the dependency of the outflow and MFD functions on the density variable $\rho_r(k)$. Specifically, we use $q_r(k)$ instead of $q_r(\rho_r(k))$ and $g_r(k)$ instead of $g_r(\rho_r(k))$. The same notion is adopted for each quantity related to flow for each class of commuters. Next, we derive superset linear constraints for the five non-convex and nonlinear constraints of Problem $P_1$ mentioned before.\par
First, constraint \eqref{mfd} can be relaxed by substituting the equality sign ``$=$'' with the inequality sign ``$\leq$'' yielding
\begin{align}
	\label{traffic2}
	g_r(\rho_r(k)) \leq&~ c_{r,l}\rho_r(k) + b_{r,l}, \quad \rho_r(k)\in[\bar{\rho}_{r,l}, \bar{\rho}_{r,l+1}],~l=1,\ldots,N+1.
\end{align}

\begin{figure}[t]
	\centering
	\includegraphics[width=9cm]{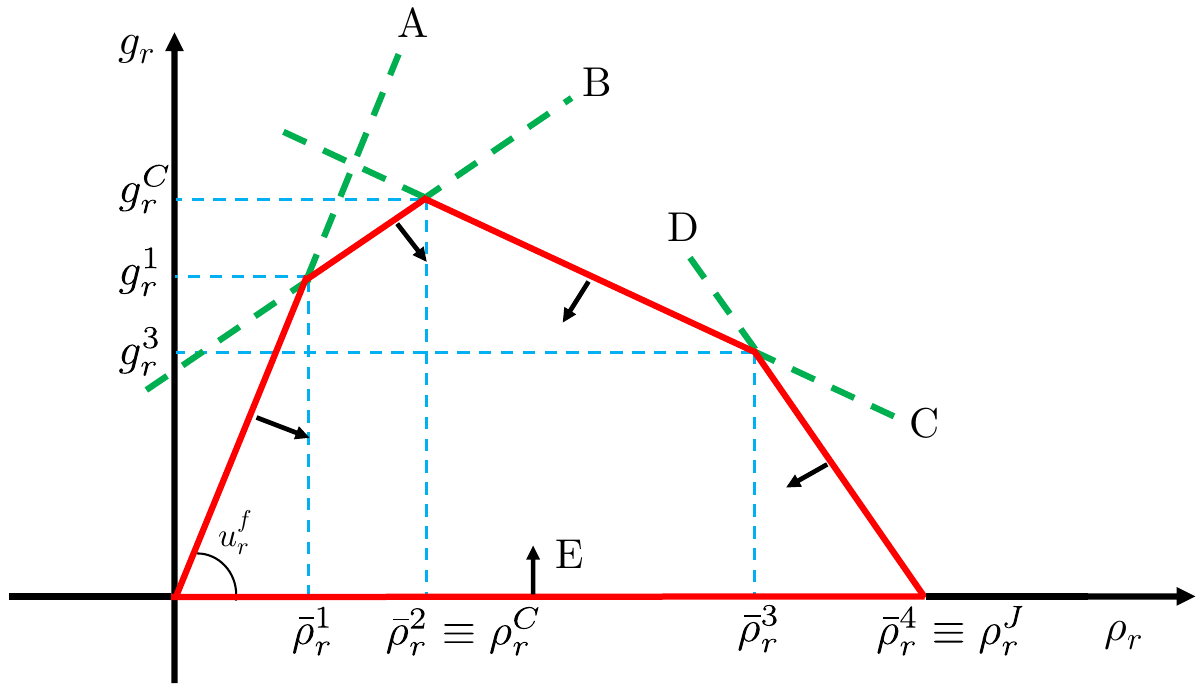} 
	\caption{An example of the feasibility domain for piecewise linear MFDs with $N=4$ linear segments: i) the exact constraint \eqref{mfd} is captured with the red solid lines, ii) the relaxed constraint \eqref{traffic2} corresponding to the l-th linear segment ($l=1,\ldots,N+1$), is represented with the cutting planes A,B,C,D,E.}
	\label{fig:Ddemw1}
\end{figure}

\noindent Constraint \eqref{traffic2} produces a convex feasibility domain for $\{g_r(k),\rho_r(k)\}$ as shown in Fig. \ref{fig:Ddemw1}.

%\textcolor{blue}{\textbf{26/03/2026 shade the area and change in thesis as well}}\\

Second, constraints \eqref{upoptos} and \eqref{secondtt} are handled together. Similar to \eqref{mfd}, constraint \eqref{secondtt} can be relaxed into the following two inequalities
\begin{align}
	\label{rrcewr}
	\tilde{q}_{o,r,j,s}^{S}(k) &\leq q_{o,r,j,s}^{S}(k),\\
	\label{wwwgbfsx}
	\tilde{q}_{o,r,j,s}^{S}(k) &\leq \frac{C_{r,j}(\rho_j(k))}{q_{r,j}(k)}q_{o,r,j,s}^S(\rho_r(k)).
\end{align}

\noindent As \eqref{wwwgbfsx} is still non-convex, we further relax constraint \eqref{wwwgbfsx} by taking the sum over all $\tilde{q}_{o,r,j,s}^S(k)$ for $o\in\mathcal{O}$, $s\in\mathcal{S}_b, b\in\mathcal{B}$ as follows,
\begin{align}
	\label{wwwwDgg}
	\sum_{o\in\mathcal{O}}\sum_{b\in\mathcal{B}}\sum_{s\in\mathcal{S}_b}\tilde{q}_{o,r,j,s}^S(k) \leq C_{r,j}(\rho_j(k))\frac{\sum_{o\in\mathcal{O}}\sum_{b\in\mathcal{B}}\sum_{s\in\mathcal{S}_b}q_{o,r,j,s}^S(k)}{q_{r,j}(k)}.
\end{align}

\noindent Following a similar procedure for non-convex constraint \eqref{secondtt}, yields 
\begin{align}
	\label{rrcegdgwr}
	\tilde{q}_{o,r,j,d}^{W}(k) &\leq q_{o,r,j,d}^{W}(k),\\
	\label{cgtfdnr}
	\sum_{o\in\mathcal{O}}\sum_{d\in\mathcal{D}}\tilde{q}_{o,r,j,d}^W(k) &\leq C_{r,j}(\rho_j(k))\frac{\sum_{o\in\mathcal{O}}\sum_{d\in\mathcal{D}}q_{o,r,j,d}^W(k)}{q_{r,j}(k)}.
\end{align}

\noindent When summing inequalities \eqref{wwwwDgg} and \eqref{cgtfdnr} by parts, the nominator of the resulting right-hand side expression is canceled out with the denominator term of the resulting expression, $q_{r,j}(k)$ (see Eq. \eqref{ffsarrraga}) and hence we obtain
\begin{align}
	\label{cgtfdnr_over_again}
	\sum_{o\in\mathcal{O}}\sum_{b\in\mathcal{B}}\sum_{s\in\mathcal{S}_b}\tilde{q}_{o,r,j,s}^S(k) + \sum_{o\in\mathcal{O}}\sum_{d\in\mathcal{D}}\tilde{q}_{o,r,j,d}^W(k) \leq C_{r,j}(\rho_j(k)).
\end{align}

\noindent Now, we shift our attention to constraint \eqref{upoptos}. This constraint can be rewritten as
\begin{align}
	\label{rrvrrrr}
	C_{r,j}(\rho_j(k)) = \min\Bigg(C_{r,j}^{\text{MAX}}, \frac{C_{r,j}^{\text{MAX}}}{1-\alpha_{r,j}}\Big(1-\frac{\rho_j(k)}{\rho_j^J}\Big)\Bigg),
\end{align}

\noindent which has the same form with constraint \eqref{traffic2}. Based on that, constraint \eqref{cgtfdnr_over_again} can be further relaxed into the following two linear constraints
\begin{align}
	\label{porta}
	\sum_{o\in\mathcal{O}}\sum_{b\in\mathcal{B}}\sum_{s\in\mathcal{S}_b}\tilde{q}_{o,r,j,s}^S(k) + \sum_{o\in\mathcal{O}}\sum_{d\in\mathcal{D}}\tilde{q}_{o,r,j,d}^W(k) &\leq C_{r,j}^{\textrm{MAX}}, \\
	\label{toixos}
	\sum_{o\in\mathcal{O}}\sum_{b\in\mathcal{B}}\sum_{s\in\mathcal{S}_b}\tilde{q}_{o,r,j,s}^S(k) + \sum_{o\in\mathcal{O}}\sum_{d\in\mathcal{D}}\tilde{q}_{o,r,j,d}^W(k) &\leq \displaystyle\frac{C_{r,j}^{\textrm{MAX}}}{1-\alpha_{r,j}}\Bigg(1-\frac{\rho_j(k)}{\rho_j^J}\Bigg).
\end{align}

\noindent To summarize, constraints \eqref{secondtt} and \eqref{secondtt1} are relaxed into the linear constraints \eqref{rrcewr}, \eqref{rrcegdgwr}, \eqref{porta} and \eqref{toixos}. Let us now examine constraints \eqref{inssss}, \eqref{naruto}, which involve the product of two variables, namely the density $\rho_{o,r,s}^S(k)$ with the speed at region $r$, $u_r(k)$ for constraint \eqref{inssss} and the density $\rho_{o,r,d}^W(k)$ with speed $u_r(k)$ for constraint \eqref{naruto}. Since $u_r(k)\leq u_r^f$, for each permissible value of density $\rho_r(k)$, where $u_r^f$ denotes the free-flow speed of region $r\in\mathcal{R}$, the constraints \eqref{inssss}, \eqref{naruto}, respectively, can be relaxed into
\begin{align}
	\label{porgta}
	q_{o,r,s}^S(k) &\leq \displaystyle u_r^f\rho_{o,r,s}^S(k)\frac{L_r}{l_r}, \\
	\label{toixgos}
	q_{o,r,d}^W(k) &\leq \displaystyle u_r^f\rho_{o,r,d}^W(k)\frac{L_r}{l_r},
\end{align}

\noindent which produce a convex feasibility domain for $\{q_{o,r,s}^S(k), u_r(k), \rho_{o,r,s}^S(k)\}$ in constraint \eqref{inssss} and for $\{q_{o,r,d}^W(k), u_r(k), \rho_{o,r,d}^W(k)\}$ in constraint \eqref{naruto}, each being a superset of the nonconvex feasibility domain of constraints \eqref{inssss} and \eqref{naruto}, respectively.

Beyond the relaxation of the non-convex constraints, an additional hard constraint is incorporated in the relaxed optimization problem to ensure the operational consistency of the outer-loop procedure. Since school start time decisions are implemented through the demand projection mechanism described in Section~\ref{sec:demand_projection}, any variation of $\tilde{\tau}_s^{\lambda}$ between successive outer-loop iterations $\lambda$, induces a corresponding temporal translation of the school-related demand distributions. To promote a smooth evolution of the projected demand (no abrupt changes in successive outer-loop iterations) and enhance the efficiency of the search process, we impose an additional \textit{hard constraint} within optimization problem $P_2$
that limits how much the start time of each school can vary across successive outer-loop iterations $\lambda>0$. This constraint helps avoiding large fluctuations in the decision variables and ensures a more controlled and consistent adjustment of school start times. To this end, we define
\begin{equation}
	\label{minimal_disruption}
	|\tilde{\tau}_s^{\lambda} - \tilde{\tau}_s^{\lambda-1}| \leq V, \quad \forall s \in \mathcal{S}_b, \forall b \in \mathcal{B}, \lambda > 0,
\end{equation}

\noindent where $\tilde{\tau}_s^{\lambda}$ is the shifted start time of school $s$, updated within the iteration $\lambda$ of the Upper-Level module. Constraint \eqref{minimal_disruption} implies that, between consecutive iterations, the demand distribution associated with each school can be translated by at most one discrete shifting interval $V$. Hence, the projected school-related demand profiles evolve gradually across iterations while remaining fully consistent with our proposed demand projection mechanism (see Section \ref{sec:demand_projection}).\par

In short, problem $P_1$ can be transformed into an MILP by replacing constraints \eqref{upoptos}, \eqref{inssss}, \eqref{naruto}, \eqref{secondtt} and \eqref{secondtt1} with constraints \eqref{traffic2} - \eqref{rrcewr}, \eqref{rrcegdgwr}, \eqref{porta} - \eqref{toixgos} and with the additional incorporation of the hard-constraint \eqref{minimal_disruption}. In mathematical programming terms, this problem can be expressed as
\begin{subequations}
	\label{rfrwga}
	\begin{align}
		\label{tr}
		(P_2) \quad \underset{\xi_{m,s}, \forall m,s}{\text{Minimize}} &~ \displaystyle J_{\textrm{TTS}} = T\sum_{k\in\mathcal{K}}\Big(S_a(k) - S_b(k)\Big) \\
		\textrm{Subject To:}&~~ \text{Constraints}~ \eqref{poeee}, \eqref{rrvtvv} - \eqref{fhshshuj},  \eqref{queen} - \eqref{interesting}, \eqref{arxontas1} - \eqref{QUO1},~ \eqref{bin} -  \eqref{binary},\nonumber\\	
		&\quad \quad \quad \quad \quad  ~\eqref{shift_lhhg} - \eqref{projection_works}, \eqref{epsi},~\eqref{traffic2} - \eqref{rrcewr}, \eqref{rrcegdgwr}, \eqref{porta} - \eqref{minimal_disruption}, \nonumber\\
		\textrm{Initialization:}&~~\eqref{init},\nonumber\\
		\textrm{Input:}&~~\eqref{input}\nonumber.
		%\text{Variables:}&~~\eqref{variables}.\nonumber
	\end{align} 
\end{subequations}

\noindent Problem $P_2$ is a Mixed-Integer Linear Program (MILP) and provides a \textit{lower-bound} to the optimal objective value acquired from the solution of Problem $P_1$.  

%\textcolor{blue}{\textbf{22/03/2026 POSSIBLE COMMENT FROM REVIEWER -  how do I guarantee that the solution of optimization problem $P_2$ is always feasible?}}\\

%\textcolor{blue}{\textbf{22/03/2026 POSSIBLE COMMENT FROM REVIEWER - The authors should clarify whether this constraint is part of the optimization problem's feasible set or merely a heuristic regularizer. If it is the former, the optimal solution is constrained in a way that may not be operationally motivated. If it is the latter, its interaction with the convergence proof should be discussed.}}\\

%The bounded change constraint in \eqref{minimal_disruption} is consistent with the concept of bounded rationality and empirical observations of gradual adaptation in commuter behaviour \citep{DI2016142}. Similar bounded update rules have been adopted in dynamic traffic assignment models, where the step-size or allowable deviation is explicitly limited between iterations to promote stability and ensure realistic behavioural responses \citep{SAENZROYO2023254}. In our case, parameter $V$ acts analogously by constraining the variation of school start times, thus capturing the practical limitations in implementing schedule changes.

\begin{algorithm}[t]
	\algsetup{linenosize=\tiny}
	\footnotesize
	%{\fontsize{6}{6}\selectfont
		\begin{algorithmic}[1]
			\STATE \textbf{Input}: $\epsilon, \lambda, V, \tau_s, \forall s\in\mathcal{S}_b, \forall b\in\mathcal{B}$, $\bar{t}_d$, $\bar{t}_d^{flex}$$, \forall d\in\mathcal{D}$.
			\STATE \quad \quad \quad Traffic Network Parameters, External Demand same as Algorithm \ref{alg:user_equilibrium}.			
			\STATE Solve optimization problem $P_2$ by inserting $\epsilon$ as the right hand side of constraint \eqref{epsi}. 
			\STATE \textbf{Output}: Decision vector $\boldsymbol{\xi}$, from which the vector of school start times $\boldmath{\tau}^{\text{LB}}$$(\epsilon,\lambda)$ is retrieved, accompanied by the corresponding objective values during iteration $\lambda$, $J_{\text{TTS}}^{\text{LB},\lambda}$ and  $J_{\text{STC}}^{\text{LB},\lambda}$ using Eqs. \eqref{fixed} and \eqref{QUO1}, respectively.
		\end{algorithmic}
		%	}
	\caption{Lower Bound (LB) Algorithm}
	\label{alg:seq1gsg}
\end{algorithm}	

Algorithm \ref{alg:seq1gsg}, termed as LB\footnote{The Lower Bound Algorithm and the Lower-Level module are completely different aspects. Notice that the Lower Bound algorithm provides a lower-bound to the optimal value stemming from optimization problem $P_1$, while the Lower-Level module (see Section \ref{sec:lower_level_controller}), describes the framework that derives the equilibrium-based demand vectors of commuters given as input the school start time vector $\boldsymbol{\xi}$.}, summarizes the relaxation procedure to obtain a lower-bound to the optimal solution stemming from problem $P_1$ during iteration $\lambda$. The algorithm takes the same inputs with those given in Algorithm \ref{alg:user_equilibrium} with the addition of the outer-loop iteration counter for the Upper-Level module, $\lambda$ and the discrete duration of each school shifting interval, $V$. Algorithm \ref{alg:seq1gsg} solves optimization problem $P_2$ for a specific value of $\epsilon$ and determines the vector of school start times $\boldmath{\tau}^{\text{LB}}$$(\epsilon,\lambda)$ that satisfies the imposed constraints. The solution of problem $P_2$ provides as output the school start time vector $\boldmath{\tau}^{\text{LB}}$$(\epsilon,\lambda)$ together with the corresponding objective values $J_{\text{TTS}}^{\text{LB},\lambda}$ and $J_{\text{STC}}^{\text{LB},\lambda}$, computed using Eqs.~\eqref{fixed} and \eqref{QUO1}, respectively (Line 4).

\subsection{Solution to problem $P_1$ via a MILP-based approach}
\label{sec:upper}

This section proposes the \Coupled\ Algorithm, termed \Coupledacronym. The purpose of the \Coupledacronym\ algorithm, outlined in Algorithm \ref{alg:upper_bound}, is to bridge the LB Algorithm with the AUE Algorithm to compute an \textit{upper-bound} to the optimal value of $P_1$ and by extension, an \textit{upper-bound} to the bilevel problem \BilevelOptwithepsilon\ for a fixed value of $\epsilon$. First, the LB Algorithm solves $P_2$ and returns a \textit{lower-bound} to the optimal value of $P_1$. However, the school start time vector $\boldmath{\tau}^{\text{LB}}$$(\epsilon,\lambda)$ produced by $P_2$ is derived under relaxed (linearised) traffic constraints and an approximated demand projection. The \Coupledacronym\ Algorithm addresses this gap: by evaluating $\boldmath{\tau}^{\text{LB}}$$(\epsilon,\lambda)$ through the nonlinear traffic dynamics (see \eqref{mfd} - \eqref{QUO1}) and the DUE approximation of Algorithm \ApproximateAcronym, it produces a feasible solution to $P_1$ (this problem considers the original nonlinear and nonconvex traffic dynamics dictated by Eqs. \eqref{mfd} - \eqref{QUO1}) and, by extension, a feasible solution to the bilevel problem \BilevelOptwithepsilon. \par 

Starting from the initial demand vectors and a convergence tolerance $\delta_2$, the algorithm proceeds as follows: At each outer iteration 
$\lambda$, the Upper-Level problem $P_2$ is solved to obtain the decision vector $\boldsymbol{\xi}$ and the corresponding shifted school start times $\boldsymbol{\tau}^{\mathrm{LB}}(\epsilon,\lambda)$ (Lines 1–10). In the first iteration, the demand inputs correspond to the initialized demand, whereas in subsequent iterations they are replaced by the equilibrium demand obtained from the previous iteration. The resulting vector $\boldsymbol{\tau}^{\mathrm{LB}}(\epsilon,\lambda)$ is then evaluated via the Lower-Level module by executing Algorithm~\ref{alg:user_equilibrium}, which simulates the nonlinear traffic dynamics and returns equilibrium demand vectors together with the performance metrics $J_{\text{TTS}}^{\Iota(\lambda)}, J_{\text{STC}}^{\Iota(\lambda)}$ (Line 11). Convergence is assessed by monitoring the relative change in TTS between successive outer-loop iterations (Lines 12–15); if it falls below $\delta_2$, the procedure terminates, otherwise the updated equilibrium demand is fed back to the Upper-Level problem and the loop continues. After termination, we obtain as output the TTS and STC values, respectively (Line 18). \par

%When we feed the vector $\boldmath{\tau}^{\text{LB}}$$(\epsilon,\lambda)$ to the nonlinear traffic system dictated by Eqs. \eqref{mfd} - \eqref{QUO1}, the traffic states, i.e., density $\rho_r(k)$ and flow $q_r(k)$ remain within their admissible bounds. Consequently, the school start time decision $\boldmath{\tau}^{\text{LB}}$$(\epsilon,\lambda)$ determined from Algorithm \ref{alg:seq1gsg} that results in the Upper-Level decision vector $\boldsymbol{\xi}$ will also be feasible to the original MINLP problem $P_1$.

%It is worth noting that the TTS value obtained from the LB Algorithm, $J_{\text{TTS}}^{\text{LB}}(\epsilon)$, once the vector $\boldmath{\tau}^{\text{LB}}$$(\epsilon,\lambda)$ is retrieved, corresponds to the relaxed constraints of the traffic system (as incorporated in problem $P_2$). On the contrary, the respective TTS value obtained from the \Coupledacronym\ Algorithm, $J_{\text{TTS}}^{\Coupledacronym}(\epsilon)$, after evaluating the same vector $\boldmath{\tau}^{\text{LB}}$$(\epsilon,\lambda)$, reflects the original nonlinear and nonconvex traffic dynamics defined by Eqs.~\eqref{mfd}--\eqref{QUO1} via the aid of Algorithm \ref{alg:user_equilibrium}. Hence, there can be a discrepancy between these two quantities. 

\begin{algorithm}[t]
	\algsetup{linenosize=\tiny}
	\footnotesize
	\begin{algorithmic}[1]
		\STATE \textbf{Input}: Same as Algorithm \ref{alg:seq1gsg}, tolerance parameter $\delta_2$.
		\STATE \quad \quad \quad Traffic Network Parameters and External Demand same as Algorithm \ref{alg:seq1gsg}.
		\STATE \textbf{Initialization:} $\lambda = 0$, $\Delta = +\infty$, $\mathbf{d}_{o,s}^{S}, \mathbf{d}_{o,d}^{W}$.
		\WHILE{$\Delta \ge \delta_2$}
		\IF{$\lambda=0$}
		\STATE Execute Algorithm \ref{alg:seq1gsg} providing as input to the optimization problem $P_2$ the initial demand vectors $\mathbf{d}_{o,s}^{S}, \mathbf{d}_{o,d}^{W}$.
		\ELSE
		\STATE Execute Algorithm \ref{alg:seq1gsg} providing as input to the optimization problem $P_2$ the equilibrium-based demand vectors $\mathbf{\hat{d}}_{o,s}^{S},\mathbf{\hat{d}}_{o,d}^{W}$ found from previous outer-loop iteration $\lambda-1$.
		\ENDIF
		\STATE \textbf{Output}: Vector of school start times $\boldmath{\tau}^{\text{LB}}$$(\epsilon,\lambda)$. 
		\STATE Execute Algorithm \ref{alg:user_equilibrium} using the vector of school start times $\boldmath{\tau}^{\text{LB}}$$(\epsilon,\lambda)$ and store
		the equilibrium-based demand vectors $ \mathbf{\hat{d}}_{o,s}^{S,\Iota(\lambda)}, \mathbf{\hat{d}}_{o,d}^{W,\Iota(\lambda)}$ and the TTS, STC values, $J_{\text{TTS}}^{\Iota(\lambda)}$, $J_{\text{STC}}^{\Iota(\lambda)}$ using Eqs. \eqref{fixed} and \eqref{QUO1}, respectively, for the iteration $\Iota(\lambda)$ that achieves convergence of Algorithm \ref{alg:user_equilibrium} and for the current outer-loop iteration $\lambda$.
		\IF{$\lambda > 0$}
		\STATE Compute $\Delta = \frac{|J_{\text{TTS}}^{\Iota(\lambda)} - J_{\text{TTS}}^{\Iota(\lambda-1)}|}{J_{\text{TTS}}^{\Iota(\lambda-1)}}$.
		\ENDIF
		\STATE $\lambda = \lambda + 1$
		\ENDWHILE
		\STATE Set $\Lambda = \lambda-1$.
		\STATE  \textbf{Output}: Compute TTS value \text{$J_{\text{TTS}}^{\Coupledacronym}(\epsilon)
			= J_{\text{TTS}}^{\Iota(\Lambda)}$} and subsequently compute STC value $J_{\text{STC}}^{\Coupledacronym}(\epsilon)$ using Eq. \eqref{QUO1}. 
	\end{algorithmic}
	\caption{\Coupled\ (\Coupledacronym) Algorithm}
	\label{alg:upper_bound}
\end{algorithm}

\subsection{Pareto Front Generation}
\label{sec:pareto}

Algorithm \ref{alg:seq1hj} describes the Pareto Front Generation (PFG) process, whose role is to investigate the trade-off between network efficiency and schedule disruption induced by school start time policies. While the \Coupledacronym\ algorithm produces a feasible upper-bound for a fixed value of $\epsilon$, the decision maker is ultimately interested in understanding how system performance improves as larger adjustments to school start times become permissible. The PFG algorithm therefore acts as an outer exploration layer around \Coupledacronym\ algorithm. By systematically varying the allowable total shift in school start times, the algorithm generates a collection of solutions (based on the considered values of $\epsilon$ stored in set $\mathcal{E}$) forming a Pareto front that quantifies the compromise between TTS and STC metrics. The Pareto front enables the decision planner to evaluate how much schedule disruptions is required to achieve specific congestion reduction benefits.

\begin{algorithm}[t]
	\algsetup{linenosize=\tiny}
	\footnotesize
	\begin{algorithmic}[1]
		\STATE \textbf{Input}: $\mathcal{E}$.
		\STATE \quad \quad \quad Traffic Network Parameters, External Demand same as Algorithm \ref{alg:seq1gsg}.
		\FOR{$\epsilon\in\mathcal{E}$}	
		\STATE Execute Algorithm \ref{alg:upper_bound} and store the TTS value, $J_{\text{TTS}}^{\Coupledacronym}(\epsilon)$ and the STC value, $J_{\text{STC}}^{\Coupledacronym}(\epsilon)$.
		\ENDFOR
		\STATE \textbf{Output}: Pareto Front between $J_{\text{TTS}}^{\Coupledacronym}(\epsilon)$ and  $J_{\text{STC}}^{\Coupledacronym}(\epsilon), \forall \epsilon$.
	\end{algorithmic}
	\caption{Pareto Front Generation (PFG) Algorithm}
	\label{alg:seq1hj}
\end{algorithm}

\section{Simulation Results}
\label{ch:simulation}

%\textcolor{blue}{\textbf{18/08/2025 Do a sensitivity analysis for i) low, ii) moderate and iii) high traffic demand entering the network?}}\\

%\textcolor{blue}{\textbf{18/08/2025 To validate the scalability of the proposed optimization approach consider a larger network consisting of a city center with 4 regions and its periphery with 12 regions (in total 16 regions) as in 2024 paper of Geroliminis?}}\\

%\textcolor{blue}{\textbf{04/08/2025 put the partitioned network of San Francisco as in the paper of Geroliminis?}}\\

%\textcolor{blue}{\textbf{04/08/2025 In his 2025 paper Antoniou also uses Gaussian-like shaped demand for the departure time choice model - cite this paper and refer to that for my demand profile}}\\

%\textcolor{blue}{\textbf{19/02/2025 Shall I say that the split ratios $\gamma_{rjd}$ are derived using Dijkstra's algorithm for K-shortest paths (distance-based, $K=3$ for this paper)?}}\\

%\textcolor{blue}{\textbf{04/05/2025 Are we interested in the scalability in terms of the number of regions or the number of schools? Solve the same optimization program for the same architecture considering 30, 50, 70, 90, 120 schools and do a plot like in the transportmetrica paper}}\\

%\textcolor{blue}{\textbf{13/09/2025 say that a similar network partitioning in 5 regions has been performed in \cite{PAIPURI2020102709} for the Dallas network, USA}}\\

\begin{figure}[t]
	\centering
	\begin{subfigure}[t]{0.48\linewidth}
		\includegraphics[width=\linewidth]{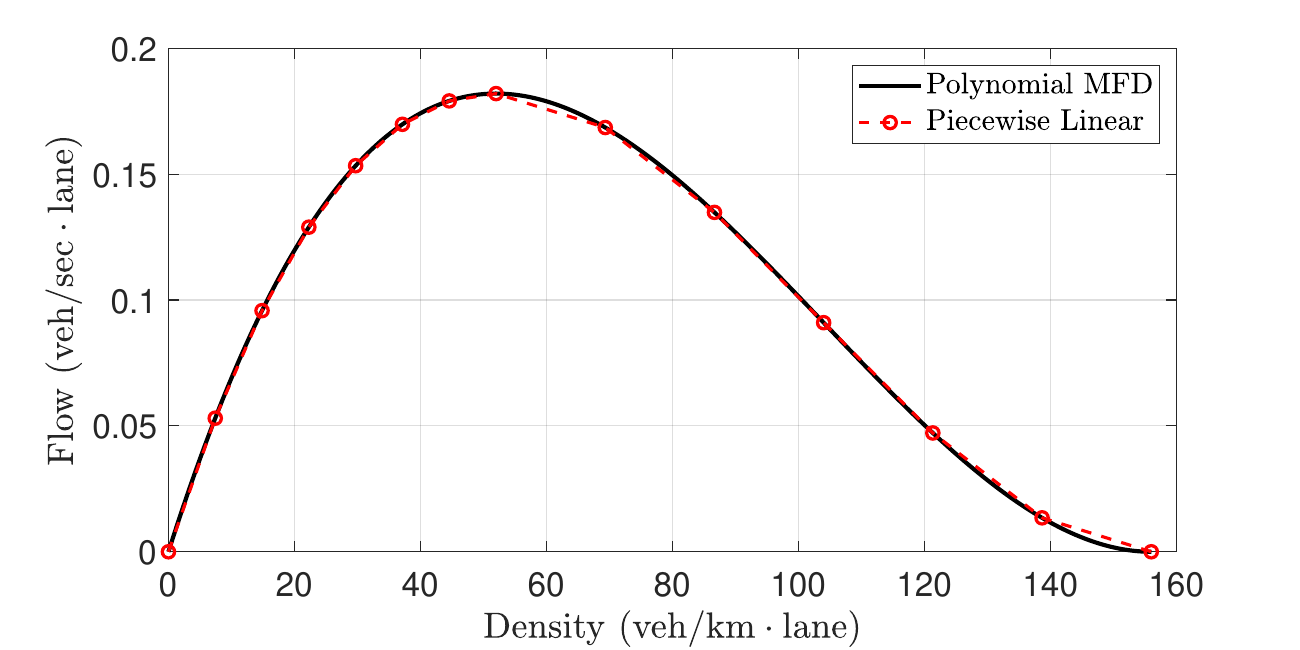}
		\caption{Real and piecewise approximation of the MFD of the city center with $N=15$ segments.}
		\label{fig:MFD_city_compare}
	\end{subfigure}
	\hfill
	\begin{subfigure}[t]{0.48\linewidth}
		\includegraphics[width=\linewidth]{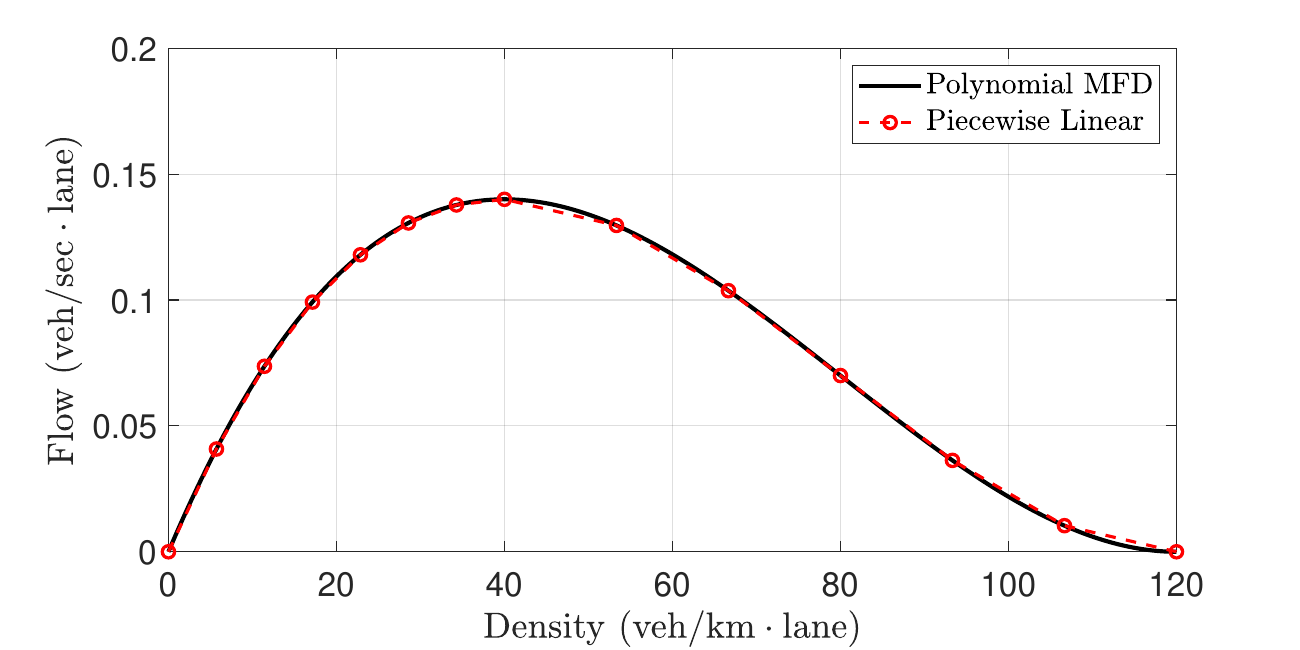}
		\caption{Real and piecewise approximation of the MFD of the periphery with $N=15$ segments.}
		\label{fig:MFD_periphery_compare}
	\end{subfigure}
	\caption{Real and Piece-wise linear approximation of the MFD of regions within (a) the city center and (b) the periphery.}
	\label{fig:MFDs_piecewise}
\end{figure}

To evaluate the performance of the proposed school start time selection methodology, we consider a case study of a five-region synthetic network mimicking a city center and its periphery (Fig. \ref{fig:distea1}). Every region is regarded as an origin, an intermediate destination for home-school-work commuters and a final destination for home-work commuters. The network configuration shown in Fig. \ref{fig:distea1} was also deployed in \cite{doi:10.1287/trsc.2023.0091} for a perimeter control application. %Inspired by \cite{DAGANZO2008771}, we treat the area of downtown San Francisco as if it could be decomposed into sets of homogeneous 1-lane links, similar within each city, e.g., through the visualization of multi-lane links as side-by-side juxtapositions of 1-lane links.
The MFD of each region $r \in \mathcal{R}$ is initially described by a smooth third-order polynomial of traffic density $\rho_r$ as
\begin{equation}
	g_r(\rho_r) = A_r' \rho_r^3 + B_r' \rho_r^2 + C_r'\rho_r, \quad \rho_r \in [0, \rho_{r}^J],
\end{equation}

\noindent where the coefficients $A_r', B_r', C_r'$ are fitted such that $g_r$ reaches its maximum value at the critical density $\rho_{r}^C$, $g_r^C= g_r(\rho_r^C)$, and satisfies $g_r(0) = 0$ and $g_r(\rho_{r}^J) = 0$. The periphery regions share an identical MFD as presented in Fig. \ref{fig:MFD_periphery_compare}, adopted from \cite{geroliminis2008existence} 
with flow 0.14 veh/(sec$\cdot$lane) corresponding to critical density 40 veh/(km$\cdot$lane), and jam density 120 veh/(km$\cdot$lane). Following a similar reasoning to \cite{hajiahmadi2013optimal}, it is assumed that the size of region 1 (city center) is different from the periphery regions, hence the MFD of the city center is the periphery MFDs multiplied by a coefficient (1.3) as shown in Fig. \ref{fig:MFD_city_compare}. Thus, the parameters of the third-order polynomial would be $A_1' = 3.2367\times 10^{-7}$, $B_1' = -1.009\times10^{-4}$, $C_1' = 7.8750\times10^{-3}$ and $A_r' = 5.4698\times 10^{-7}$, $B_r' = -1.3125\times10^{-4}$ and $C_r' = 7.8750\times10^{-3}, r=2,\ldots,5$. \par 

For computational tractability purposes in dynamic traffic modeling and control, we approximate the smooth MFD, $g_r(\rho_r)$ of each region $r$ by a piecewise linear function with $l=1,\ldots,N+1$ segments, where the density breakpoints are denoted by $\rho_{r}^{l}$. These breakpoints are generated via linear interpolation in the two intervals dictating the operation of the MFD: $[0, \rho_r^C]$ and $[\rho_r^C, \rho_r^J]$, with a denser sampling around the critical density $\rho_r^C$. This piecewise approximation enables the embedding of nonlinear MFD behavior into linearized models like the one we introduce in this chapter, rendering it suitable for traffic control applications, which is implemented using MATLAB's built-in linear algebra and numerical computation functions. Specifically, polynomial coefficients are determined by solving a system of linear equations, resulting in a continuous piecewise linear function that approximates $g_r(\rho_r)$.\par 

Fig.~\ref{fig:MFDs_piecewise} illustrates the considered third-order MFDs for the city center and the periphery regions (with red), along with their piecewise linear approximation using $N=15$ segments (with black), respectively. The piecewise-linear approximation is visually almost indistinguishable from the original polynomial MFDs, indicating a very small approximation error over the entire density range. This is expected because the third-order MFDs are smooth and unimodal, which makes them particularly well suited for piecewise-linear representation. Using $N=15$ segments provides a dense discretization of the flow–density relationship, ensuring that the approximation preserves the key characteristics of the original curves, namely the critical density, capacity, and congestion branch slope.\par

The traffic parameters associated with each region $r\in\mathcal{R}$ are displayed in Table \ref{table_example}. The average trip length for vehicles in the periphery is $l_r=3.6$ km, $r=2,\ldots,5$; we set $l_1 = 1.5l_2$. We consider that $|\mathcal{S}_1| = 15, |\mathcal{S}_2|=11, |\mathcal{S}_3| = 7, |\mathcal{S}_4| = 9$ and $|\mathcal{S}_5| = 8$, resulting in fifty schools in total, i.e., $\zeta=50$. We have also included 5800 workplaces scattered across all the regions in the network\footnote{For an average-sized city of approximately 300,000 inhabitants, the number of schools and workplaces can be reasonably estimated based on demographic and employment statistics. Following Eurostat and OECD data, about 18$\%$ of the population are school-aged children, leading to roughly 50 schools assuming an average school size of 300 students per school. Similarly, considering an employment rate of 60$\%$ among working-age adults and an average workplace size of 25 employees, the city would contain approximately 5800 workplaces \citep{eurostat_population_employment, oecd_education}.}.\par

The simulation horizon is set from [6:00 - 9:40] AM, and the simulation time step is set equal to $T = 30$ s. The demand loading procedure mimics a 2-hour peak period, followed by a half-hour off-peak period for cool-down. The initial start time of the schools is set to $\tau_{s} =$ 7:40 AM, $\forall s\in\mathcal{S}_b, \forall b\in\mathcal{B}$, while the fixed work start time is set to $\bar{t}_d = $ 7:20 AM and the flexible work start time range is set to $[\bar{t}_d^l, \bar{t}_d^u] = $[7:50, 9:20] AM among all the workplaces located in destination $d\in\mathcal{D}$. For simplicity, in the case of flexible work schedules, the representative work arrival time $\bar{t}^{flex}_d$ is assumed to correspond to the midpoint of the allowable interval, i.e., $\bar{t}^{flex}_d = (\bar{t}_d^l + \bar{t}_d^u)/2$. This assumption reflects an average behavioral response within the feasible arrival range and facilitates analytical tractability. \par

We consider $V=10$ min discretization intervals for the different school start times, while the maximum shifting time for each school is set equal to 40 minutes (forward or backward shift), corresponding to a value $M=4$. Furthermore, the start time for each school is allowed take nine different readings, i.e., $(2M+1 = 9)$, four for the backward candidate shifts, four for the forward candidate shifts and one for the no-shifting case. For commuters of class W (with fixed or flexible work hours), the shadow values of travel time and schedule delay are identical. This is because, despite the flexibility in departure time for the latter group, the overall opportunity cost of commuting and the disutility associated with travel time are largely determined by the work activity itself, which is typically fixed in location and duration. Flexible schedules may allow modest adjustments in arrival time, but the marginal benefit of such adjustments is relatively small compared to the intrinsic cost of travel and the economic value of time lost en route. Therefore, both fixed and flexible-hour commuters experience comparable shadow values for travel time and schedule penalties, reflecting the underlying importance of work commitments in shaping departure time choices~\citep{noland1995travel,ettema2003modeling}. Hence, the shadow values used in the departure choice model are set to $\alpha^S =8, ~\beta^S = 5, ~\gamma^S = 16$ and $\alpha^W = \alpha^{W,flex} = 14,~ \beta^W = \beta^{W,flex} = 7, ~\gamma^W = \gamma^{W,flex} = 20$ all measured in (\euro/hour).\par

Following \cite{doi:10.1177/0361198105191700122}, we set the background traffic (cumulative demand heading directly to the workplaces) to be equal to 75$\%$ of the trips in the system, while the remaining 25$\%$ represents the cumulative number of vehicle commuters traveling at schools first, before heading to their workplace. Furthermore, for work-related commuters, we assume that the total demand $d_{o,d}^W$ at each destination $d \in \mathcal{D}$ is divided into two groups: a fixed work schedule cohort comprising $\eta_d = 0.8$ (i.e., $80\%$ of the total work-related demand) and a flexible work schedule cohort comprising the remaining $20\%$. Lastly, the values for the split ratios $\theta_{r,j,d}, \forall r\in\mathcal{R}, \forall j\in\mathcal{J}_r^-, \forall d\in\mathcal{D}$ were chosen to ensure that vehicles always progress toward their destination without forming loops, i.e., all regional paths are acyclic, consistent with the assumption we adopted in Section \ref{sec:path_identification}.\\

%In general, within the overall solution framework (see Section \ref{sec:solution}), we develop three distinct solution approaches for the Upper-Level optimization problem $P_1$. In all cases, the corresponding Lower-Level problem, which models commuters’ departure time choices is consistently solved using Algorithm~\ref{alg:user_equilibrium}.

In general, within the overall solution framework (see Section \ref{sec:solution}), we develop and compare three distinct solution algorithms. Each algorithm embeds the solution of the Upper-Level optimization problem $P_1$ within an alternating bilevel optimization procedure, where the Upper-Level decisions are iteratively updated in conjunction with the Lower-Level commuter departure-time equilibrium, which is consistently computed using Algorithm~\ref{alg:user_equilibrium}. The three algorithms differ in the optimization strategy adopted to solve the Upper-Level module.

\begin{enumerate}
	\item \textbf{ABO-MILP}:
	In this algorithm, the Upper-Level optimization problem is solved using a relaxation of the original nonlinear and nonconvex traffic dynamics, leading to the MILP problem $P_2$ described in Section \ref{sec:relaxation}. The resulting MILP is solved using the Gurobi mathematical programming solver (version 12.0.3) \citep{gurobi}. 
	
	%\item \textbf{Lower-bound approach (MILP)}: This method arises from a relaxation of the original nonlinear and nonconvex traffic dynamics, leading to the MILP $P_2$, as described in Section \ref{sec:relaxation}. The problem $P_2$ is solved using the Gurobi mathematical programming solver (version 12.0.3) \citep{gurobi}.
	
	\item \textbf{ABO-SS}:  
	This algorithm solves the Upper-Level optimization problem $P_4$ (see Appendix \ref{sec:derivative_free} for the formulation and solution approach) that will be solved using a derivative-free solver designed to handle the original nonlinear and nonconvex formulation without relying on gradient information, for benchmarking purposes. Problem $P_4$ relies on the \textit{Surrogate derivative-free optimization solver} implemented through the Matlab package \citep{matlab}.
	
	%\item \textbf{Derivative-free approach - Nonlinear Program (NLP)}: This corresponds to problem $P_4$ (see Appendix \ref{sec:derivative_free} for the formulation and solution approach) that will be solved using a derivative-free solver designed to handle the original nonlinear and nonconvex formulation without relying on gradient information, for benchmarking purposes. Problem $P_4$ relies on the \textit{Surrogate derivative-free solver} implemented through the Matlab optimization package \citep{matlab}.
	
	\item \textbf{ABO-ES}:  
	This algorithm constitutes an exact enumeration method that yields the global optimal solution to $P_1$ for a selected value of $\epsilon$, introduced for benchmarking purposes (see Appendix \ref{sec:es_contribution} for the details). Since this is an exhaustive approach, its computational complexity increases exponentially as $\epsilon$ grows, rendering its use prohibitive for larger values of $\epsilon$. For the Exhaustive Search part, we use the Parallel Programming toolbox of Matlab \citep{matlab_parallel_toolbox}.
	
	%\item \textbf{Exhaustive Search approach}: An exact enumeration method that yields the global optimal solution to $P_1$ for a selected value of $\epsilon$, introduced for benchmarking purposes (see Appendix \ref{sec:es_contribution} for the details). Since this is an exhaustive approach, its computational complexity increases exponentially as $\epsilon$ grows, rendering its use prohibitive for larger values of $\epsilon$. For the Exhaustive Search part, we use the Parallel Programming toolbox of Matlab \citep{matlab_parallel_toolbox}.
\end{enumerate}

To evaluate these three developed approaches mentioned above, the following two Relative Deviation (RD) metrics are used for different values of parameter $\epsilon$ belonging to set $\mathcal{E}$
\begin{align*}
	\text{RD}_1(\epsilon) &=\frac{J_{\text{TTS}}^{\Coupledacronym}(\epsilon) - J_{\text{TTS}}^{\ExhaustiveAcronym,*}(\epsilon)}{J_{\text{TTS}}^{\ExhaustiveAcronym,*}(\epsilon)}\times 100\%,\\
	\text{RD}_2(\epsilon) &= \frac{|J_{\text{TTS}}^{\Coupledacronym}(\epsilon) - J_{\text{TTS}}^{\SurrogateAcronym}(\epsilon)|}{\min\Big(J_{\text{TTS}}^{\Coupledacronym}(\epsilon), J_{\text{TTS}}^{\SurrogateAcronym}(\epsilon)\Big)}\times 100\%,
\end{align*}

\noindent where $J_{\text{TTS}}^{\ExhaustiveAcronym,*}(\epsilon)$ denotes the optimal objective value of Total Time Spent when executing the \Exhaustive\ Algorithm that we developed for benchmarking purposes (see Algorithm \ref{alg:serqt} in Appendix \ref{sec:es_contribution}), whereas $J_{\text{TTS}}^{\Coupledacronym}(\epsilon)$, $J_{\text{TTS}}^{\SurrogateAcronym}(\epsilon)$ denote the upper-bound Total Time Spent values acquired from the execution of \Coupledacronym\ (see Algorithm \ref{alg:upper_bound}) and \SurrogateAcronym\ (SS standing for Surrogate Solver, corresponding to Algorithm \ref{alg:seq1gsg_appendix} in Appendix \ref{sec:derivative_free}) Algorithms, using every time as optimization engine the relaxed Upper-Level optimization problem $P_2$ and the Upper-Level optimization problem $P_4$, respectively. $\text{RD}_1(\epsilon)$,  $\text{RD}_2(\epsilon)$ show how good our feasible solution is, as it measures the relative percentage deviation between the TTS value, obtained from the \Coupledacronym\ Algorithm, the \SurrogateAcronym\ Algorithm and the respective optimal value, obtained from the \ExhaustiveAcronym\ Algorithm. All experiments were conducted on a desktop computer with an Intel i9-10980XE  processor (3.0 GHz) with 64 GB of DDR4 RAM running 64-bit Windows 10.\par

%\textcolor{blue}{\textbf{ 21/03/2026 POSSIBLE COMMENT REVIEWER - The CF algorithm solves a relaxed MILP (P2) while the EUB algorithm uses a surrogate derivative-free NLP solver limited to 20 iterations. The derivative-free solver is clearly at a disadvantage by design in terms of allowed computational budget. The comparison in Figure 12 and Table 3 therefore partially reflects an algorithmic budget asymmetry rather than an inherent structural advantage of the MILP relaxation. A fairer comparison would equalize computational budgets. Solve for a specific $\epsilon$ one instance with 20, 50, 150 iterations and see the difference in performance}}\\

%\textcolor{blue}{\textbf{ 21/03/2026 POSSIBLE COMMENT REVIEWER - do a pareto front only for EUB and state the $\epsilon$ values for different values of iterations and motivate the use of 20 iterations initially used in the subsequent experiments}}\\

%\textcolor{blue}{\textbf{21/03/2026 POSSIBLE REVIEWER COMMENT - Readers need to know whether P2 is being solved to proven optimality or to a Gurobi-imposed gap. This significantly affects interpretation of the "lower bound" claim}}\\

As an indication of the problem's complexity, the resulting MILP formulation (Problem $P_2$), derived from the parameter specifications above, involves approximately 6 million constraints and 36 million decision variables, including 450 binary variables. This scale of formulation places the problem well within the class of large-scale, computationally challenging optimization problems. \par 

\begin{table}[t]
	\caption{Traffic parameters used in the case study.}
	\label{table_example}
	\centering
	\footnotesize
	\begin{tabular}{@{}lccccccccc@{}}
		\toprule
		Region & $u_r^f$ & $\rho_r^C$ & $\rho_r^J$ & $g_r^C$ & $q_r^C$ & $C_{r,j}^{\text{MAX}}$ & $\alpha_{r,j}$ & $L_r$ & $l_r$ \\
		& (km/h) & (veh/km$\cdot$lane) & (veh/km$\cdot$lane) & (veh/h$\cdot$lane) & (veh/h$\cdot$lane) & (veh/h$\cdot$lane) &  & (km) & (km) \\
		\midrule
		1 & 25.71 & 52 & 156 & 655.20 & 6673.30 & 6673.30 & 0.33 & 55.00 & 5.40 \\
		2 & 25.71 & 40 & 120 & 504.00 & 5460.00 & 5460.00 & 0.33 & 39.00 & 3.60 \\
		3 & 25.71 & 40 & 120 & 504.00 & 5460.00 & 5460.00 & 0.33 & 39.00 & 3.60 \\
		4 & 25.71 & 40 & 120 & 504.00 & 5460.00 & 5460.00 & 0.33 & 39.00 & 3.60 \\
		5 & 25.71 & 40 & 120 & 504.00 & 5460.00 & 5460.00 & 0.33 & 39.00 & 3.60 \\
		\bottomrule
	\end{tabular}
\end{table}

%\textcolor{blue}{\textbf{12/08/2025 instead of saying without UE, say without departure choice for $\epsilon=0$}}\\

%\textcolor{blue}{\textbf{12/08/2025 say that there is a trend for commuters to depart earlier as the penalty of early arrival is smaller than the penalty of late arrival}}\\

We underline that when considering $\zeta$ schools and $2M+1$ candidate start times for each school, then the total number of different school start time vectors in an exhaustive search setting is $(2M+1)^{\zeta} = 9^{50}$. In principle, evaluating the impact of each candidate pair of school start times of a 50-element vector is a computationally prohibitive task due to the exponential growth of the search space. Here, we can only compute the optimal school start time vector $\boldmath{\tau}^{\ExhaustiveAcronym,*}$ when the value of $\epsilon$ is kept relatively small. To deal with this issue, for higher values of $\epsilon$ we will resort to the the derivative-free approach for benchmarking purposes (see Appendix \ref{sec:derivative_free}) that does not pose such a limitation.\par

\begin{table}[t]
	\centering
	\footnotesize
	\caption{Total number of candidate vectors of school start time $\boldsymbol{\tilde{\tau}}^{\ExhaustiveAcronym}$ for different values of $\epsilon$ utilizing the ABO-ES Algorithm. 
		Values marked as ``N/A'' were omitted due to the intractability of exhaustive enumeration for values of $\epsilon \geq 30$.}
	\label{table_examfrdple2}
	\begin{tabular}{@{}m{1cm} m{3cm} m{2cm} m{2cm}@{}}
		\toprule
		$\epsilon$ & Total combinations & $J_{\text{TTS}}^{\ExhaustiveAcronym,*}$ & $J_{\text{STC}}^{\ExhaustiveAcronym,*}$ \\
		(min) &  & (veh h) & (min) \\
		\midrule
		\multirow{1}{*}{0}   & 1         & 18449      & 0   \\
		\multirow{1}{*}{10}  & 101       & 18111      & 10  \\
		\multirow{1}{*}{20}  & 5101      & 17686      & 20  \\
		\multirow{1}{*}{30}  & 171,801   & N/A        & N/A \\
		\multirow{1}{*}{40}  & 4,341,801 & N/A        & N/A \\
		\multirow{1}{*}{50}  & 83,500,020 & N/A        & N/A \\
		\bottomrule
	\end{tabular}
\end{table}

In Table \ref{table_examfrdple2} we show the total number of all feasible school start time pairs for some selected values of $\epsilon$ along with the acquired smallest TTS value from the \ExhaustiveAcronym\ algorithm. Indeed, even for a small value, i.e., $\epsilon=40$, we end up with 4,341,801 school start time combinations, rendering its use for higher values of $\epsilon$ prohibitive; still even for small values of $\epsilon$ we can get valuable insights as to how our solution stemming from the \Coupledacronym\ algorithm performs against the global optimal solution found from the \ExhaustiveAcronym\ algorithm.\par 

\begin{figure}[t]
	\centering
	\includegraphics[width=0.7\columnwidth]{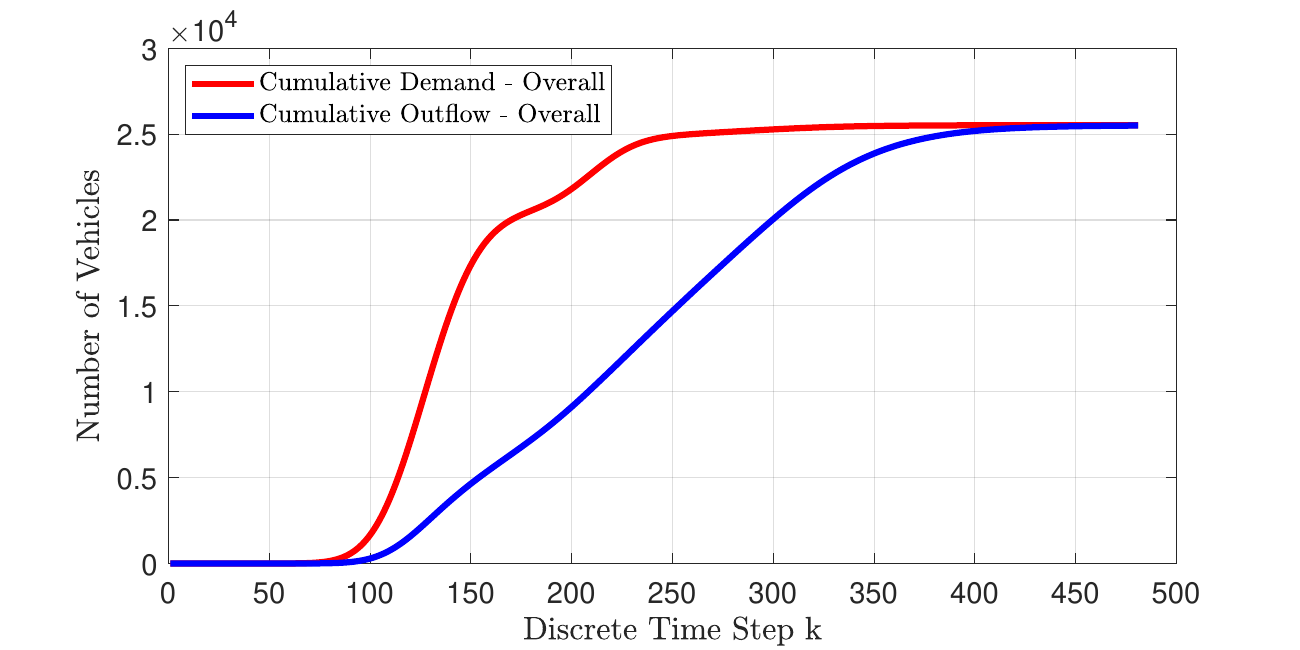} 
	\caption{Cumulative Demand and Cumulative Outflow that enters and exits the network, respectively, when we give as input to the traffic system the initial start time $\tau_s$ for every school $s$, i.e., $\epsilon=0$ min after executing the AUE Algorithm.}
	\label{fig:cumulative_1}
\end{figure}

%\textcolor{blue}{\textbf{13/03/2026 when I put acronym of algorithm in caption it appears without italic. Wrong for captions. Correct it!}}\\

\begin{figure}[t]
	\begin{minipage}{0.52\linewidth}
		\includegraphics[width=\linewidth]{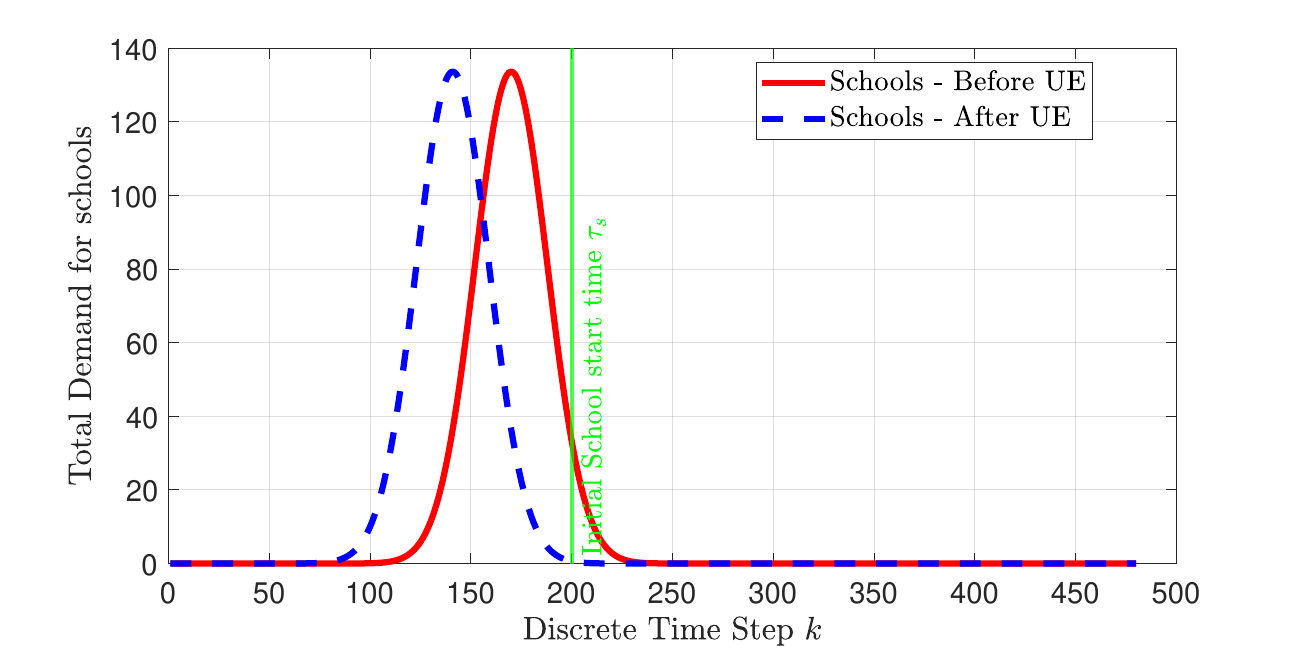}
		\caption*{(a) Distributions before and after the execution of the AUE algorithm for class S when we give as input to the traffic system the initial start time $\tau_s$\\ for every school $s$, i.e., $\epsilon=0$ min.}
		\label{fig:distribution_schools_UE}
	\end{minipage}
	\hfill
	\begin{minipage}{0.52\linewidth}
		\includegraphics[width=\linewidth]{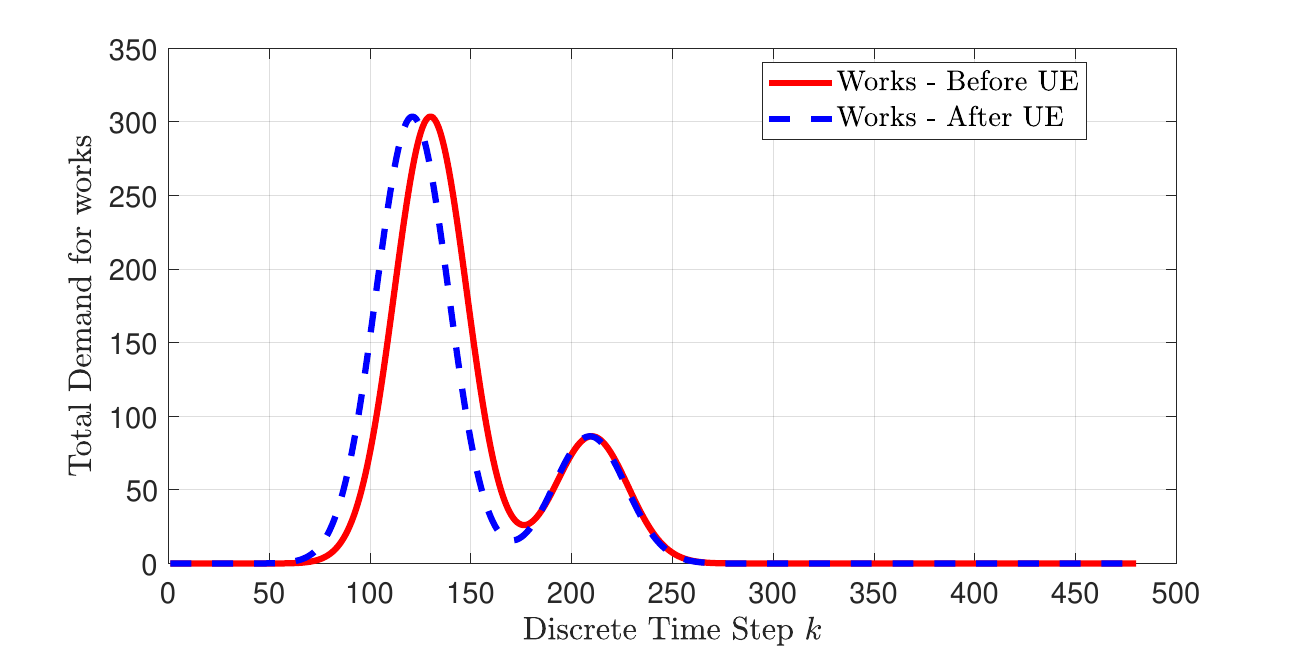}
		\caption*{(b) Distributions before and after the execution of the AUE algorithm\\ for class W when we give as input to the traffic system the initial start time $\tau_s$ for every school $s$, i.e., $\epsilon=0$ min.}
		\label{fig:distribution_works_UE}
	\end{minipage}%
	\hfill
	
	\caption{Distributions before and after the execution of the AUE algorithm for (a) class S and (b) class W, without considering any school start time shift ($\epsilon = 0$ min).}
	\label{fig:distributions_after_UE}
\end{figure}

\begin{figure}[h!t]
	\centering
	\begin{subfigure}{0.48\linewidth}
		\centering
		\includegraphics[width=\linewidth]{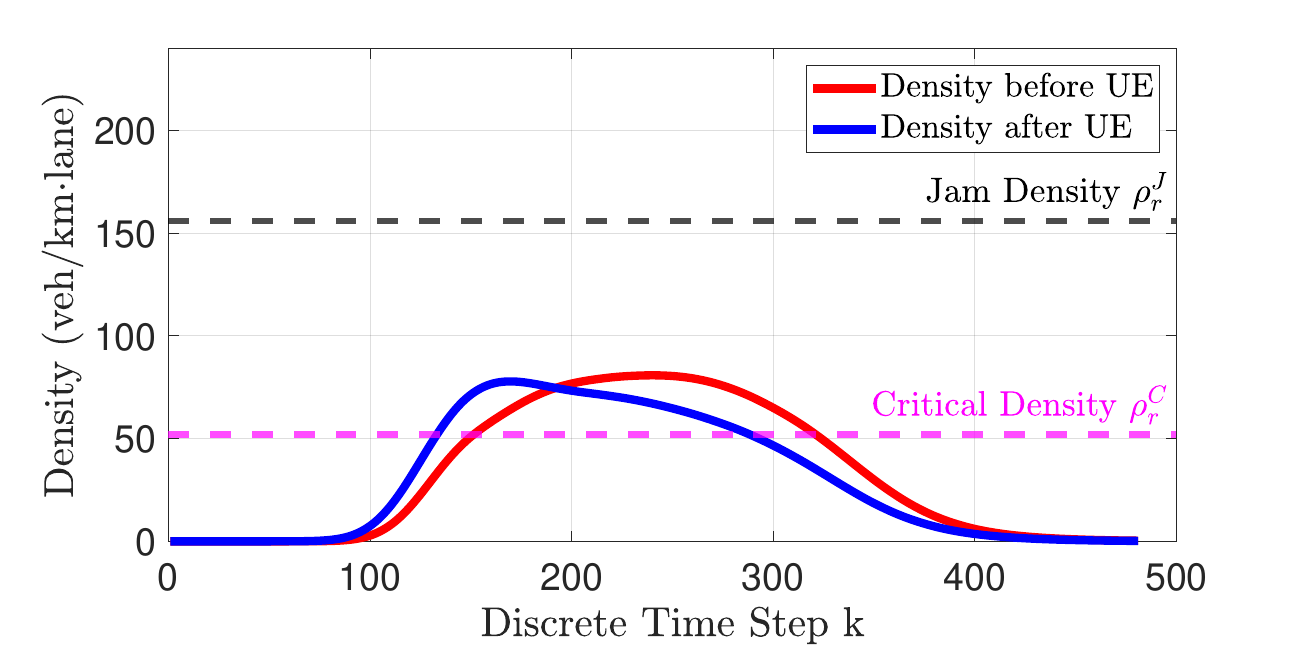}
		\caption{Density evolution in Region 1.}
		\label{fig:Dfdem1}
	\end{subfigure}\hfill
	\begin{subfigure}{0.48\linewidth}
		\centering
		\includegraphics[width=\linewidth]{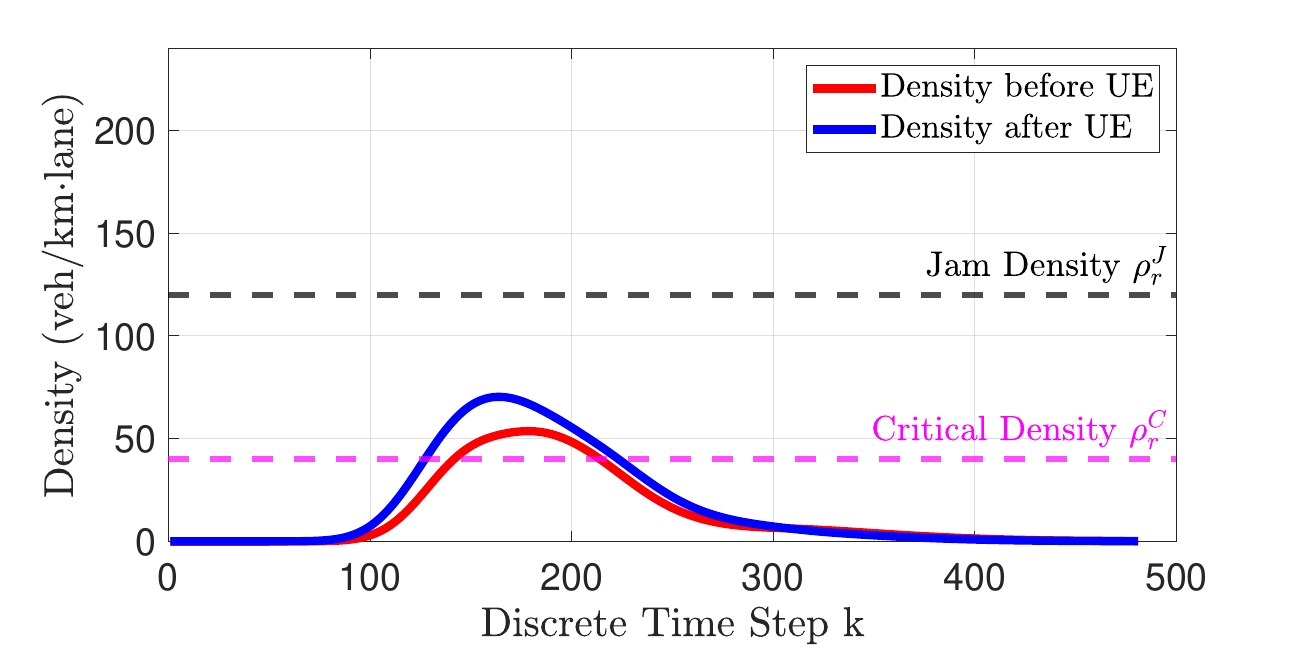}
		\caption{Density evolution in Region 2.}
		\label{fig:Dedfm2}
	\end{subfigure}
	
	\vspace{0.3cm}
	
	\begin{subfigure}{0.48\linewidth}
		\centering
		\includegraphics[width=\linewidth]{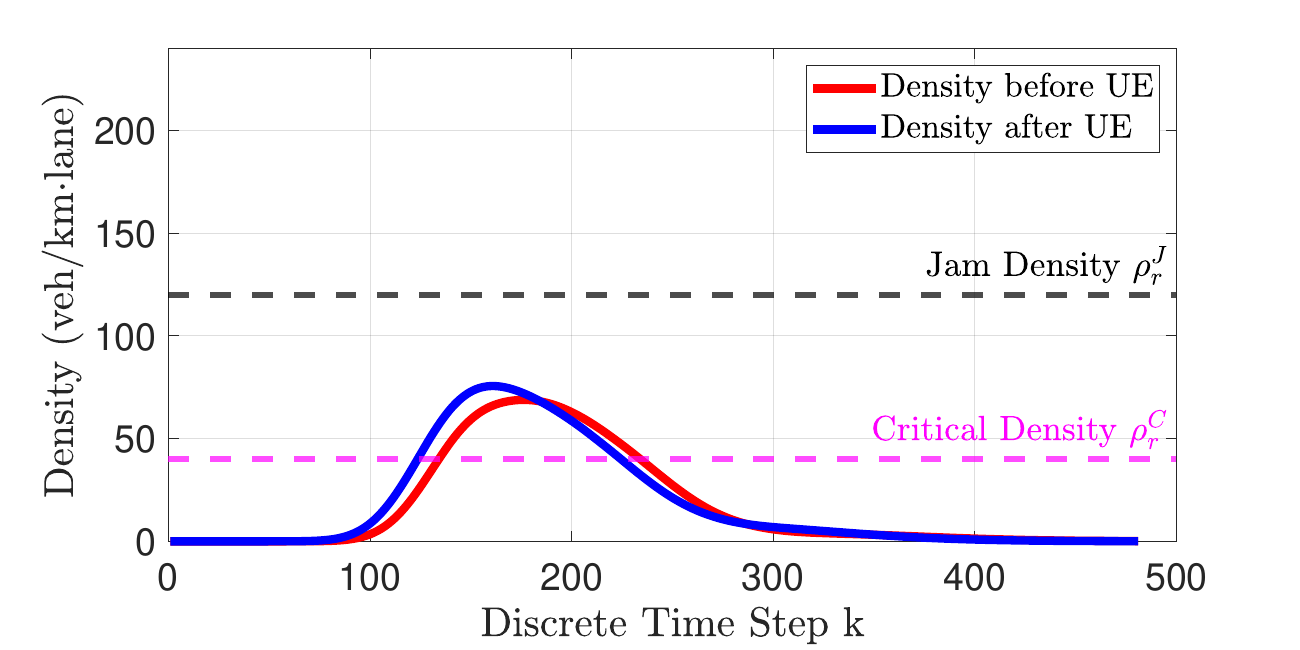}
		\caption{Density evolution in Region 3.}
		\label{fig:Defm3}
	\end{subfigure}\hfill
	\begin{subfigure}{0.48\linewidth}
		\centering
		\includegraphics[width=\linewidth]{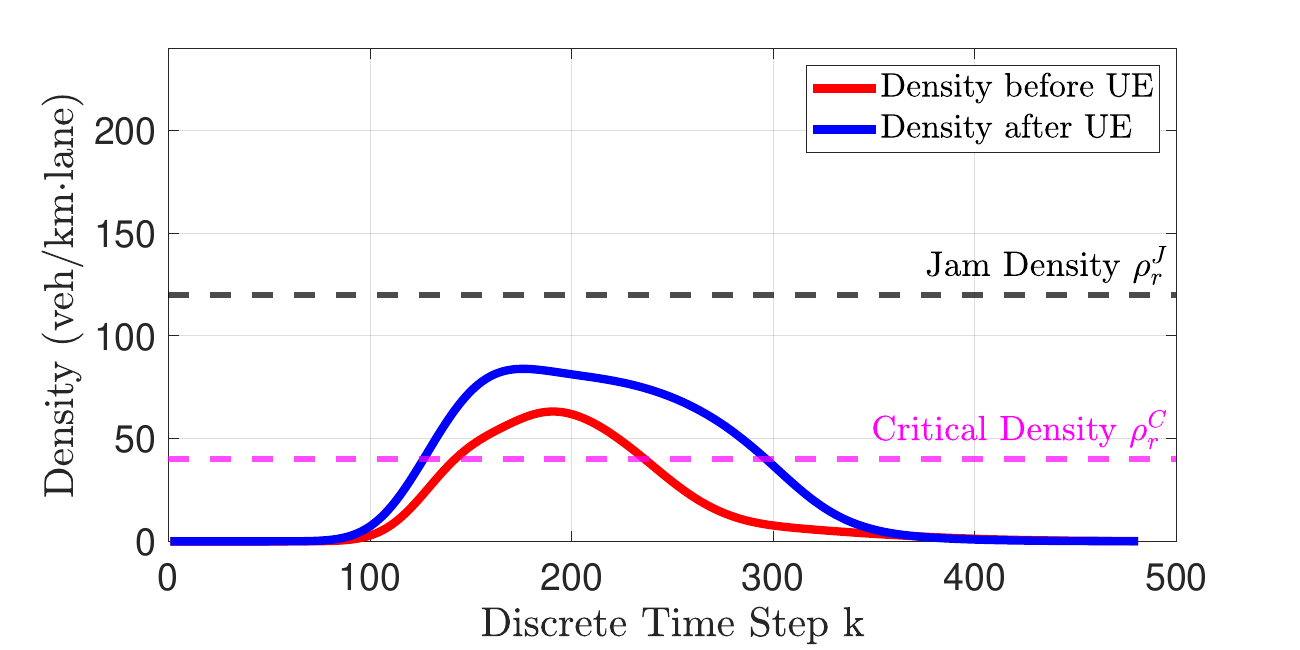}
		\caption{Density evolution in Region 4.}
		\label{fig:Demf4}
	\end{subfigure}
	
	\vspace{0.3cm}
	
	\begin{subfigure}{0.48\linewidth}
		\centering
		\includegraphics[width=\linewidth]{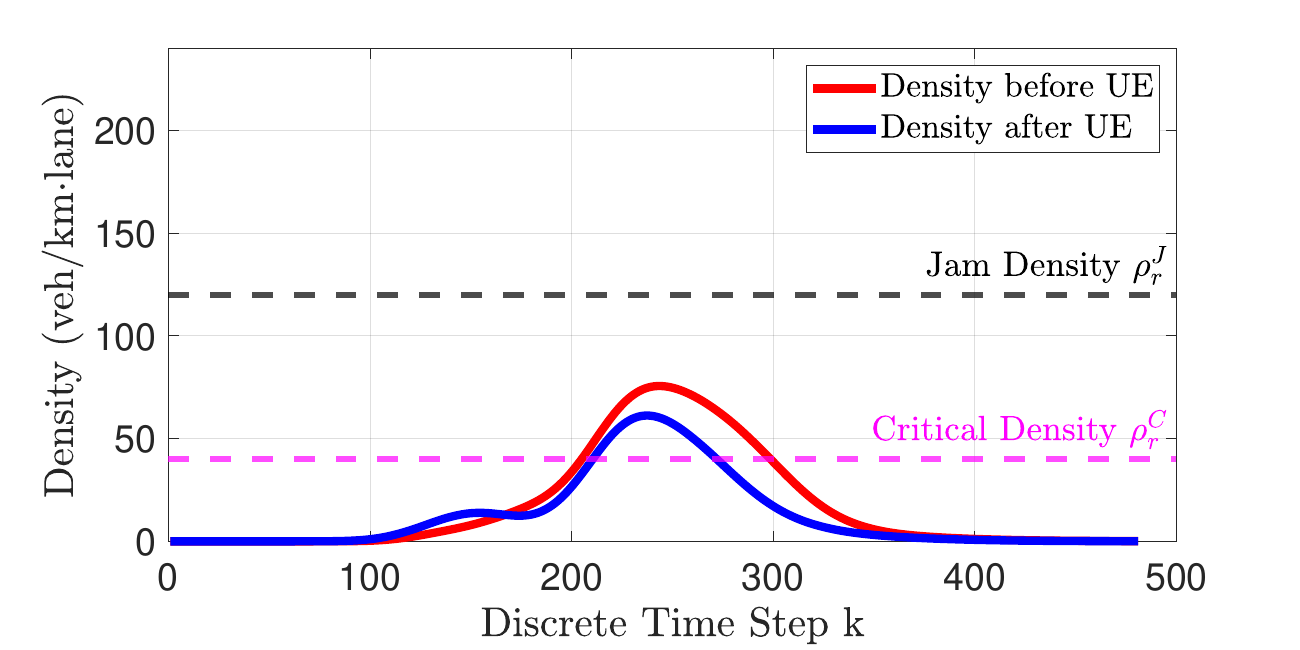}
		\caption{Density evolution in Region 5.}
		\label{fig:Defm5}
	\end{subfigure}
	
	\caption{Density evolution before and after the execution of the AUE algorithm for the case where no school start time shifting takes place, i.e., $\epsilon = 0$.}
	\label{fig:dens_evolf1}
\end{figure}

In Fig. \ref{fig:cumulative_1} we illustrate the cumulative demand versus the cumulative outflow that enters and exits the network, respectively, when no school start shifting is performed ($\epsilon=0$) and after executing the \ApproximateAcronym\ algorithm, indicating the emergence of very high congestion. To better illustrate the departure choice behaviour of commuters as an outcome of the \ApproximateAcronym\ algorithm, in Fig. \ref{fig:distributions_after_UE} we exhibit the overall distribution of demand for each class before and after the execution of the \ApproximateAcronym\ algorithm for the case with $\epsilon=0$ min. This figure suggests that even without a school start time shift taking place in the network, commuters of both classes have indeed an incentive to alter their initial departure times aiming to minimize their perceived travel cost. To further support this observation, Fig. \ref{fig:dens_evolf1} depicts the evolution of density over the entire morning commute period before and after the execution of the \ApproximateAcronym\ algorithm for the no school start time shifting case. The formation of a pattern in the behaviour of commuters can be observed. They choose to depart earlier in view of the initial school start time $\tau_s$ and the work start times (fixed and flexible), $\bar{t}_d, \bar{t}_d^{flex}$, respectively. This can be attributed to the fact that the shadow value for early arrival is smaller than the one imposed for late arrival. Interestingly, a portion of school-related commuters are willing to depart after the initial school start time $\tau_s$, as they regard that doing so would be beneficial to them. From this figure it is also apparent that properly bridging the departure time choice of commuters in view of the selected school start time and the work start times as well, is critical, as it can have considerable impact in the traffic dynamics, especially during the the morning rush hour.

Table \ref{table_examrple2E} shows the performance of our proposed solution methodology, when considering the \Coupledacronym\ algorithm (that makes use of relaxed optimization problem $P_2$) and the \SurrogateAcronym\ Algorithm (that makes use of optimization problem $P_4$ defined in Appendix \ref{sec:derivative_free}), respectively, compared to the \ExhaustiveAcronym\ algorithm for different values of $\epsilon$. We can observe that in the instances that the \ExhaustiveAcronym\ algorithm was able to retrieve a tractable solution, both relative deviation terms are below 3$\%$. An interesting finding is that the TTS value derived from the \SurrogateAcronym\ Algorithm, i.e., $J_{\text{TTS}}^{\SurrogateAcronym}$ perfectly matches the solution obtained from the \ExhaustiveAcronym\ algorithm. It is important to highlight that for $\epsilon\geq 20$, where the \ExhaustiveAcronym\ algorithm cannot provide the global optimal solution due to exponential increase in complexity, the upper-bound TTS values $J_{\text{TTS}}^{\Coupledacronym}$ and $J_{\text{TTS}}^{\SurrogateAcronym}$ are quite close. In fact, for $\epsilon\geq 100$, our proposed \Coupledacronym\ algorithm consistently outperforms the \SurrogateAcronym\ algorithm (that makes use of a derivative-free optimization solver).\par 

As observed in Tables \ref{table_examfrdple2} and \ref{table_examrple2E}, the \ExhaustiveAcronym\ algorithm fails to provide results (indicated as N/A) as the parameter $\epsilon$ increases. This is due to the exponential growth of the feasible region, which renders brute-force evaluation computationally prohibitive even for medium-scale networks. On the contrary, our proposed \Coupledacronym\ and \SurrogateAcronym\ algorithms, respectively, consistently converge to high-quality solutions across all tested scenarios. These 'N/A' instances underscore the necessity of our optimization-based approaches; without them, the bi-level optimization of school start times in a realistic multi-region MFD context would be numerically intractable for city-scale planning.\par

\begin{table}[t]
	\centering
	\footnotesize
	\caption{TTS and STC values derived from different values of $\epsilon$ stemming from the ABO-ES, ABO-MILP and ABO-SS algorithms, respectively.}
	\label{table_examrple2E}
	\begin{tabular}{@{}m{0.8cm} m{1.2cm} m{1.2cm} m{1.2cm} m{1.2cm} m{1.2cm} m{1.2cm} m{1.2cm} m{1.2cm}@{}}
		\toprule
		$\epsilon$ & $J_{\text{TTS}}^{\ExhaustiveAcronym,*}$ & $J_{\text{STC}}^{\ExhaustiveAcronym,*}$ & $J_{\text{TTS}}^{\Coupledacronym}$ & $J_{\text{STC}}^{\Coupledacronym}$ & $J_{\text{TTS}}^{\SurrogateAcronym}$ & $J_{\text{STC}}^{\SurrogateAcronym}$ & RD$_1$ & RD$_2$ \\
		(min) & (veh h) & (min) & (veh h) & (min) & (veh h) & (min) & (\%) & (\%) \\
		\midrule
		\multirow{1}{*}{0}   & 18449 & 0   & 18449 & 0   & 18449 & 0   & 0    & 0    \\
		\multirow{1}{*}{10}  & 18111 & 10  & 18367 & 10  & 18111 & 10  & 1.41 & 1.41 \\
		\multirow{1}{*}{20}  & 17686 & 20  & 18171 & 20  & 17786 & 20  & 2.16 & 2.74 \\
		\multirow{1}{*}{50}  & N/A   & N/A & 17467 & 50  & 17387 & 50  & N/A  & 0.46 \\
		\multirow{1}{*}{100} & N/A   & N/A & 16973 & 100 & 17178 & 100 & N/A  & 1.21 \\
		\multirow{1}{*}{200} & N/A   & N/A & 16917 & 200 & 17047 & 110 & N/A  & 0.77 \\
		\multirow{1}{*}{300} & N/A   & N/A & 16145 & 290 & 17135 & 260 & N/A  & 6.13 \\
		\multirow{1}{*}{600} & N/A   & N/A & 15275 & 600 & 16494 & 560 & N/A  & 7.98 \\
		\bottomrule
	\end{tabular}
\end{table}

\begin{figure}[h!t]
	\centering
	\includegraphics[width=13cm]{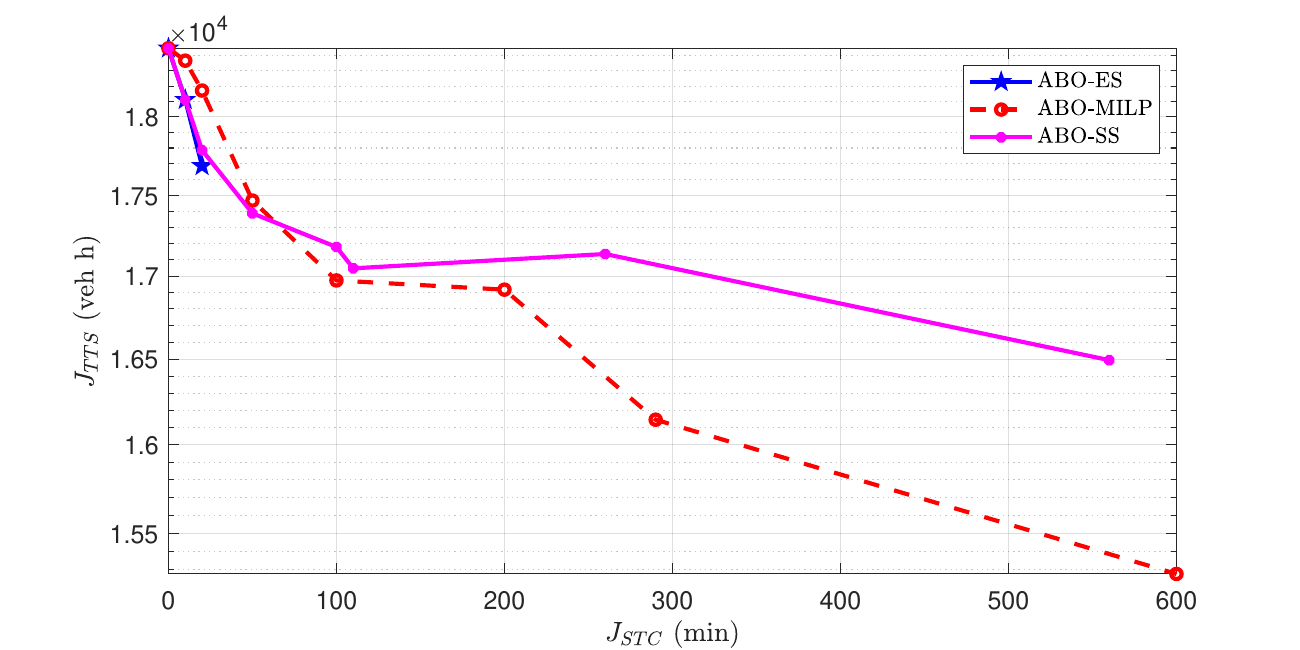} 
	\caption{Pareto Front - Trade-off in the values of the objective functions of interest obtained from the ABO-ES, ABO-MILP and ABO-SS Algorithms, respectively.}
	\label{fig:portabrlre}
\end{figure}

Fig. \ref{fig:portabrlre} illustrates the Pareto front between the two objectives obtained from the \ExhaustiveAcronym\ and PFG algorithms, respectively. For computational tractability purposes, we choose to execute the \ExhaustiveAcronym\ algorithm up to $\epsilon=20$. For comparison purposes, the (i) \Coupledacronym\ algorithm (via the BPFG algorithm) and (ii) the \SurrogateAcronym\ algorithm are executed independently for some selected values of $\epsilon$ within the range $[0, 600]$. Indeed, Table \ref{table_examrple2E} indicates that as $\epsilon$ grows, the corresponding TTS values obtained from the \ExhaustiveAcronym\ and \Coupledacronym\ algorithms, respectively, monotonically decrease, where for $\epsilon=600$ min, we can obtain a 17.20$\%$ reduction in the value of TTS. In practical terms, a value of $\epsilon=600$ min indicates an average of 12 minute start time shift applied to every school among the 50 schools considered in the network. This indicates a very small disruption in the operation of the traffic system. A similar behaviour is displayed from the \SurrogateAcronym\ algorithm, but in some cases the TTS value can even increase as $\epsilon$ grows. \par 

\begin{figure}[h!t]
	\centering
	\includegraphics[width=13cm]{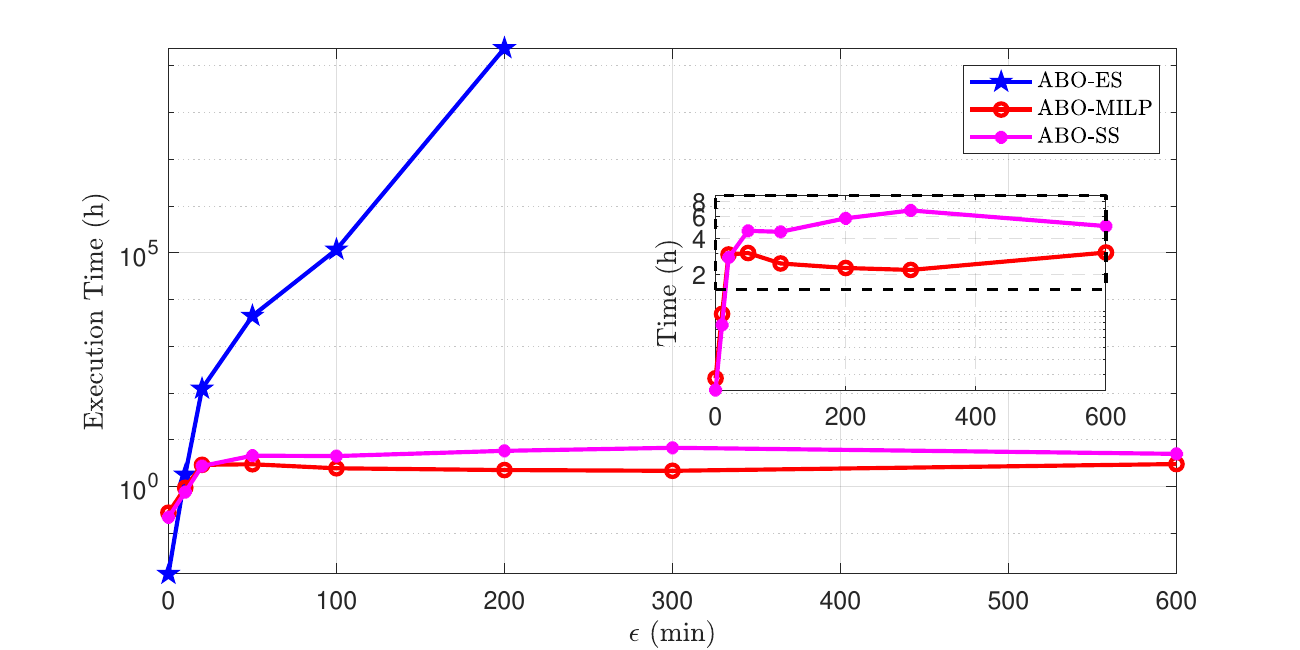} 
	\caption{Execution time needed to obtain the school start time vectors for some selected values of $\epsilon$ derived from the ABO-ES, ABO-MILP and ABO-SS Algorithms, respectively.}
	\label{fig:portabrltre}
\end{figure}

%\textcolor{blue}{\textbf{18/09/2025 Timotheou: Quantify how many times faster is our proposed method compared to the derivative-free based method from Fig. 12. }}\\

Fig. \ref{fig:portabrltre} depicts the execution time needed (in log-scale y-axis) to obtain the school start time vectors derived from the \ExhaustiveAcronym\, \Coupledacronym\ and \SurrogateAcronym\ Algorithms, respectively, for some selected values of $\epsilon$. This figure exhibits that the execution time needed to obtain the global optimal solution (from the \ExhaustiveAcronym\ Algorithm) increases exponentially, since the algorithm has to compute each permissible combination of school start times depending on the value of $\epsilon$. For example, when considering $\epsilon=40$ (min), the total number of different school start time vectors is equal to 4,341,801. In this case, the execution time needed from the \ExhaustiveAcronym\ Algorithm to analytically evaluate the impact of each candidate vector of school start times to the TTS including the departure time choice behaviour of commuters, dictated by the User Equilibrium Algorithm is approximately $10^5$ hours. On the other hand, the time needed to obtain the vector of school start times from the \Coupledacronym\ and \SurrogateAcronym\ Algorithms is marginally above the one hour mark for each selected value of $\epsilon$. Another critical finding stemming from this figure is that the execution time of the \Coupledacronym\ algorithm remains, on average, approximately the same. Note that our MILP-based optimization method (based on problem $P_2$), stemming from the \Coupledacronym\ Algorithm (shown with the red dashed line) is on average 1.65 times faster than the derivative-free one (based on problem $P_4$ obtained from the \SurrogateAcronym\ algorithm, shown with the purple solid line). This finding indicates that our \Coupledacronym\ algorithm is superior both in terms of solution quality and execution time. \par

\begin{figure}[t]
	\begin{minipage}[t]{0.51\linewidth}
		\centering
		\includegraphics[width=\linewidth]{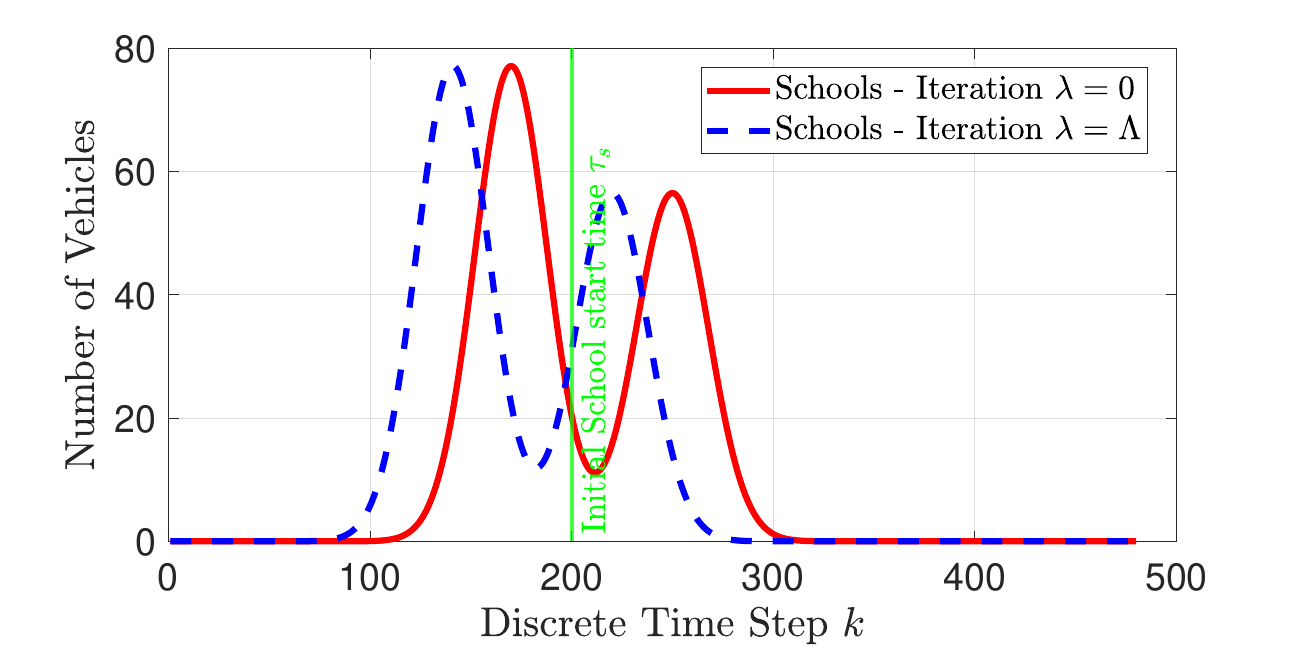}
		\caption*{(a) Equilibrium-based distribution for school-related commuters (class~S) for $\epsilon=600$ min after executing the ABO-MILP algorithm, shown for the first outer-loop iteration ($\lambda = 0$) and for the iteration that achieves convergence ($\lambda = \Lambda$).}
		\label{fig:distribution_schools_UE_shifted}
	\end{minipage}
	\hfill
	\begin{minipage}[t]{0.51\linewidth}
		\centering
		\includegraphics[width=\linewidth]{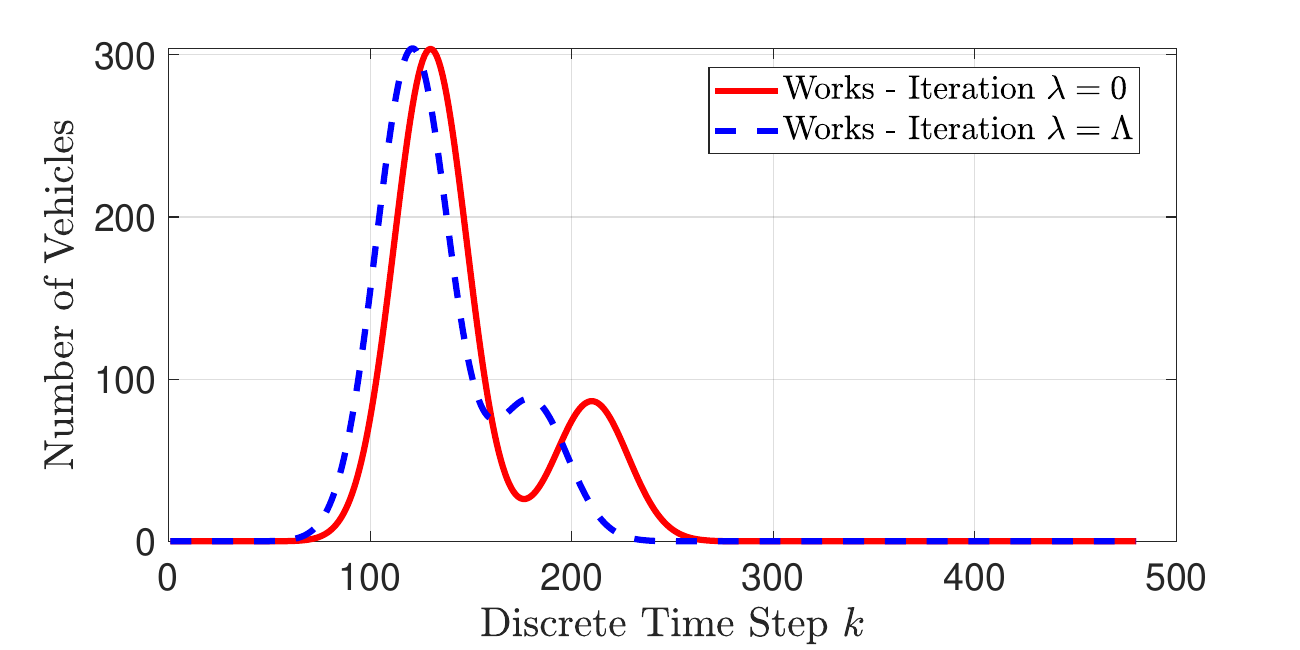}
		\caption*{(b) Equilibrium-based distribution for work-related commuters (class~W) for $\epsilon=600$ min after executing the ABO-MILP algorithm, shown for the first outer-loop iteration ($\lambda = 0$) and for the iteration that achieves convergence ($\lambda = \Lambda$).}
		\label{fig:distribution_works_UE_shifted}
	\end{minipage}
	\caption{Equilibrium-based demand distributions for school and work-related commuters for the case with $\epsilon=600$ min, when executing the ABO-MILP algorithm, comparing the first outer-loop iteration ($\lambda = 0$) and the iteration at which the algorithm converges ($\lambda = \Lambda$).}
	\label{fig:distributions_after_UE_shifted}
\end{figure}

\begin{figure}[h!t]
	
	\begin{subfigure}{0.48\linewidth}
		\centering
		\includegraphics[width=\linewidth]{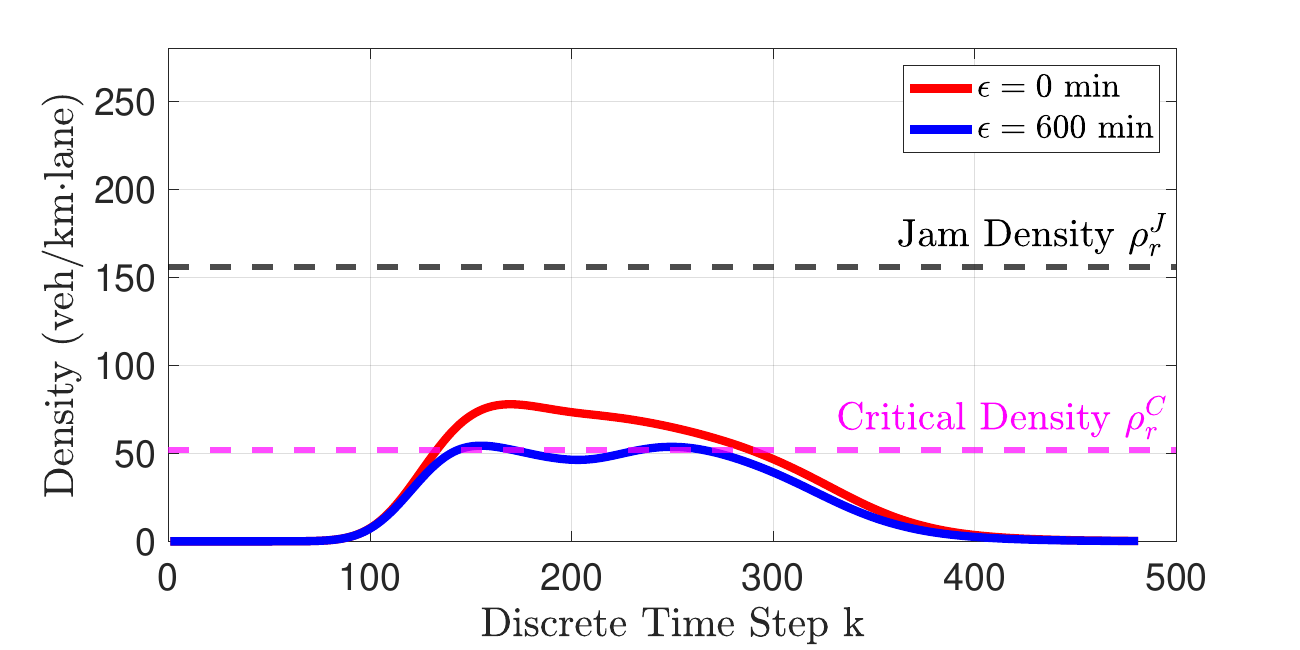}
		\caption*{(a) Density evolution in Region 1.}
		\label{fig:density_UE_1}
	\end{subfigure}
	\hfill
	\begin{subfigure}{0.48\linewidth}
		\centering
		\includegraphics[width=\linewidth]{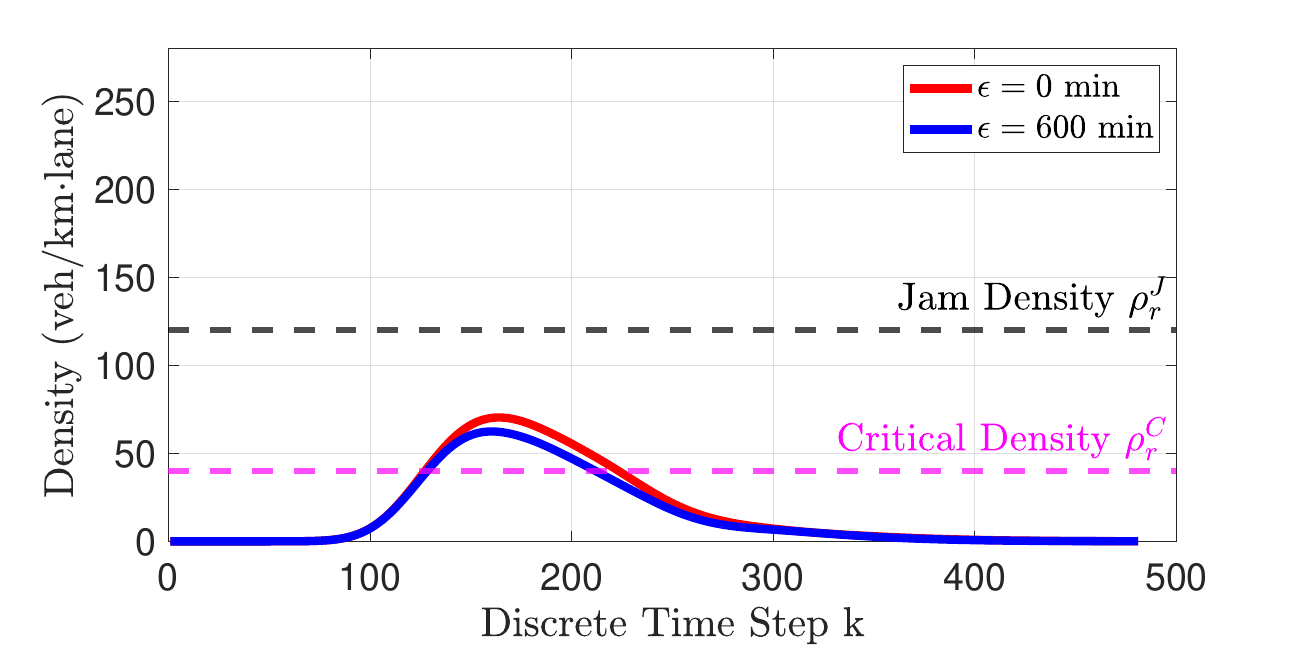}
		\caption*{(b) Density evolution in Region 2.}
		\label{fig:density_UE_2}
	\end{subfigure}%
	\hfill
	
	\vspace{0.3cm}
	
	\begin{subfigure}{0.48\linewidth}
		\includegraphics[width=\linewidth]{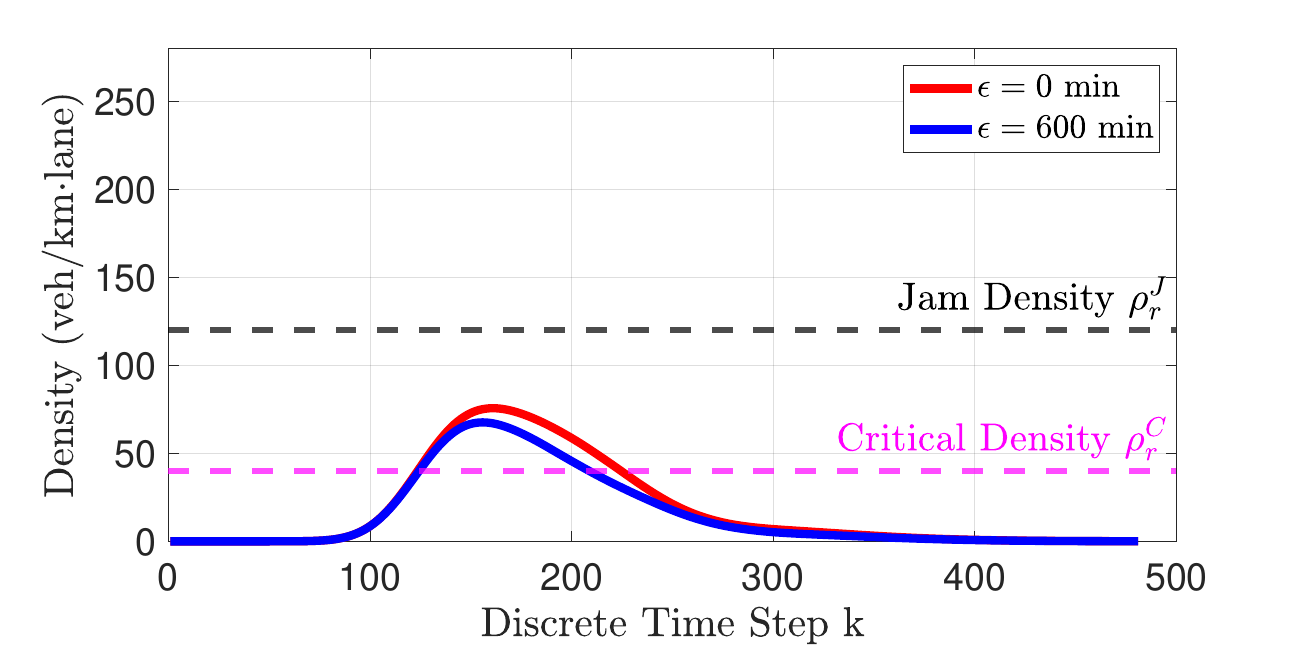}
		\caption*{(c) Density evolution in Region 3.}
		\label{fig:density_UE_3}
	\end{subfigure}
	\hfill
	\begin{subfigure}{0.48\linewidth}
		\includegraphics[width=\linewidth]{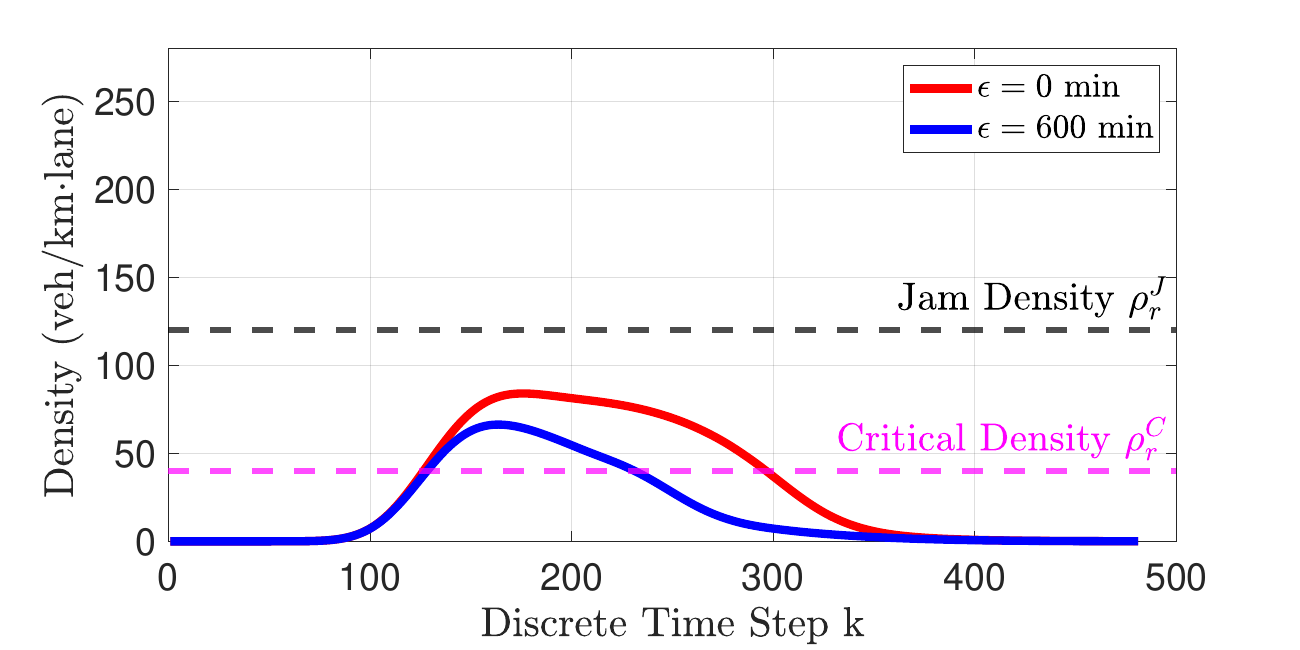}
		\caption*{(d) Density evolution in Region 4.}
		\label{fig:density_UE_4}
	\end{subfigure}
	
	\vspace{0.3cm}
	
	\begin{subfigure}{\linewidth}
		\centering
		\includegraphics[width=0.48\linewidth]{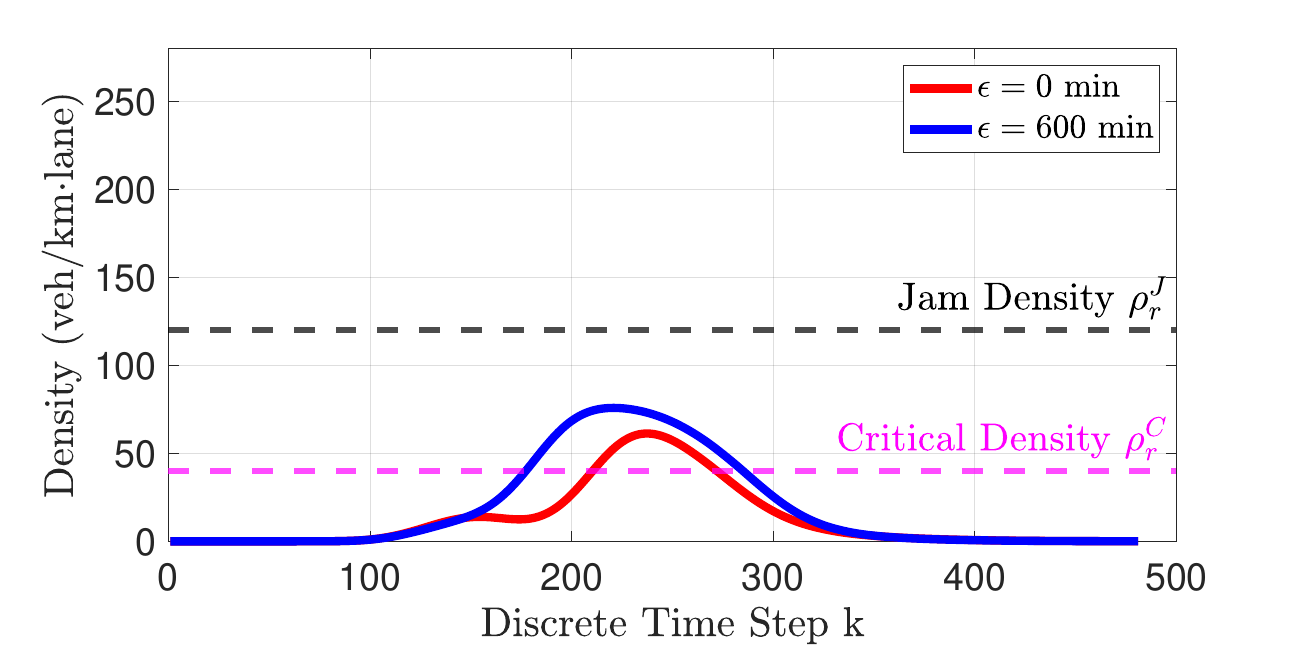}
		
		\vspace{0.2em}
		
		{\centering \caption*{(e) Density evolution in Region 5.}}
		
		\label{fig:density_UE_5}
	\end{subfigure}
	
	\caption{Density evolution when executing ABO-MILP algorithm using as input the school start times corresponding to the values ($\epsilon=0$) and ($\epsilon=600$) min, respectively.}
	\label{fig:dens_evolf1_new}
\end{figure}

Fig.~\ref{fig:distributions_after_UE_shifted} depicts the demand distributions for the two classes of commuters for the case where $\epsilon=600$ min, after executing the \Coupledacronym\ algorithm, comparing the first outer-loop iteration ($\lambda = 0$) and the iteration at which the algorithm converges ($\lambda = \Lambda$). Given these two distributions of demand, Fig. \ref{fig:dens_evolf1_new} illustrates the impact of the school start time change taking into account the departure time preferences of commuters to the traffic system. We can observe that the traffic conditions are improved, leading to smaller peak values of density compared to the no-shifting case, i.e., $\epsilon=0$.

\begin{figure}[t]
	\centering
	\includegraphics[width=13cm]{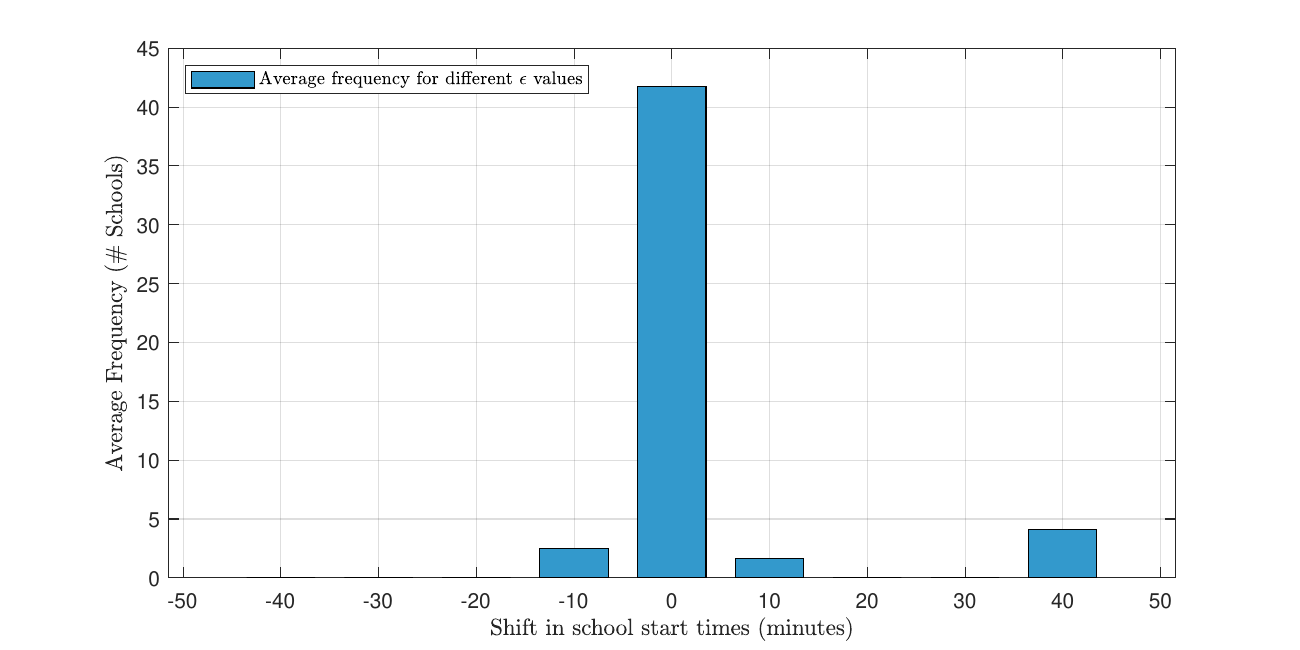} 
	\caption{Average shifting performed to the school start times among all the considered values of $\epsilon$ with respect to the ABO-MILP Algorithm, shown in Table \ref{table_examrple2E}.}
	\label{fig:frequency_distribution}
\end{figure}

Next, Fig.~\ref{fig:frequency_distribution} illustrates the distribution of the average school start time shifts across all the considered values of $\epsilon$ when executing \Coupledacronym\ algorithm. The horizontal axis represents the shift in minutes relative to the baseline schedule, where negative values indicate an earlier start and positive values indicate a delayed start. The vertical axis shows the average frequency with which each shift level was derived. The results reveal that the majority of schools retain their baseline start time (0~minutes shift), while moderate positive and negative shifts (e.g., $\pm$10~minutes) also occur with notable frequency. Larger shifts are comparatively less frequent, indicating that the optimization tends to favor adjustments that balance schedule flexibility with minimal disruption to the existing school timetable. This distribution provides a concise overview of how the proposed framework reallocates start times under varying levels of allowable flexibility.

\subsection{Sensitivity Analysis to the parameters of the MFD}

\begin{figure}[h!t]
	\centering
	\begin{minipage}{0.48\linewidth}
		\includegraphics[width=\linewidth]{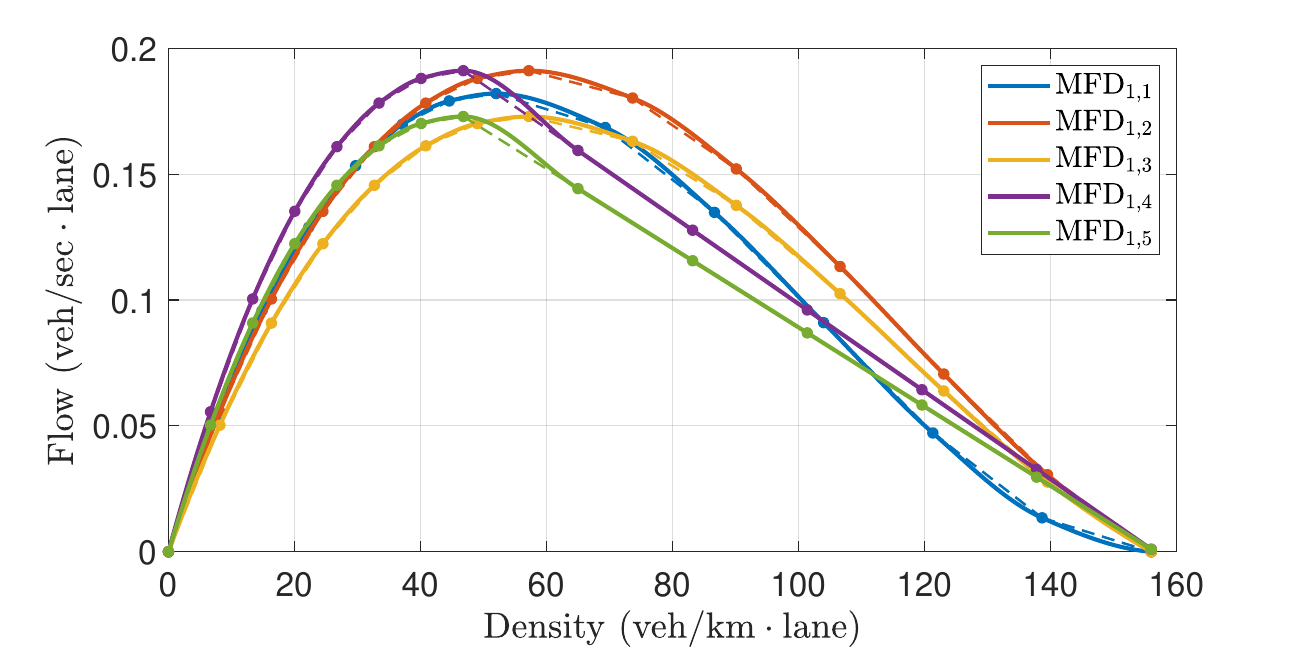}
		\caption*{(a) Real and Piece-wise approximation for noisy MFD of the city center with $N=15$ segments.}
		\label{fig:noisyMFD_1}
	\end{minipage}
	\hfill
	\begin{minipage}{0.48\linewidth}
		\includegraphics[width=\linewidth]{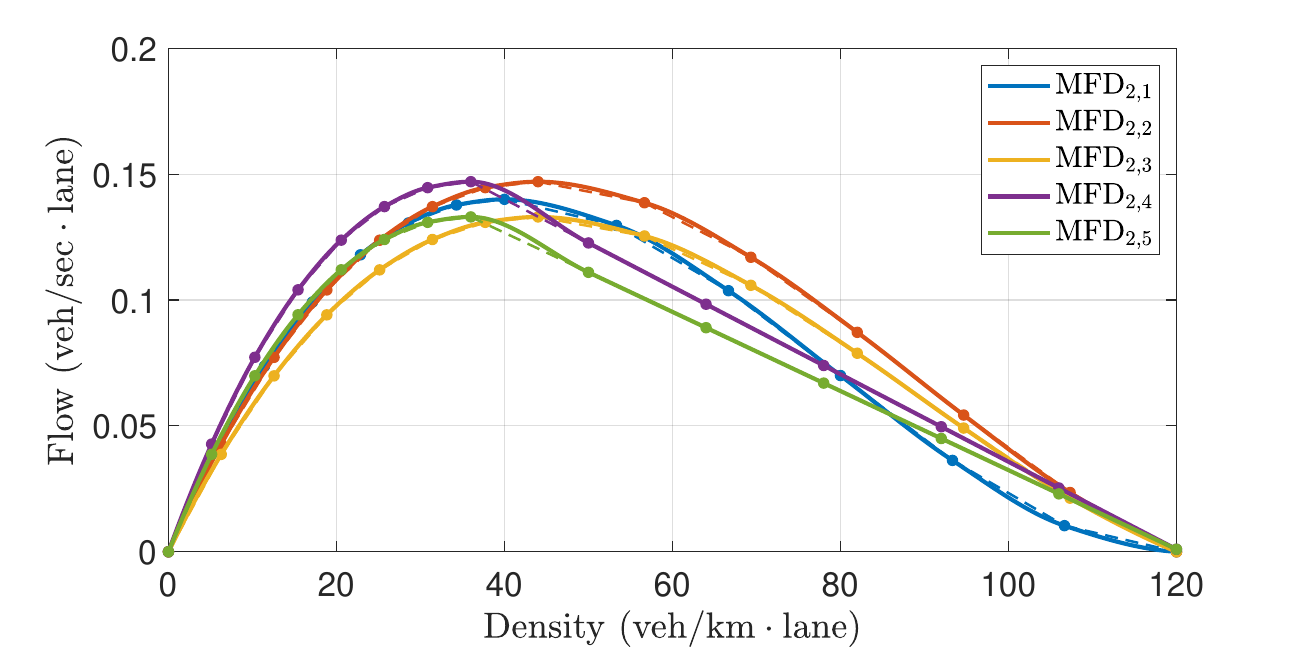}
		\caption*{(b) Real and Piece-wise approximation for noisy MFD of the periphery with $N=15$ segments.}
		\label{fig:noisyMFD_2}
	\end{minipage}
	\caption{Noisy MFD for the city center region and the periphery regions, respectively.}
	\label{fig:noisy}
\end{figure}

In this section, we evaluate the performance of our solution approach, when we add noise in the MFDs associated with every region. We consider as input the network parameters and demand distributions that resulted in the cumulative curves shown in Fig \ref{fig:cumulative_1}. The percentages of the deviations are $ \pm10\%$ and $\pm 5\%$ for the critical density $\rho_r^C$ and maximum outflow generated by the MFD, $g_r^C$, respectively. For the sake of clarity, $\text{MFD}_{r,1}, r=1,\ldots,5$ denotes the MFD without the incorporation of noise, as shown in Fig. \ref{fig:noisy}. Note that for every noisy MFD expressed as a smooth third-order polynomial in this figure, we pursue the approximation procedure described at the beginning of Section \ref{ch:simulation} to obtain the corresponding piece-wise linear segments. In Eqs. \eqref{TTS_red}-\eqref{TTS_red_sur} below, we define the TTS reduction (compared to the no-control case, i.e., $\epsilon=0$ min) attained when executing the \Coupledacronym\ and \SurrogateAcronym\ algorithms, respectively, for different values of $\epsilon$
\begin{align}
	\label{TTS_red}
	\text{TTS reduction}_1(\epsilon) &= \frac{J_{\text{TTS}}^{\text{NC}} - J_{\text{TTS}}^{\Coupledacronym}(\epsilon)}{J_{\text{TTS}}^{\text{NC}}}\times 100\%,\\
	\label{TTS_red_sur}
	\text{TTS reduction}_2(\epsilon) &= \frac{J_{\text{TTS}}^{\text{NC}} - J_{\text{TTS}}^{\SurrogateAcronym}(\epsilon)}{J_{\text{TTS}}^{\text{NC}}}\times 100\%, 
\end{align}

\noindent where the obtained TTS value for the no-control case (without any school start time change, i.e., $\epsilon=0$) is denoted with $J_{\text{TTS}}^{\text{NC}}$ (veh h), NC standing for No-Control. We reiterate that the selected values of $\epsilon$ render the execution of the \ExhaustiveAcronym\ algorithm in this part of the experiments prohibitive. To this end, for each percentage deviation pair shown in Table \ref{tab:uncertain_MFD}, we assess the effectiveness of our proposed optimization algorithms, considering three values of overall permissible change in the start time of schools, i.e., $\epsilon = \{100,200,300\}$ min.

\begin{table}[t]
	\centering
	\footnotesize
	\caption{TTS values derived from the ABO-MILP and ABO-SS Algorithms for different variants of the MFDs with a cumulative demand of $2.5 \times 10^4$ vehicles entering the network.}
	\label{tab:uncertain_MFD}
	\begin{tabular}{@{}m{0.8cm} m{3.2cm} m{1.7cm} c c c c c c@{}}
		\toprule
		Case & Deviation from $(\rho_r^C, g_r^C)$ & $J_{\text{TTS}}^{\text{NC}}$ & $\epsilon$ & $J_{\text{TTS}}^{\SurrogateAcronym}$ & $J_{\text{TTS}}^{\Coupledacronym}$ & RD$_2$ & TTS red\_1 & TTS red\_2 \\
		& (\%) & (veh h) & (min) & (veh h) & (veh h) & (\%) & (\%) & (\%) \\
		\midrule
		\multirow{3}{*}{1} & \multirow{3}{*}{Baseline}  & \multirow{3}{*}{18449} & 100 & 17178 & 16973 & 1.21 & 8.00  & 6.89\\
		&                              &                        & 200 & 17047 & 16917 & 0.77 & 8.30  & 7.60 \\
		&                              &                        & 300 & 17135 & 16145 & 6.13 & 12.49  & 7.12 \\
		\midrule
		\multirow{3}{*}{2} & \multirow{3}{*}{(+10$\%$, +5$\%$)}  & \multirow{3}{*}{16753} & 100 & 16540 & 16392 & 0.90 & 2.15 & 1.27 \\
		&                              &                        & 200 & 16141 & 15891 & 1.57 & 5.15 & 3.65 \\
		&                              &                        & 300 & 16169 & 15451 & 4.65 & 7.77 & 3.49 \\
		\midrule
		\multirow{3}{*}{3} & \multirow{3}{*}{(+10$\%$, -5$\%$)}  & \multirow{3}{*}{18677} & 100 & 18212 & 18327 & 0.63 & 1.87 & 2.49 \\
		&                              &                        & 200 & 18190 & 18098 & 0.51 & 3.10 & 2.61 \\
		&                              &                        & 300 & 17913 & 17415 & 2.86 & 6.76 & 4.09 \\
		\midrule
		\multirow{3}{*}{4} & \multirow{3}{*}{(-10$\%$, +5$\%$)}  & \multirow{3}{*}{19121} & 100 & 17266 & 17315 & 0.28 & 9.45 & 9.70 \\
		&                              &                        & 200 & 17979 & 17074 & 5.30 & 10.71 & 5.97 \\
		&                              &                        & 300 & 17313 & 16607 & 4.25 & 13.15 & 9.46 \\
		\midrule
		\multirow{3}{*}{5} & \multirow{3}{*}{(-10$\%$, -5$\%$)}  & \multirow{3}{*}{21655} & 100 & 19537 & 19462    & 0.38 & 10.13 & 9.78 \\
		&                              &                        & 200 & 19499 & 19711    & 1.09 & 8.98 & 9.96 \\
		&                              &                        & 300 & 19086 & 18898    & 0.99 & 12.73 & 11.86 \\
		\bottomrule
	\end{tabular}
\end{table}

From Table \ref{tab:uncertain_MFD} we can observe that as we decrease the values of critical density and maximum outflow of the MFD, then the TTS values observed for the no-control case, i.e.,  $J_{\text{TTS}}^{\text{NC}}$, increase significantly. In particular, for the combination ($-10\%, -5\%$), we get a $\frac{21655-18449}{21655}=$14.80$\%$ increase in the TTS value compared to the baseline case (without noise added in the MFD) due to the reduction of the critical density and maximum outflow in every regional MFD, respectively. Moreover, the TTS values derived from the \Coupledacronym\ algorithm almost always outperform the respective values obtained from the \SurrogateAcronym\ algorithm. Lastly, the TTS values, $J_{\text{TTS}}^{\Coupledacronym}$ and $J_{\text{TTS}}^{\SurrogateAcronym}$ yield on average a $7.66\%$ and $5.46\%$ reduction in TTS, respectively.

\subsection{Sensitivity Analysis of the impact of different scheduling preferences on the traffic system}

\begin{table}[t]
	\caption{TTS values for different scheduling parameters when executing the AUE algorithm giving as input the initial school start time $\tau_s, \forall s\in\mathcal{S}_b, \forall b\in\mathcal{B}$, i.e., $\epsilon = 0$ min and the fixed and flexible work start times $\bar{t}_d, \bar{t}_d^{flex}, \forall d\in\mathcal{D}$, respectively.}
	\label{tab:sensitivity_scheduling}
	\centering
	\footnotesize
	% 12 columns: 1 (Category) + 1 (Instance) + 9 params + 1 TTS
	\begin{tabular}{@{}m{2.6cm} c c c c c c c c c c c@{}}
		\toprule
		Category & Case & $\alpha^S$ & $\beta^S$ & $\gamma^S$ & $\alpha^W$ & $\beta^W$ & $\gamma^W$ & $\alpha^{W,flex}$ & $\beta^{W,flex}$ & $\gamma^{W,flex}$ & $J_{\text{TTS}}^{\text{NC}}$ (veh h) \\
		\midrule
		
		\multirow{1}{*}{Baseline}
		& -- & 8 & 5 & 16 & 14 & 7 & 20 & 14 & 7 & 20 & 18449 \\
		\midrule
		
		\multirow{4}{*}{Flex Adjustments}
		& 1 & 8  & 5 & 16 & 14 & 7 & 20 & 11 & 5  & 10 & 18104 \\
		& 2 & 8  & 5 & 16 & 14 & 7 & 20 & 16 & 8  & 30 & 18334 \\
		& 3 & 8  & 5 & 16 & 14 & 7 & 20 & 3  & 12 & 12 & 18401 \\
		& 4 & 8  & 5 & 16 & 14 & 7 & 20 & 25 & 12 & 12 & 17842 \\
		\midrule
		
		\multirow{4}{*}{Symmetric Scale Effects}
		& 5 & 6  & 3  & 12 & 8  & 5 & 16 & 8  & 5  & 16 & 18749 \\
		& 6 & 10 & 7  & 30 & 16 & 5 & 24 & 16 & 5  & 24 & 16953 \\
		& 7 & 30 & 10 & 20 & 30 & 8 & 20 & 30 & 8  & 20 & 16659 \\
		& 8 & 1  & 30 & 12 & 1  & 30 & 16 & 1  & 30 & 16 & 20015 \\
		\midrule
		
		\multirow{4}{*}{Mixed Effects}
		& 9  & 4  & 15 & 6  & 20 & 30 & 6  & 1  & 9  & 30 & 19623 \\
		& 10 & 15 & 3  & 30 & 3  & 10 & 30 & 30 & 10 & 3  & 17743 \\
		& 11 & 30 & 1  & 12 & 1  & 12 & 30 & 12 & 30 & 1  & 17722 \\
		& 12 & 12 & 6  & 1  & 3  & 1  & 12 & 30 & 1  & 12 & 16570 \\
		\bottomrule
	\end{tabular}
\end{table}

\begin{table}[t]
	\caption{TTS values derived from the ABO-MILP and ABO-SS Algorithms for the scheduling parameters stemming from Cases 8 and 9 of Table \ref{tab:sensitivity_scheduling}.}
	\label{tab:scheduling_values}
	\centering
	\footnotesize
	\begin{tabular}{@{}m{1.5cm} m{1.7cm} c c c c c c@{}}
		\toprule
		Case & $J_{\text{TTS}}^{\text{NC}}$ & $\epsilon$ & $J_{\text{TTS}}^{\SurrogateAcronym}$ & $J_{\text{TTS}}^{\Coupledacronym}$ & RD$_2$ & TTS red\_1 & TTS red\_2 \\ 
		& (veh h) & (min) & (veh h) & (veh h) & (\%) & (\%) & (\%) \\ \midrule
		
		\multirow{3}{*}{Baseline} 
		& \multirow{3}{*}{18449} & 100 & 17178 & 16973 & 1.21 & 8.12 & 6.89 \\
		&                         & 200 & 17047 & 16917 & 0.77 & 8.30 & 7.60 \\
		&                         & 300 & 17135 & 16145 & 6.13 & 12.49 & 7.12 \\ \midrule
		
		\multirow{3}{*}{8} 
		& \multirow{3}{*}{20015} & 100 & 18236 & 18714 & 2.36 & 6.50 & 8.89 \\
		&                         & 200 & 18064 & 17571 & 2.81 & 12.21 & 9.75 \\
		&                         & 300 & 17555 & 16940 & 3.63 & 15.36 & 12.29 \\ \midrule
		
		\multirow{3}{*}{9} 
		& \multirow{3}{*}{19623} & 100 & 17618 & 18696 & 6.12 & 4.72 & 10.22 \\
		&                         & 200 & 17610 & 17480 & 0.74 & 10.92 & 10.26 \\
		&                         & 300 & 16761 & 16453 & 1.87 & 16.15 & 14.58 \\
		
		\bottomrule
	\end{tabular}
\end{table}

In this section, we perform a sensitivity analysis to assess the influence of commuters’ scheduling preferences on the overall traffic performance. The analysis focuses on variations of the $\alpha$-$\beta$-$\gamma$ parameters that characterize the trade-offs between travel time, early arrival and late arrival, respectively, for each class of commuters. By systematically modifying these parameters, in Table \ref{tab:sensitivity_scheduling} we show a set of distinct scenarios that represent different behavioral profiles, ranging from highly schedule-constrained commuters (high $\beta$ values) to more flexible ones (higher $\gamma$ values). Each scenario aims to capture how heterogeneity in departure-time preferences affects the resulting network equilibrium and the TTS. Adjustments to flexible work preferences (Cases~1-4) indicate that increasing sensitivity to travel time and schedule variability can reduce TTS by shifting departures to less congested periods. Scaling the scheduling values symmetrically across all trip types (Cases~5-8) shows a nonlinear effect: moderate increases in these values can improve network performance by smoothing demand, while higher values of penalties for early or late arrivals can exacerbate congestion. Mixed scaling scenarios (Cases~9-12) demonstrate that carefully combining heterogeneous preferences across trip types can further reduce TTS (e.g., Case~12 achieves 16570~veh$\cdot$h), whereas uncoordinated heterogeneity may increase congestion. Overall, the findings shown in this table, indicate that the departure patterns for school trips are largely constrained by the fixed schedule. While variations in flexible work trip preferences slightly influence total system travel time, the analysis confirms that, under the current school schedule (no school start change), school trip parameters set the reference conditions for network equilibrium.\par

Based on the sensitivity analysis, we selected the two cases from Table \ref{tab:sensitivity_scheduling} that resulted in the highest values of TTS as reference points to evaluate the effectiveness of our proposed solution methodology. By focusing on these scenarios, which represent the least favorable network conditions under the baseline school start time schedule ($\epsilon = 0$), we plan to assess how the proposed adjustments to school departure patterns can improve traffic dynamics and reduce congestion. To this end, we consider the Cases 8 and 9 for further investigation.

First, in the baseline case (see Table \ref{tab:scheduling_values}), $J_{\text{TTS}}^{\Coupledacronym}$ consistently outperforms $J_{\text{TTS}}^{\SurrogateAcronym}$ across all values of $\epsilon$, with the performance gap widening as $\epsilon$ increases (most pronounced at $\epsilon=300$). This indicates that the \Coupledacronym\ Algorithm is particularly effective when the scheduling flexibility is higher. While the \SurrogateAcronym\ Algorithm provides competitive or even superior solutions for low values of $\epsilon$ (notably in Cases 8 and 9), the \Coupledacronym\ Algorithm becomes dominant as scheduling flexibility increases, offering both higher efficiency and greater TTS reductions. Another interesting finding is that for both Cases 8 and 9, RD$_2$ decreases as $\epsilon$ increases, which shows that the discrepancy between $J_{\text{TTS}}^{\SurrogateAcronym}$ and $J_{\text{TTS}}^{\Coupledacronym}$ narrows at higher flexibility levels, albeit with some mild fluctuations. Lastly, in both approaches for both cases, we can observe that greater scheduling flexibility consistently leads to larger total time spent savings. 

It is also worth commenting on the scheduling parameters associated with Cases 8 and 9 shown in Table~\ref{tab:sensitivity_scheduling}. 
Case~8 corresponds to extreme weights where both schools and workplaces are highly penalized for early arrivals ($\beta^S=\beta^W=30$) while the penalties for travel time and lateness remain comparatively low. 
This configuration explains why the solutions are sensitive to the choice of optimization framework: the surrogate solver $J_{\text{TTS}}^{\SurrogateAcronym}$ performs better for $\epsilon=100$, but as flexibility grows ($\epsilon>100$), the relaxed optimization program $J_{\text{TTS}}^{\Coupledacronym}$ better balances the strong earliness penalties, leading to superior outcomes. On the other hand, Case~9 represents a more mixed configuration where schools face relatively high lateness penalties ($\gamma^S=6$) and workplaces face very large penalties on both earliness ($\beta^W=30$) and lateness ($\gamma^W=6$), while the flexible workplace parameters put additional emphasis on lateness ($\gamma^{W,flex}=30$). 
This setup creates competing objectives that again favor the surrogate solver at small $\epsilon$ (where flexibility is limited and the solver can quickly adjust to extreme penalties), but as $\epsilon$ increases, the MILP-based optimization algorithm can retrieve a solution of higher quality.  

\begin{figure}[h!t]
	
	\begin{subfigure}{0.48\linewidth}
		\includegraphics[width=\linewidth]{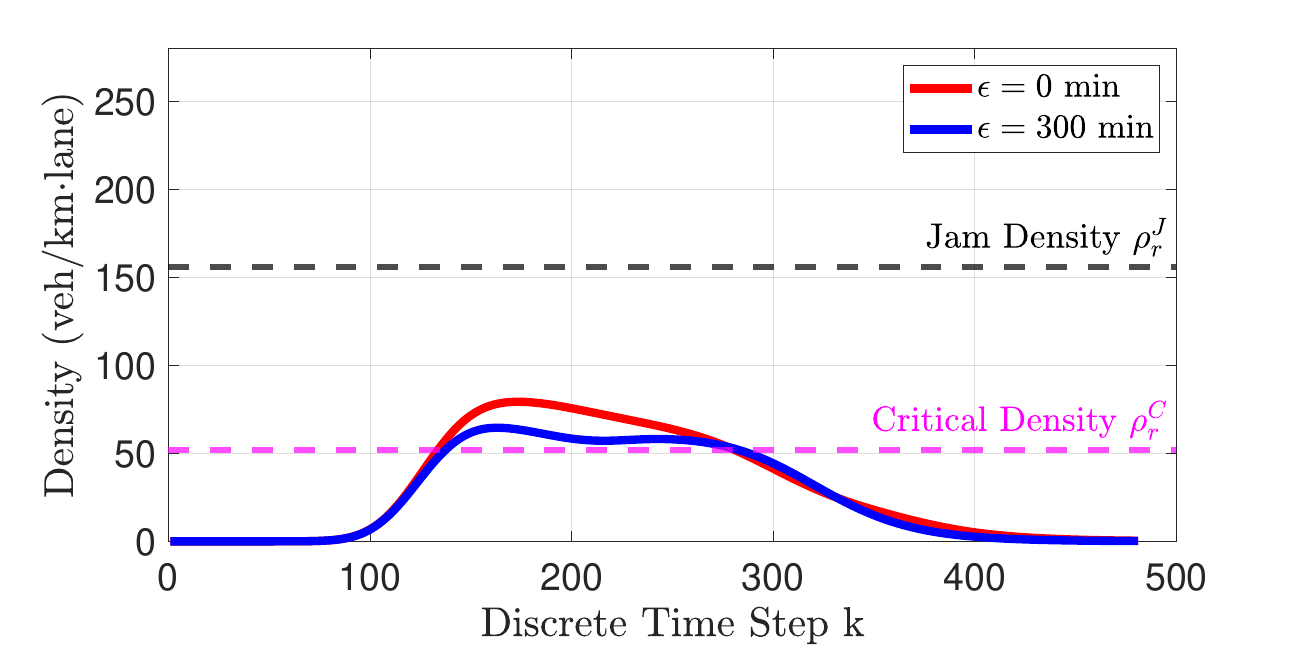}
		\caption*{(a) Density evolution in Region 1.}
		\label{fig:density_UE_shadow_1}
	\end{subfigure}
	\hfill
	\begin{subfigure}{0.48\linewidth}
		\includegraphics[width=\linewidth]{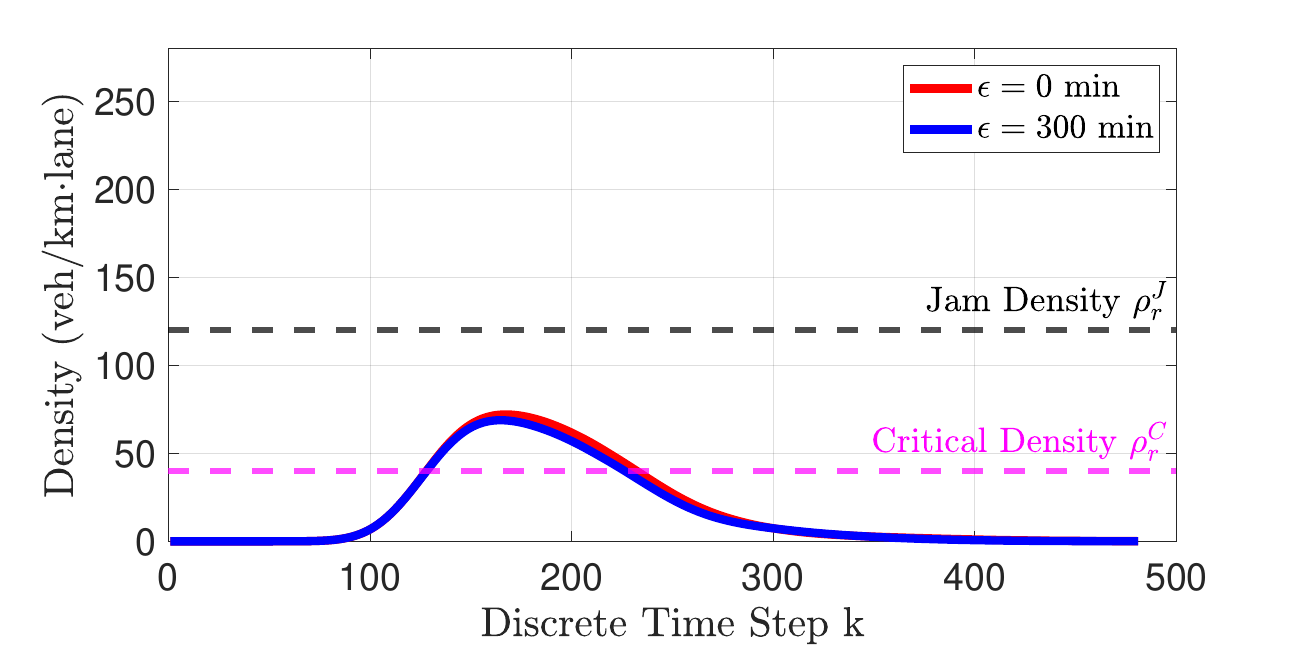}
		\caption*{(b) Density evolution in Region 2.}
		\label{fig:density_UE_shadow_2}
	\end{subfigure}%
	\hfill
	
	\vspace{0.3cm}
	
	\begin{subfigure}{0.48\linewidth}
		\includegraphics[width=\linewidth]{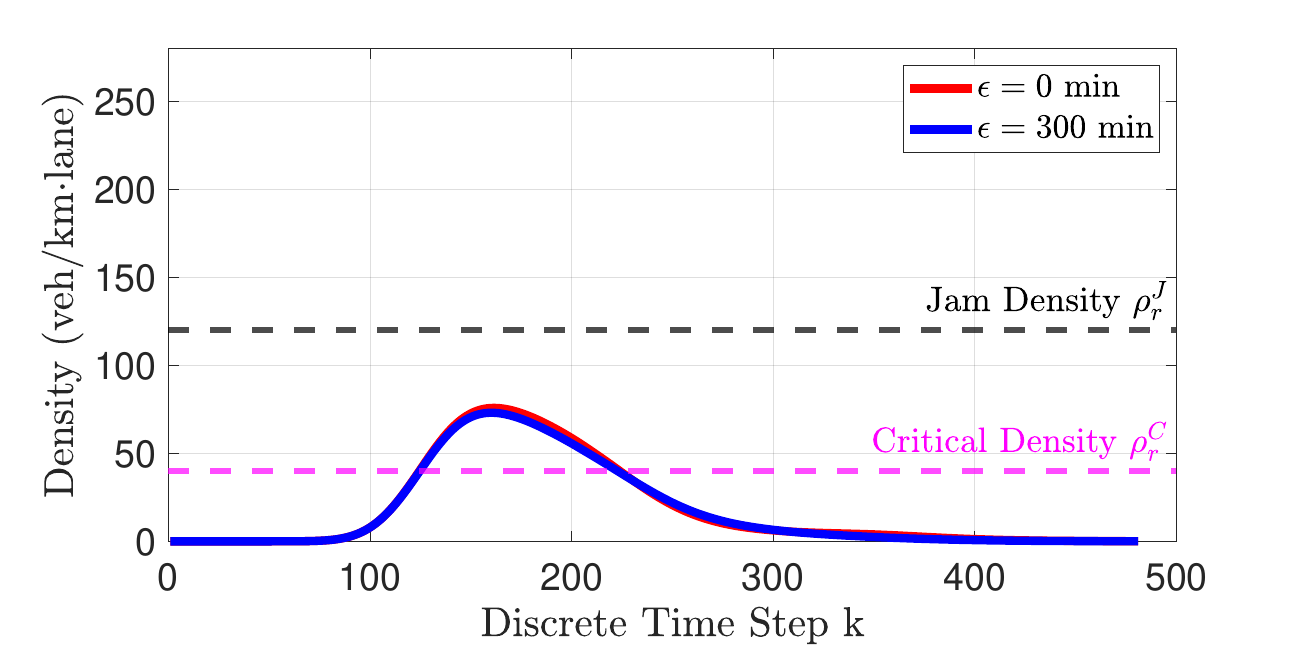}
		\caption*{(c) Density evolution in Region 3.}
		\label{fig:density_UE_shadow_3}
	\end{subfigure}
	\hfill
	\begin{subfigure}{0.48\linewidth}
		\includegraphics[width=\linewidth]{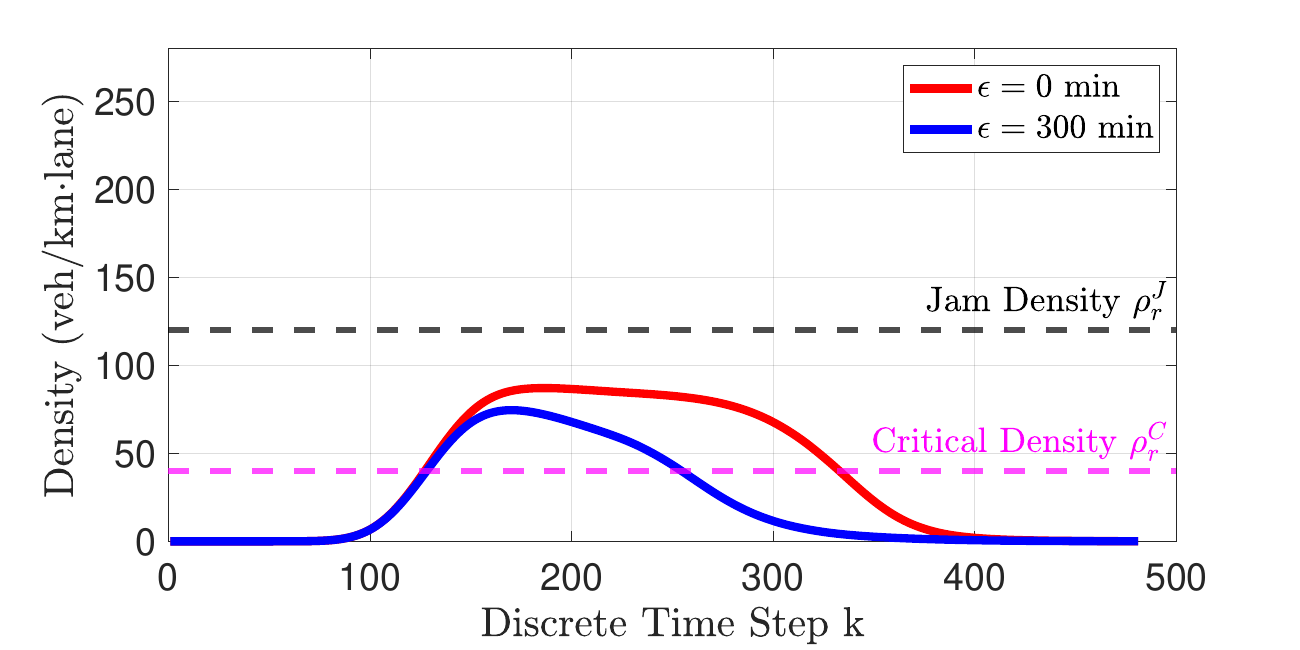}
		\caption*{(d) Density evolution in Region 4.}
		\label{fig:density_UE_shadow_4}
	\end{subfigure}
	
	\vspace{0.3cm}
	
	\begin{subfigure}{\linewidth}
		\centering
		\includegraphics[width=0.48\linewidth]{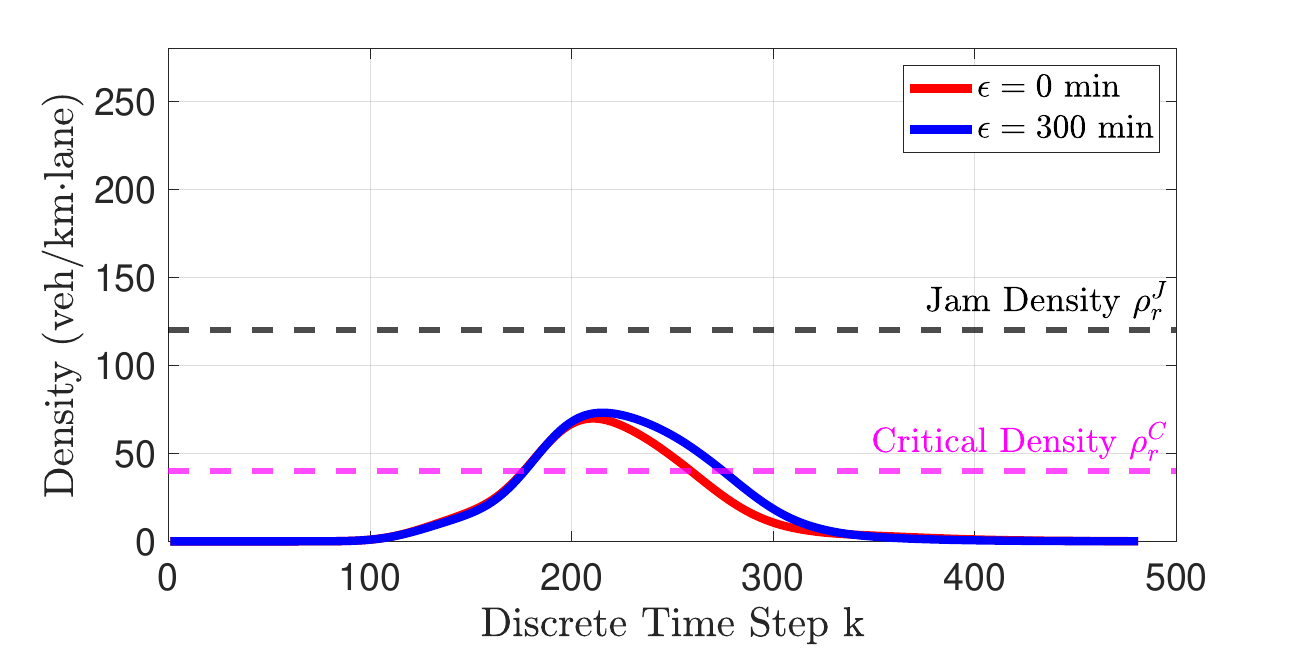}
		
		\vspace{0.2em}
		
		{\centering \caption*{(e) Density evolution in Region 5.}}
		
		\label{fig:density_UE_shadow_5}
	\end{subfigure}
	
	\caption{Density evolution when executing ABO-MILP algorithm for values ($\epsilon=0$) and ($\epsilon=300$) min, respectively for the scheduling preferences of commuters with respect to Case 8 of Table \ref{tab:sensitivity_scheduling}.}
	\label{fig:dens_evolf1_new_shadow}
\end{figure}

Overall, the comparison shows that the relative performance of the two solution frameworks is not only a function of $\epsilon$, but is also strongly influenced by the underlying trade-offs encoded in the scheduling preferences of each class of commuters $(\alpha, \beta, \gamma)$. Fig. \ref{fig:dens_evolf1_new_shadow} displays the evolution of density for each region in the network for the cases where i) no school start time change takes place ($\epsilon=0$) and ii) when the permissible overall school start time change is set to $\epsilon=300$ minutes with respect to the scheduling parameters of Case 8 illustrated in Table \ref{tab:sensitivity_scheduling}, stemming from the \Coupledacronym\ Algorithm.

\section{Conclusions}
\label{ch:conclusions}

This paper presents a methodology for the coordinated selection of school start times in multi-region urban networks, explicitly accounting for traffic dynamics as described by macroscopic fundamental diagrams. The framework bridges the departure time choice behavior of both school and work commuters with the school start times and the fixed and flexible work start times, respectively. It captures the interaction between school-related and work-related trips, emphasizing that congestion emerges from the combination of multiple demand classes rather than from school demand alone. Within this context, a nonconvex optimization problem is formulated to retrieve school start times that mitigate congestion while balancing travel efficiency and schedule delay costs.

The proposed solution approach replaces nonlinear constraints with convex approximate ones, enabling the acquisition of both lower and upper-bounds, while providing guarantees on optimality gaps. Benchmarking against the \ExhaustiveAcronym\ Algorithm demonstrates that our methodology can achieve near-optimal solutions at a fraction of the computational cost for the cases where the \ExhaustiveAcronym\ algorithm is tractable.

A key methodological contribution of this work is the development of a tailored solution approach for the resulting bilevel optimization problem. The Upper-Level and Lower-Level problems are tightly coupled: school start time decisions influence commuters’ departure-time choices, while the resulting equilibrium demand patterns determine the traffic dynamics that shape the Upper-Level objective. To capture this mutual feedback, we propose an iterative framework that alternates between the two levels. For a given set of school start times, the Lower-Level module computes the equilibrium departure-time profiles of both commuter classes through a deterministic dynamic user equilibrium algorithm. The resulting demand profiles are then supplied to the Upper-Level module, which updates the school start times by solving the corresponding optimization problem using one of the proposed solution approaches. This alternating procedure enables consistent coordination between scheduling decisions and behavioral responses, while maintaining computational tractability despite the nonlinear and nonconvex nature of the original bilevel formulation.

The experiments show that the system performance is strongly influenced by the choice of scheduling preferences of commuters ($\alpha$-$\beta$-$\gamma$) and the permissible school start time change, $\epsilon$. Notably, the comparative effectiveness of optimization methods varies with these operational conditions: while the lower-bound stemming from the relaxed optimization framework typically outperforms the corresponding solution derived from the surrogate solver, specific parameter combinations and low $\epsilon$ values can lead to superior performance from the derivative-free approach. These results highlight the importance of jointly considering behavioral preferences, scheduling flexibility, and network congestion when designing policy interventions.

From a practical standpoint, the framework provides transportation authorities with a structured tool to explore trade-offs between total time spent and allowable shifts in school start times. The Pareto front derived in the experiments illustrate how different $\epsilon$ choices influence network performance, offering actionable guidance for policy decisions. Coordination between the school start times and the departure time choice aspect of both classes of commuters, coupled with MFD-informed traffic dynamics, consistently reduces peak congestion and improves network resilience.

Overall, the study confirms that carefully designing adjustments to school start times can mitigate congestion and enhance network efficiency. Nevertheless, the proposed framework relies on several modeling assumptions. In particular, route choice behavior is not modeled explicitly, and travelers are assumed to respond to congestion solely through adjustments in their departure times. Moreover, the MFD-based regional representation abstracts from detailed network topology and signal control operations. Future research could extend the framework by incorporating endogenous route choice behavior and stochastic travel times. In addition, further validation of the proposed methodology will involve integrating a microsimulation environment for a real-world traffic network using traffic data from heterogeneous sources, thereby providing empirical evidence of the potential effectiveness of the proposed approach in practical urban traffic settings.

\section{Acknowledgements}

This work is supported by the European Union (i. ERC, URANUS, No. 101088124, and ii. Horizon 2020 Teaming, KIOS CoE, No. 739551), and the Government of the Republic of Cyprus through the Deputy Ministry of Research, Innovation, and Digital Strategy. Views and opinions expressed are however those of the author(s) only and do not necessarily reflect those of the European Union or the European Research Council Executive Agency. Neither the European Union nor the granting authority can be held responsible for them. Antonios Georgantas is partially supported by the Foundation for Education and European Culture, Founders Nicos \& Lydia Tricha, Greece.

\begin{appendices}
	
	\section{Alternating Bilevel Optimization - Exhaustive Search (ABO-ES) Algorithm}
	\label{sec:es_contribution}
	
	This section proposes the \Exhaustive\ algorithm, termed \ExhaustiveAcronym\, for determining the optimal vector of school start times that leads to the acquisition of the smallest value of the Total Time Spent metric for an appropriate discretization step $V$ and threshold $\epsilon$. The algorithm presented in this section will serve as a benchmark to compare against our main solution approach (see Sections \ref{sec:reformulation} - \ref{sec:pareto}). The \ExhaustiveAcronym\ algorithm, outlined in Algorithm \ref{alg:serqt}, is an iterative procedure that evaluates the impact of all the feasible school start time combinations to the TTS metric, ensuring that the overall change in the school start times does not exceed a threshold $\epsilon$ (min). \par

	\begin{algorithm}[t]
		\algsetup{linenosize=\tiny}
		\footnotesize
		%{\fontsize{6}{6}\selectfont
			\begin{algorithmic}[1]
				\STATE \textbf{Input}: $\tau_s, \forall s\in\mathcal{S}_b, \forall b\in\mathcal{B}, \bar{t}_d, \bar{t}_d^{flex}, \forall d\in\mathcal{D}, V, K,~ \mathcal{E}, \zeta$.
				\STATE \quad \quad \quad Traffic Network Parameters, External Demand same as Algorithm \ref{alg:seq1gsg}.		
				%			\STATE  \quad \quad \quad Current State: $\rho_r(k),~ \rho_r^S(k),~ \rho_r^W(k),~r\in\mathcal{R},~k\in\mathcal{K},~ \rho_{o,r,s}^S(k),~ \rho_{o,r,d}^W(k),~ r\in\mathcal{R},~ o\in\mathcal{O},~ s\in\mathcal{S}_b, b\in\mathcal{B}, d\in\mathcal{D}, k\in\mathcal{K}$.
				%			\STATE \quad \quad \quad  Traffic Network Parameters: $u_r^f, \rho_r^C, \rho_r^J, q_r^C, g_r^C, L_r, l_r, r\in\mathcal{R}, C_{r,j}^{\text{MAX}}, \alpha_{r,j}, r\in\mathcal{R}, j\in\mathcal{J}_r^-$, $\gamma_{r,j,d}(k), r\in\mathcal{R},  j\in\mathcal{J}_r^-,  k\in\mathcal{K}, d\in\mathcal{D}, \gamma_{r,j,b}(k), r\in\mathcal{R},  j\in\mathcal{J}_r^-,  k\in\mathcal{K}, b\in\mathcal{B}$.
				%			\STATE \quad \quad \quad External Demand: $d_{o,s}^S(k),  \forall o\in\mathcal{O}, \forall b\in\mathcal{B}, \forall s\in\mathcal{S}_b$, $\forall k\in\mathcal{K}, d_{o,d}^W(k), \forall o\in\mathcal{O}, \forall d\in\mathcal{D}, \forall k\in\mathcal{K}, d_{o,s,m}^S(k), \forall o\in\mathcal{O},$\\
				%			\quad \quad \quad $ \forall b\in\mathcal{B}, \forall s\in\mathcal{S}_b,
				%			\forall m\in\mathcal{M}, \forall k\in\mathcal{K}$.
				\FOR{$\kappa=0: |\mathcal{E}|/V$}	
				\STATE Generate set $\mathcal{T}$ of all school start time vector combinations $\boldmath{\tau}_{\omega}^{\ExhaustiveAcronym}$ =  $[\tilde{\tau}_1^{\ExhaustiveAcronym}, \tilde{\tau}_2^{\ExhaustiveAcronym},\ldots,\tilde{\tau}_{\zeta}^{\ExhaustiveAcronym}]$ such that $J_{\text{STC}}\leq \epsilon = \kappa V$ using the formula shown in Eq. \eqref{pkjhgh_bidirectional}.
				\FOR{$\boldmath{\tau}_{\omega}^{\ExhaustiveAcronym}$ $\in\mathcal{T}$}
				\STATE Execute Algorithm \ref{alg:user_equilibrium} when school start time vector $\boldmath{\tau}_{\omega}^{\ExhaustiveAcronym}$ and demand vectors $\mathbf{d}_{o,s}^S, \mathbf{d}_{o,d}^W$ are given as input.
				\STATE \textbf{Output}: TTS and STC values for the iteration $\Iota$ that achieves convergence of Algorithm \ref{alg:user_equilibrium}, $J_{\text{TTS}}^{\ExhaustiveAcronym,I}(\omega)$ and $J_{\text{STC}}^{\ExhaustiveAcronym,I}(\omega)$ using Eqs. \eqref{fixed} and \eqref{QUO1}, respectively.
				%\STATE Simulate traffic system utilizing Eqs. \eqref{shift_lhhg} - \eqref{QUO1} when vector $\boldmath{\tilde{\tau}}^{\text{ES}}$ is given as input.
				\ENDFOR
				\STATE Select $\boldmath{\tau}_{\kappa}^{\ExhaustiveAcronym,*}$ that yields the smallest $J_{\text{TTS}}^{\ExhaustiveAcronym, I}(\omega)$ value among all school start time combinations. 
				\STATE \textbf{Output}: $J_{\text{TTS}}^{\ExhaustiveAcronym,*}(\kappa), J_{\text{STC}}^{\ExhaustiveAcronym,*}(\kappa)$ based on vector $\boldmath{\tau}_{\kappa}^{\ExhaustiveAcronym,*}$ obtained from Line 11.
				\ENDFOR
				
			\end{algorithmic}
			%	}
		\caption{\Exhaustive\ (\ExhaustiveAcronym) Algorithm}
		\label{alg:serqt}
	\end{algorithm}
	
	First, Algorithm \ref{alg:serqt} takes as input the start time of each school, $\tau_s$, a set of values of $\epsilon$ stored in set $\mathcal{E}$ and $\zeta$ as the summation of the cardinality of every set $\mathcal{S}_b$, i.e., $\zeta = \sum_{b\in\mathcal{B}}|\mathcal{S}_b|$. Note that the discretization step $V$ is the same for threshold $\epsilon$ and for the school start time $\tau_s$, respectively. Furthermore, the algorithm uses the same traffic network parameters and external demand with those given in Algorithm \ref{alg:seq1gsg} (Lines 1-2). \par
	
	Later, we generate set $\mathcal{T}$ that contains all the school start time combinations for which the summation of all the shifted school start times does not exceed threshold $\epsilon$, $J_{\text{STC}}\leq \epsilon = \kappa V$ (Line 4). It is worth reiterating that in the developed framework we allow for both forward and backward shifts in the school start time.
	Then in the inner-loop iterative procedure of the algorithm, we execute Algorithm \ref{alg:user_equilibrium} giving as input the vector of school start times, $\boldmath{\tau}_{\omega}^{\ExhaustiveAcronym}$ and the demand vectors $\mathbf{d}_{o,s}^S, \mathbf{d}_{o,d}^W$ (Lines 5-8). During each inner-loop iteration we obtain as output the TTS and STC values associated with the vector $\boldmath{\tau}_{\omega}^{\ExhaustiveAcronym}$, respectively (Line 7). 
	Then we choose the vector $\boldmath{\tau}_{\kappa}^{\ExhaustiveAcronym,*}$ that yields the smallest value of TTS among all the school start time vector combinations (considered in the inner-loop iterative procedure) (Line 9). Thus, Algorithm \ref{alg:serqt} provides as output during each outer-loop iteration, the optimal value of TTS, $J_{\text{TTS}}^{\ExhaustiveAcronym,*}(\kappa)$ and the associated value of STC metric, $J_{\text{STC}}^{\ExhaustiveAcronym,*}(\kappa)$ (Line 10). In the next outer-loop iteration, the permissible overall change in the start time of schools increases by $V$ minutes, since $\epsilon= \kappa V$, and the above-described process is repeated to obtain the $J_{\text{STC}}^{\ExhaustiveAcronym,*}(\kappa)$ value for all $\kappa$ of interest.\\
	
	\noindent\underline{Derivation of the number of school start time vectors $\hat{P}(\kappa)$ within set $\mathcal{T}$ under bidirectional shifts}\\
	
	Consider the variable $x_{s}\in\mathbb{Z},~ s\in\mathcal{S}_b, b\in\mathcal{B}$ standing for the number of discrete intervals that we shift the start time of school $s$, such that $-M \leq x_s \leq M, M\geq 0$.  We make the convention that when $x_{s}=0, \forall s\in\mathcal{S}_b, \forall b\in\mathcal{B}$, then the overall change in the school start times is equal to zero, while when $x_{s}= M$, then the school start time takes the value, $\tau_{s} + MV$ and when $x_{s}= -M$, the school start time takes the value, $\tau_{s} - MV$. In the case that we do not allow any overall change in the start time of the schools, i.e., $\epsilon=0$, then we only get one valid combination, it being the pair containing the initial start time of each school $s$, i.e., $x_{s} = 0, s\in\mathcal{S}_b, b\in\mathcal{B}$ with $\boldmath{\tau}_{\epsilon}^{\ExhaustiveAcronym,*} = $ $[\tau_1, \tau_2,\ldots,\tau_{\zeta}]$.\par
	
	We seek the total number $\Breve{P}(\kappa)$ of integer solution vectors $\mathbf{x} = [x_{1}, x_{2}, \ldots, x_{\zeta}]$ for a specific discretization step $V$ that satisfies the condition
	\begin{equation}
		\label{bidirectional}
		J_{\text{STC}} = \epsilon=\kappa V \Rightarrow \sum_{s=1}^{\zeta} |x_s| = \epsilon, \quad \kappa \in  \Big\{0, 1, \ldots, \frac{\epsilon^{\text{MAX}}}{V} \Big\}.
	\end{equation}
	
	\noindent We define the absolute shift vector $\mathbf{y} = [y_1, y_2, \ldots, y_\zeta]$ where $y_s = |x_s| \in [0, M]$. Hence, the condition in \eqref{bidirectional} can be rewritten as
	\begin{equation}
		\label{bidirectional2}
		\sum_{s=1}^{\zeta} |x_s| = \sum_{s=1}^{\zeta} y_s = \epsilon = \kappa V.
	\end{equation}
	
	\noindent Each vector $\mathbf{y} \in \{1, \ldots, M\}^{\zeta}$ with $\sum_{s=1}^{\zeta} y_s = \kappa V$ corresponds to the number of both signed vectors of $\mathbf{x}$, it being equal to $2^{\digamma(\mathbf{y})}$, where $\digamma(\mathbf{y})$ denotes the number of non-zero entries in vector $\mathbf{y}$ \citep{guichard2017introduction}. To count them, we proceed by summing over the number $\sigma$ of non-zero components:
	
	\begin{itemize}
		\item $\binom{\zeta}{\sigma}$ ways to choose which $\sigma$ positions in $\mathbf{y}$ are non-zero.
		\item For each such $\sigma$-subset, the number of integer compositions of $\kappa$ into $\sigma$ positive integers $y_s \in [0, M]$ is given by the inclusion–exclusion formula:
		\begin{equation}
			\label{bidirectional3}
			V_{\sigma}(\kappa, M) = \sum_{s=1}^{\sigma} (-1)^s \binom{\sigma}{s} \binom{\kappa - s(M+1) - 1}{\sigma - 1}.
		\end{equation}
		\item Each such composition has $2^{\sigma}$ possible signed versions for $\mathbf{x}$.
	\end{itemize}
	
	\noindent Hence, the total number $\Breve{P}(\kappa)$ of feasible school start time vectors that satisfy condition \eqref{bidirectional2} is
	\begin{align}
		\label{eq:prune_bidirectional}
		\Breve{P}(\kappa) = \sum_{\sigma=1}^{\min(\zeta, \kappa)} \binom{\zeta}{\sigma} \cdot 2^{\sigma} V_{\sigma}(\kappa, M), 
	\end{align}
	
	\noindent where the term $\binom{n}{c}$ denotes the binomial coefficient on $n$ and $c$; we consider that $\binom{n}{c} = 0$ if $n<c$.
	
	The analysis above indicates that the total number of schools start times combinations when $J_{\text{STC}} = \epsilon=\kappa V$ is equal to $\Breve{P}(\kappa)$. To derive the total number of schools start times combinations when $J_{\text{STC}}\leq  \epsilon=\kappa V$, $\hat{P}(\kappa)$, we need to simply add the terms $\hat{P}(0), \hat{P}(1), \ldots, \hat{P}(\kappa)$, yielding
	\begin{equation}
		\label{pkjhgh_bidirectional}
		\hat{P}(\kappa) = \sum_{\kappa'\in\{0,1,\ldots, \kappa\}}\Breve{P}(\kappa').
	\end{equation}
	
	\noindent Therefore, the total number of different school start time vectors that are included within set $\mathcal{T}$ for each iteration of the \ExhaustiveAcronym\ Algorithm (Line 6) is equal to $\hat{P}(\kappa)$, such that $\epsilon=\kappa V$, retrieved from Eq. \eqref{pkjhgh_bidirectional}.
	
	\section{Solution to problem $P_1$ via a derivative-free Surrogate Solver (SS)-based approach}
	\label{sec:derivative_free}
	
	\begin{figure}[t]
		\centering
		\includegraphics[width=1.0\columnwidth]{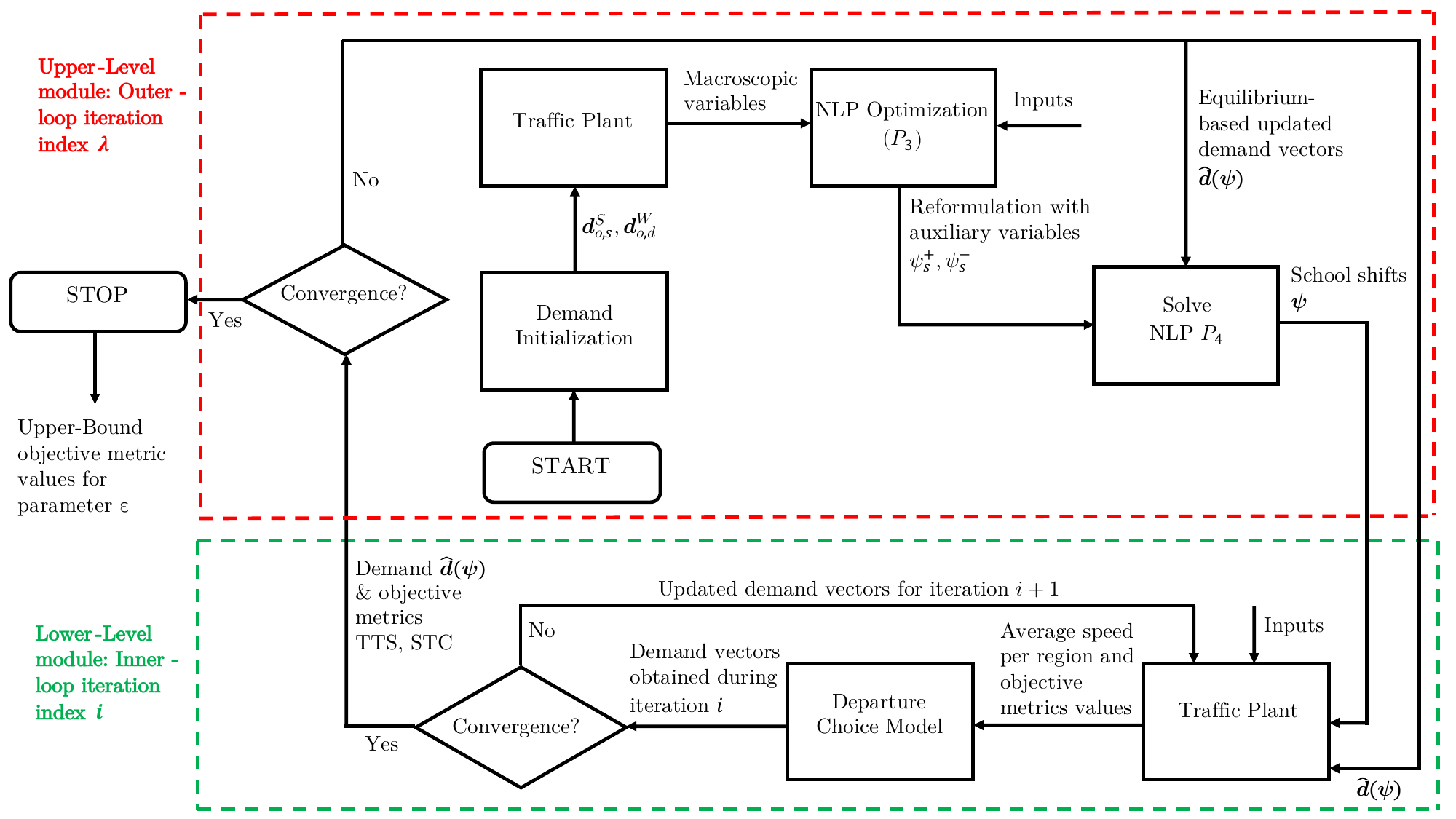} 
		\caption{Flowchart of the entire school start time selection solution approach for a selected value of overall permissible school start time change $\epsilon$, based on Fig. \ref{fig:flowchart} that is adapted for the derivative-free optimization architecture.}
		\label{fig:flowchart_continuous_case}
	\end{figure}
	
	This section presents a methodology to obtain a solution (i.e., school start time vector and corresponding objective values for TTS and STC) with respect to the Upper-Level optimization problem $P_1$ introduced in Section \ref{sec:solution}. The proposed solution is obtained through the development of the \Surrogate\ (\SurrogateAcronym) algorithm. The overall structure of the alternating bilevel solution methodology is preserved as depicted in Fig. \ref{fig:flowchart}: the Lower-Level module remains identical to that described in Section~\ref{sec:lower_level_controller}, while the Upper-Level module is adapted to accommodate the derivative-free optimization architecture. For clarity, the complete solution methodology presented in this section is summarized in Fig.~\ref{fig:flowchart_continuous_case}.
	
	The combinatorial complexity introduced by the binary decision variables, coupled with the nonlinear and nonconvex nature of the traffic dynamics with respect to problem $P_1$, substantially undermine the performance of standard nonlinear solvers, rendering it difficult to achieve a good quality solution in our setup \citep{bestuzheva2023global}. Conversely, derivative-free optimization solvers are inherently more robust in such highly discontinuous, and mixed-variable contexts. Indeed, these approaches have demonstrated favorable performance on challenging nonconvex MINLPs \citep{rios2013derivative}. Thus, we resort to a derivative-free solver as a practical benchmark against our MINLP-based formulation (problem $P_1$).\par
	
	The formulation pursued here is structurally identical with the one followed in Section~\ref{sec:reformulation} for problem $P_1$, with one key difference. The binary decision variables $\xi_{m,s}$, which explicitly determine the school start times in $P_1$ are replaced by the continuous decision variables $\psi_s\in\mathbb{R}$, $s\in\mathcal{S}_b, b\in\mathcal{B}$. The variable $\psi_s$ represents a continuous temporal shift of the start time of school $s$, expressed in units of simulation time steps. Hence, $\psi_s$ is not restricted to integer values and allows the optimization problem to explore 
	the time axis (in terms of shifting) continuously. Accordingly, the expression $\xi_{m,s}|m|V$ in \eqref{QUO1} is replaced by $\psi_s$ and the $\epsilon$-constraint in \eqref{epsi} is replaced by
	\begin{equation}
		\label{bint}
		\sum_{b\in\mathcal{B}}\sum_{s\in\mathcal{S}_b}|\psi_s| \leq \epsilon.
	\end{equation}
	
	\subsection{Continuous demand projection mechanism}
	
	Let $\boldsymbol{\psi}$ denote the Upper-Level decision vector that contains the continuous school start time shifts $\psi_s$. Because $\psi_s$ is allowed to take real values, the demand projection mechanism must operate in continuous time and translate the baseline equilibrium demand accordingly. In the continuous case, the demand projection mechanism is expressed as a set of functional constraints that link the Upper-Level decision vector $\boldsymbol{\psi}$ to the time-dependent demand entering problem $P_3$, i.e., $\hat{\mathbf{d}}(\boldsymbol{\psi})$. The behavioural rationale and modeling assumptions are identical to those introduced in Section~\ref{sec:demand_projection}; the present section provides the continuous counterpart of that mechanism.
	
	In analogy with the discrete case, the school-related demand
	entering the network is required to correspond to a temporal translation of the baseline equilibrium demand
	$\hat{d}^{S}_{o,s}(k)$, obtained as output from Algorithm \ref{alg:user_equilibrium}. Formally, for each origin $o$, school $s$, and time step $k\in\mathcal{K}$, the school-related demand entering the network under the continuous shift decision 
	$\psi_s$ is defined as
	\begin{equation}
		\label{shift_continuous}
		d_{o,s}^{S}(k;\boldsymbol{\psi}) = d_{o,s,\psi_s}^{S}(k),
	\end{equation}
	
	\noindent where the translated demand profile is given by
	\begin{equation}
		\label{mapping_continuous}
		d_{o,s,\psi_s}^{S}(k)=
		\begin{cases}
			\hat{d}^{S}_{o,s}(k-\psi_s), & \psi_s>0,\; \psi_s < k \le K,\\[3pt]
			0, & \psi_s>0,\; 1 \le k \le \psi_s,\\[3pt]
			\hat{d}^{S}_{o,s}(k), & \psi_s = 0,\\[3pt]
			\hat{d}^{S}_{o,s}(k-\psi_s), & \psi_s<0,\; 1 \le k \le K+\psi_s,\\[3pt]
			0, & \psi_s<0,\; K+\psi_s < k \le K.
		\end{cases}
	\end{equation}
	
	\noindent Similarly to the discrete formulation, school start time decisions are assumed to affect only commuters of class S. Hence, the demand associated with work-related commuters (class W) is treated as exogenous and fixed to the equilibrium-based demand obtained from Algorithm~\ref{alg:user_equilibrium}. Hence,
	\begin{equation}
		\label{projection_W_cont}
		d_{o,d}^{W}(k;\boldsymbol{\psi}) = \hat{d}_{o,d}^{W}(k),
		\qquad
		\forall o\in\mathcal{O},\;
		\forall d\in\mathcal{D},\;
		\forall k\in\mathcal{K}.
	\end{equation}
	
	\noindent Equations~\eqref{shift_continuous}-\eqref{projection_W_cont} define the \emph{continuous demand projection constraints} with respect to the problem $P_3$. Consequently, the continuous demand 
	projection mechanism modifies only the school-related demand,  which forms the demand input of optimization problem $P_3$.
	
	\subsection{Continuous Upper-Level Optimization Problem Formulation}
	\label{sec:Continuous_Upper_Optimization}
	
	The resulting Upper-Level optimization problem is then expressed as
	\begin{subequations}
		\label{rfwgfaaa}
		\begin{align}
			\label{arxontas}
			(P_3) \quad\underset{\psi_s, \forall s}{\text{Minimize}} &~ \displaystyle J_{\text{TTS}} \\ %+ \nonumber\\
			\textrm{Subject To:}&~~ \text{Demand Dynamics}~ \eqref{shift_continuous} - \eqref{projection_W_cont}, \nonumber \\
			&~~ \textrm{Traffic Dynamics}~ \eqref{mfd} - \eqref{fixed}, \nonumber\\
			& ~~\text{Constraints}~ \eqref{jam}, \eqref{minimal_disruption}, \nonumber\\
			&~~ \sum_{b\in\mathcal{B}}\sum_{s\in\mathcal{S}_b}|\psi_s| \leq \epsilon,\nonumber\\
			\label{jamt}
			&~~~\psi_s^l\leq \psi_s \leq \psi_s^u, \quad \forall b\in\mathcal{B}, \forall s\in\mathcal{S}_b, \\
			&~~~\psi_s\in\mathbb{R}, \psi_s^l\in\mathbb{R}^-, \psi_s^u\in\mathbb{R}^+,
			\forall b\in\mathcal{B}, \forall s\in\mathcal{S}_b,\nonumber\\
			\textrm{Initialization:}&~~\eqref{init},\nonumber\\
			\text{Input:} 
			\label{input_continuous_case}
			&~~ d_{o,s}^S(k;\boldsymbol{\psi}), 
			\quad \forall o\in\mathcal{O},~\forall s\in\mathcal{S}_b,~\forall b\in\mathcal{B},~\forall k\in\mathcal{K}, \nonumber\\
			&~~ d_{o,d}^W(k;\boldsymbol{\psi}), 
			\quad \forall o\in\mathcal{O},~\forall d\in\mathcal{D},~\forall k\in\mathcal{K}, \nonumber\\
			&~~\theta_{r,j,d}(k)\in[0,1],
			\quad \forall r\in\mathcal{R},~\forall j\in\mathcal{J}_r^-,~\forall k\in\mathcal{K},~\forall d\in\mathcal{D},\nonumber\\
			&~~\theta_{r,j,b}(k)\in[0,1],
			\quad \forall r\in\mathcal{R},~\forall j\in\mathcal{J}_r^-,~\forall k\in\mathcal{K},~\forall b\in\mathcal{B}.
		\end{align} 
	\end{subequations}
	
	\noindent Problem $P_3$ is a Nonlinear Program (NLP), where $\psi_s^l<0, \psi_s^u>0$ denote the smallest and largest admissible school start time shifts, expressed in time steps (backward and forward, respectively).\\
	
	\textbf{Remark (role of the continuous relaxation).} Problem $P_3$ constitutes a continuous relaxation of the original MINLP $P_1$. The optimization is performed entirely in the continuous decision space with respect to the decision variables $\psi_s$. The obtained solution therefore does not necessarily satisfy the discrete time-step structure imposed in $P_1$. To obtain a feasible solution for the original problem, the final converged solution of Algorithm~\ref{alg:seq1gsg_appendix} is projected onto the nearest multiple of the discrete interval $V$ (Lines 20-22).
	This projection is a post-processing feasibility recovery step and is not part of the optimization problem.\\
	
	Problem $P_3$ incorporates constraint \eqref{bint} expressed in absolute value form, which complicates its direct handling within the optimization process. To address this, we introduce the auxiliary non-negative variables $\psi_s^+$ and $\psi_s^-$ such that $\psi_s = \psi_s^+ - \psi_s^-, \forall s\in\mathcal{S}_b, \forall b\in\mathcal{B}$, which transform the absolute value constraint into an equivalent linear form that can be more readily integrated into the optimization program. To address this, we introduce the auxiliary nonnegative variables $\psi_s^+$ and $\psi_s^-$ such that $\psi_s = \psi_s^+ - \psi_s^-, \forall s\in\mathcal{S}_b, \forall b\in\mathcal{B}$, thereby transforming the absolute value constraints into equivalent linear forms that can be more readily integrated into the optimization program $P_3$. Hence, problem $P_3$ can be equivalently rewritten as
	\begin{subequations}
		\label{rfwgfaaa_new}
		\begin{align}
			\label{arxontas2}
			(P_4) \quad\underset{\psi_s^+, \psi_s^-, \forall s}{\text{Minimize}} &~ \displaystyle J_{\text{TTS}} \\ %+ \nonumber\\
			%& ~~\textcolor{blue}{\lambda \max(0,||J_{\text{STC}} - \eta||)} \\
			\textrm{Subject To:}&~~ \text{Demand Dynamics}~ \eqref{shift_continuous} - \eqref{projection_W_cont}, \nonumber \\
			&~~ \textrm{Traffic Dynamics}~ \eqref{mfd} - \eqref{fixed}, \nonumber\\
			& ~~\text{Constraints}~ \eqref{jam}, \eqref{minimal_disruption}, \nonumber\\
			& ~~\sum_{b \in \mathcal{B}} \sum_{s \in \mathcal{S}_b} \left(\psi_s^+ - \psi_s^- \right) \leq \epsilon, \label{eps_constraint_abs} \\
			\label{jamtT}
			&~~~\psi_s^l\leq \psi_s^+ - \psi_s^- \leq \psi_s^u, \quad \forall b\in\mathcal{B}, \forall s\in\mathcal{S}_b, \\
			&~~~\psi_s^l\in\mathbb{R}^-, \psi_s^u, \psi_s^+, \psi_s^-\in\mathbb{R}^+,
			\forall b\in\mathcal{B}, \forall s\in\mathcal{S}_b,\nonumber\\
			\textrm{Initialization:}&~~\eqref{init},\nonumber\\
			\textrm{Input:}&~~\eqref{input_continuous_case}.\nonumber
		\end{align} 
	\end{subequations}
	
	\noindent Problem $P_4$ constitutes a Nonlinear Program that is used for benchmarking purposes against the original MINLP $P_1$. As such, $P_4$ is intended for the construction of valid upper-bounds with respect to problem $P_1$.\\

\begin{algorithm}[t]
	\algsetup{linenosize=\tiny}
	\footnotesize
	%\scriptsize
	%\setlength{\baselineskip}{0.85\baselineskip}
	%{\fontsize{6}{6}\selectfont
		\begin{algorithmic}[1]
			\STATE \textbf{Input}: Value of parameter $\epsilon,$ School start time $\tau_s, \forall s\in\mathcal{S}_b, \forall b\in\mathcal{B}$, $\bar{t}_d, \bar{t}_d^{flex}, \forall d\in\mathcal{D}, \delta_2$.
			\STATE \quad \quad \quad Traffic Network Parameters, External Demand same as Algorithm \ref{alg:seq1gsg}.
			\STATE \textbf{Initialization:} $\lambda = 0$, $\Delta = +\infty$, $\mathbf{d}_{o,s}^{S}, \mathbf{d}_{o,d}^{W}$.
			\WHILE{$\Delta \ge \delta_2$}
			\IF{$\lambda=0$}
			\STATE Solve Optimization Problem $P_4$ when inserting $\epsilon$ as the right-hand side of constraint \eqref{eps_constraint_abs} and considering as input to the optimization problem $P_4$ the initial demand vectors $\mathbf{d}_{o,s}^{S}, \mathbf{d}_{o,d}^{W}$.
			\ELSE
			\STATE Solve Optimization Problem $P_4$ when inserting $\epsilon$ as the right-hand side of constraint \eqref{eps_constraint_abs} and giving as input to the optimization problem $P_4$ the equilibrium-based demand vectors $\mathbf{\hat{d}}_{o,s}^{S},\mathbf{\hat{d}}_{o,d}^{W}$ found from the previous outer-loop iteration $\lambda-1$.
			\ENDIF
			\STATE \textbf{Output}: Control decisions $\psi_s^+, \psi_s^-$ and calculate $\psi_s(\epsilon, \lambda) = \psi_s^+ - \psi_s^-, \; s \in \mathcal{S}_b, \; b \in \mathcal{B}$.
			\STATE Compute \textit{new} school start times: 
			$\tau_s + \psi_s \cdot T, \quad s \in \mathcal{S}_b, \; b \in \mathcal{B}$ and construct the school start time vector $\boldmath{\tau}^{\SurrogateAcronym}$$(\epsilon,\lambda)$.
			\STATE Execute Algorithm \ref{alg:user_equilibrium} using the vector of school start times $\boldmath{\tau}^{\SurrogateAcronym}$$(\epsilon,\lambda)$ and store
			the equilibrium-dependent demand vectors $ \mathbf{\hat{d}}_{o,s}^{S,I(\lambda)},\mathbf{\hat{d}}_{o,d}^{W,I(\lambda)}$ and the TTS, STC values, $J_{\text{TTS}}^{I(\lambda)}$, $J_{\text{STC}}^{I(\lambda)}$ using Eqs. \eqref{fixed} and \eqref{QUO1}, respectively, for the iteration $I(\lambda)$ that achieves convergence of Algorithm \ref{alg:user_equilibrium} and for the current outer-loop iteration $\lambda$.
			\IF{$\lambda > 0$}
			\STATE Compute $\Delta = |J_{\text{TTS}}^{I(\lambda)} - J_{\text{TTS}}^{I(\lambda-1)}|$.
			\ENDIF
			\STATE $\lambda = \lambda + 1$
			\ENDWHILE
			\STATE Set $\Lambda = \lambda-1$.
			\STATE Best continuous solution:
			Retrieve the corresponding optimal continuous shifts 
			$\psi_s^{*} = \psi_s(\epsilon,\Lambda)$.
			\STATE Projection onto discrete decision space:
			Project the continuous shifts onto the nearest multiple of the discrete interval $V$: $\psi_s^{\text{disc}} = \operatorname{round}_V(\psi_s^{*}), 
			\qquad \forall s\in\mathcal{S}_b,\; b\in\mathcal{B}$.
			\STATE Compute the feasible discrete school start times $
			\tau_s^{\text{disc}} = \tau_s + \psi_s^{\text{disc}} T,
			\qquad \forall s\in\mathcal{S}_b,\; b\in\mathcal{B}$.
			\STATE \textbf{Output:} Execute Algorithm~\ref{alg:user_equilibrium} using 
			$\boldsymbol{\tau}^{\text{disc}}$ and compute TTS value \text{$J_{\text{TTS}}^{\SurrogateAcronym}(\epsilon)$} and subsequently compute STC value $J_{\text{STC}}^{\SurrogateAcronym}(\epsilon)$ using Eq. \eqref{QUO1}. 
		\end{algorithmic}
		%	}
	\caption{\Surrogate\ (\SurrogateAcronym) Algorithm}
	\label{alg:seq1gsg_appendix}
\end{algorithm}
	
	Next, we develop the \Surrogate\ (\SurrogateAcronym) algorithm, whose pseudo-code is provided in Algorithm~\ref{alg:seq1gsg_appendix} to provide an upper-bound to the optimal value stemming from problem $P_1$. 
	
	\medskip
	\noindent\textbf{Remark.}  
	Algorithm~\ref{alg:seq1gsg_appendix} is structurally identical to the \Coupled\ (\Coupledacronym) algorithm, i.e., Algorithm \ref{alg:upper_bound}. Both algorithms follow the same sequential  interaction between the Upper-Level and the Lower-Level modules, alternating between:
	\begin{enumerate}
		\item Solving an Upper-Level optimization problem that generates candidate school start time shifts, and  
		\item Evaluating these decisions through the Approximate User Equilibrium Algorithm~\ref{alg:user_equilibrium}, which updates the demand vectors.
	\end{enumerate}
	
	The only difference lies in the optimization problem solved inside the outer loop.  
	While \Coupledacronym\ Algorithm solves the MILP \(P_2\) using the binary decision variables $\xi_{m,s}$, the present \SurrogateAcronym\ algorithm solves the nonlinear program \(P_4\) using decision variables \(\psi_s\). Consequently, the overall algorithmic structure, convergence logic, stopping criterion, and performance evaluation procedure remain unchanged.

	\section{Symbols, Notation and Parameters}
	\label{sec:appendixc}
	
	\counterwithin{table}{section}
	
	Table \ref{sec:table_example1} provides an overview of the symbols, notation and parameters used in this work to facilitate the reader’s understanding of the paper.
	
	\thispagestyle{empty}
	
	\begin{longtable}{@{}p{.15\textwidth}p{.70\textwidth}p{.15\textwidth}@{}}
		\caption{The symbols, notation, and parameters used in this work.}
		\label{sec:table_example1} \\
		\toprule
		Notation & Definition & Units \\
		\midrule
		\endfirsthead
		
		\multicolumn{3}{c}%
		{{\bfseries Table \thetable\ (Continued from previous page)}} \\
		\toprule
		Notation & Definition & Units \\
		\midrule
		\endhead
		
		\midrule
		\multicolumn{3}{r}{{Continued on next page}} \\
		\endfoot
		
		\bottomrule
		\endlastfoot
		
		$\mathcal{R}$ & Set of all regions in the network & - \\
		$\mathcal{O}$ & Set of origin regions in the network & - \\
		$\mathcal{B}\subseteq \mathcal{R}$ & Set containing the schools in the network & - \\
		$\mathcal{D}\subseteq \mathcal{R}$ & Set containing the destinations in the network & - \\
		$\mathcal{J}_r^-$ & Set of regions that can transmit flow to their neighbouring region $j$ & - \\
		$\mathcal{J}_r^+$ & Set of neighbouring regions directly accessible from region $r\in\mathcal{R}$ & - \\
		$\mathcal{S}_b$ & Set that contains the schools $s$ that are located in region $b\in\mathcal{B}$ & - \\
		$m$ & Number of discrete shifting intervals associated with the school start time & - \\
		$\mathcal{M}$ & Set that contains the possible shifting intervals of school start time & - \\
		$\mathcal{Y}$ & Set containing the indices corresponding to each class of commuters & - \\
		$k$ & Discrete time-step index & - \\
		$\mathcal{K}$ & Set that contains the discrete indices $k$ & - \\
		$\mathcal{T}$ & Set that contains school start time vectors & - \\
		$\zeta$ & Summation of cardinality of sets $\mathcal{S}_b, \forall b\in\mathcal{B}$ & - \\
		$\mathcal{E}$ & Set containing the values for the overall permissible change in the school start times & - \\
		$\epsilon$ & Overall permissible change in the school start times & min \\
		$T$ & Duration of each time-step & h \\
		$\Xi_U$ & Feasible set for the Upper-Level problem & - \\
		$\Xi_L$ & Feasible set for the Lower-Level problem & - \\
		$L_r$ & Total length of region $r\in\mathcal{R}$ & km \\
		$l_r$ & Average Trip length of vehicles inside region $r\in\mathcal{R}$ & km \\
		$\rho_r^C$ & Critical density of region $r\in\mathcal{R}$ & veh/km$\cdot$lane \\
		$\rho^J_{r}$ & Jam density of region $r\in\mathcal{R}$ & veh/km$\cdot$lane \\
		$C_{r,j}(\rho_j(k))$ & Inter-boundary capacity exchanged between region $r\in\mathcal{R}$ and $j\in\mathcal{J}_r^-$ at time-step $k$ & veh/h \\
		$C_{r,j}^{\text{MAX}}$ & Maximum inter-boundary capacity exchanged between region $r\in\mathcal{R}$ and $j\in\mathcal{J}_r^-$ & veh/h \\
		$\rho_{r}(k)$ & Instantaneous Density of region $r\in R$ at time-step $k$ & veh/km$\cdot$lane \\
		$\rho_{r}^y(k)$ & Instantaneous Density of region $r\in R$ with respect to commuters of class $y\in\mathcal{Y}$ at time-step $k$ & veh/km$\cdot$lane \\ 
		\midrule
		$q_{r}(k)$ & Intended outflow of region $r\in\mathcal{R}$ at time-step $k$ & veh/h \\
		$q_{r}^C$ & Capacity of region $r\in\mathcal{R}$ & veh/h \\
		$g_{r}(k)$ & Outflow generated from MFD of region $r\in\mathcal{R}$ at time-step $k$ & veh/h \\
		$g_{r}^C$ & Maximum Outflow generated from MFD of region $r\in\mathcal{R}$ & veh/h \\
		$u_{r}(k)$ & Average speed of region $r\in\mathcal{R}$ at time-step $k$ & km/h \\
		$u_{r}^f$ & Free-flow speed (slope of the first linear segment of piecewise MFD) of region $r\in\mathcal{R}$ & km/h \\
		$\rho_{o,r,s}^S(k)$ & Density of commuters of class S originating from $o\in\mathcal{O}$ and currently is in region $r\in\mathcal{R}$ that head to school $s\in\mathcal{S}_b, b\in\mathcal{B}$ at time-step $k$ & veh/km$\cdot$lane \\
		$\rho_{o,r,d}^W(k)$ & Density of commuters of class W originating from $o\in\mathcal{O}$ and currently is in region $r\in\mathcal{R}$ that head to destination $d\in\mathcal{D}$ at time-step $k$ & veh/km$\cdot$lane \\
		$q_{o,r,s}^S(k)$ & Intended transfer flow originating from $o\in\mathcal{O}$ and currently exits from region $r\in\mathcal{R}$ after having arrived at school $s\in\mathcal{S}_r, r\in\mathcal{B}$ at time-step $k$ & veh/h \\
		$q_{o,r,d}^W(k)$ & Intended transfer flow originating from $o\in\mathcal{O}$ and currently exits from region $r\in\mathcal{R}$ after having arrived at destination $d\in\mathcal{D}$ at time-step $k$ & veh/h \\
		\midrule
		$q_{r,j}(k)$ & Intended transfer flow of vehicles (in regional level) that transit from region $r\in\mathcal{R}$ at time step $k\in\mathcal{K}$ towards the directly accessible neighbouring region $j\in\mathcal{J}_r^+$ at time-step $k$ & veh/h \\
		$q_{o,r,j,s}^S(k)$ & Intended transfer flow originating from $o\in\mathcal{O}$ and currently exits from school $s\in\mathcal{S}_r$ located at region $r\in\mathcal{B}$ with respect to commuters of class S heading towards neighbouring region $j\in\mathcal{J}_r^-$ at time-step $k$ & veh/h \\
		$q_{o,r,j,d}^W(k)$ & Intended transfer flow originating from $o\in\mathcal{O}$ and currently exits from destination $d\in\mathcal{D}$ with respect to commuters of class W heading towards neighbouring region $j\in\mathcal{J}_r^-$ at time-step $k$ & veh/h \\
		$\tilde{q}_{o,r,j,s}^S(k)$ & Actual transfer flow originating from $o\in\mathcal{O}$ and currently transits from school $s\in\mathcal{S}_r$ located at region $r\in\mathcal{B}$ with respect to commuters of class S heading towards neighbouring region $j\in\mathcal{J}_r^-$ at time-step $k$ & veh/h \\
		$\tilde{q}_{o,r,j,d}^W(k)$ & Actual transfer flow originating from $o\in\mathcal{O}$ and currently transits from destination $d\in\mathcal{D}$ with respect to commuters of class W heading towards neighbouring region $j\in\mathcal{J}_r^-$ at time-step $k$ & veh/h \\
		$S_a(k)$ & Cumulative number of vehicles entering the network at time-step $k$ & veh \\
		$S_b(k)$ & Cumulative number of vehicles successfully exiting the network from destination $d$ at time-step $k$ & veh \\
		$J_{\text{TTS}}$ & Total Time Spent by all the vehicles in the network & veh h \\
		$J_{\text{STC}}$ & Overall change between the initial and the shifted start time of each school & min \\
		$d_{o,s}^{S}(k)$ & Instantaneous external demand for commuters of class S entering from $o\in O$ heading to school $s\in\mathcal{S}_b, b\in\mathcal{B}$ & veh \\
		$d_{o,s,m}^S(k)$ & Demand of commuters of Class S entering from $o\in O$ heading to school $s\in\mathcal{S}_b, b\in\mathcal{B}$ & veh \\
		$d_{o,d}^{W}(k)$ & Demand of commuters of Class W entering from region $o\in O$ heading to destination $d\in\mathcal{D}$ & veh \\
		$d_{b,d}^{S\rightarrow W}(k)$ & Demand of commuters that arrives at all schools located in region $b\in\mathcal{B}$ heading to their final destination $d\in\mathcal{D}$ & veh \\
		$\hat{d}_{o,s}^{S,i}(k)$ & Demand of commuters of class S entering from region $o\in\mathcal{O}$ heading to school $s\in\mathcal{S}_b, b\in\mathcal{B}$ during iteration $i$ of the \Approximate\ Algorithm & veh \\
		$\hat{d}_{o,d}^{W,i}(k)$ & Demand of commuters of class W entering from region $o\in\mathcal{O}$ heading to destination $d\in\mathcal{D}$ during iteration $i$ of the \Approximate\ Algorithm & veh \\
		$A_{o,s,d}$ & The ratio of vehicles originating from $o\in\mathcal{O}$ that travel from school $s\in\mathcal{S}_b, b\in\mathcal{B}$ to destination $d\in\mathcal{D}$ & - \\
		$\theta_{r,j,d}(k)$ & Split ratio of vehicles that transit from region $r\in\mathcal{R}$ at each time-step $k$ through the neighbouring region $j\in\mathcal{J}_r^+$ towards the corresponding school $s$ located in $b\in\mathcal{B}$  & - \\
		$V$ & Duration of the shifting interval for school start times & min \\
		$\varPhi$ & Range of candidate departure times associated with commuters of both classes & min \\
		$\tau_{s}$ & Initial start time of school $s\in\mathcal{S}_b, b\in\mathcal{B}$ & AM \\
		$\tilde{\tau}_{s}$ & Shifted start time of school $s\in\mathcal{S}_b, b\in\mathcal{B}$ & AM \\
		$\bar{t}_d$ & Fixed work start time of destination $d\in\mathcal{D}$ & AM \\
		$\bar{t}_d^{flex}$ & Flexible work start time of destination $d\in\mathcal{D}$ & AM \\
		$\xi_{m,s}$ & Binary variable indicating whether school $s$ starts at time $\tau_s+m V$ & - \\
		$\psi_s$ & Number of time-steps of duration $T$ (hours) that we shift the start time of school $s$ & - \\
		$t_{o,r}^{\text{SP}}$ & Duration of shortest path from origin $o$ to region $r$ & h \\
		$\Iota$ & Maximum number of iterations performed by \Approximate\ Algorithm & - \\
		$\Lambda$ & Maximum number of iterations performed by \Coupled\ Algorithm & - \\
		$\mathcal{P}_{o,d}$ & Set that contains the feasible regional paths for each origin-destination pair $(o,d)$ & - \\
		$\mathcal{Q}$ & Set that contains all the origin-destination pairs & - \\
		$TT_{p,o,d}(k)$ & Travel time of path $p\in\mathcal{P}_{o,d}$ when an OD pair departs at time-step $k$ & h \\[6pt]
		$\overline{TT}_{o,d}(k)$ & Expected travel time for pair $(o,d)$ at time-step $k$ & h \\
		$\varphi_{p,o,d}$ & Probability that a trip between origin $o$ and destination $d$ follows path $p$ & - \\
		$\alpha^S$ & Travel time cost for commuters of class S & (\euro/hour) \\
		$\beta^S$ & Early arrival cost for commuters of class S & (\euro/hour) \\
		$\gamma^S$ & Late arrival cost for commuters of class S & (\euro/hour) \\
		$\alpha^W$ & Travel time cost for commuters of class W with fixed work hours & (\euro/hour) \\
		$\beta^W$ & Early arrival cost for commuters of class W with fixed work hours & (\euro/hour) \\
		$\gamma^W$ & Late arrival cost for commuters of class W with fixed work hours & (\euro/hour) \\
		$\alpha^{W,flex}$ & Travel time cost for commuters of class W with flexible work hours & (\euro/hour) \\
		$\beta^{W,flex}$ & Early arrival cost for commuters of class W with flexible work hours & (\euro/hour) \\
		$\gamma^{W,flex}$ & Late arrival cost for commuters of class W with flexible work hours & (\euro/hour) \\
	\end{longtable}

\end{appendices}

%\textcolor{blue}{fix the references. Titles and Journals appear incorrectly with typos. Avoid Google Scholar bibtex entries.}

\bibliography{References}
\bibliographystyle{elsarticle-harv}

%\textcolor{blue}{issue with page appearing in the last page with the appendix}

\clearpage

%%%% This page is for instructions only, once the article is finalize please omit the below text before creating the final PDF

\end{document}